\documentclass{article}

\pdfoutput=1

\usepackage{htpy_in_fibs}

\usepackage{tocloft}
\title{Homotopies in Grothendieck fibrations}
\author{Joseph Helfer}
\date{}

\usepackage{fancyhdr}
\begin{document}
\maketitle

\begin{abstract}
  \noindent
  We define a natural 2-categorical structure on the base category of a large class of Grothendieck fibrations. Given any model category $\C$, we apply this construction to a fibration whose fibers are the homotopy categories of the slice categories $\C/A$, and we show that in the case $\C=\Top$, our construction applied to this fibration recovers the usual 2-category of spaces.
\end{abstract}

\vspace{5pt}
\tableofcontents
\vspace{5pt}

\section{Introduction}
The goal of this paper is to exhibit a naturally occurring 2-categorical structure on the base category of any Grothendieck fibration satisfying certain assumptions. In the motivating case of interest, the base category in question is the category of topological spaces, and our construction recovers the usual 2-category of topological spaces, continuous maps, and (homotopy classes of) homotopies.

The notion of Grothendieck fibration (as well as the essentially equivalent notion of \emph{pseudo-functor}) was first introduced \cite{grothdesc1,sga1} in order to formulate the notion of \emph{descent} (and later \cite{giraudcna}, \emph{stack}). Later, Lawvere introduced fibrations into categorical logic with his theory of \emph{hyperdoctrines} \cite{lawvereadj}. It is the latter, logical use of hyperdoctrines which is most relevant to the present work. Specifically, in an accompanying paper \cite{voltron}, we introduce ``homotopical'' semantics for first-order logic, and the 2-categorical structure introduced here is used to prove a ``homotopy-invariance'' theorem for these semantics.

We recall that a Grothendieck fibration (see \S\ref{subsec:fibs}) is a functor $\fibr{C}CB$ satisfying certain conditions which allows us to define, for each morphism $f\colon{}A\to{}B$ in $\B$, a functor $f^*\colon\fib{C}^B\to\fib{C}^A$ between the corresponding ``fibers'' of $\fib{C}$. We will be considering certain fibrations which (following \cite{folds}) we call $\wedgeq$-fibrations. Among the conditions for $\fib{C}$ to be a $\wedgeq$-fibration are that the category $\B$ have finite products and that, for each $B\in\Ob\B$, the fiber $\fib{C}^{B\times{}B}$ has an ``equality object'' $\Eq_B$, satisfying a certain universal property (the name ``equality object'' comes from the logical view of fibrations, in which the objects of the fiber $\fib{C}^B$ are viewed as predicates on the set $B$).

The 2-categorical structure in a $\wedgeq$-fibration arises as follows. Given two morphisms $f,g\colon{}A\to{}B$ in $\B$, we define a \emph{$\fib{C}$-homotopy} to be a morphism $\top_{A}\to\br{f,g}^*\Eq_B$ in $\fib{C}^A$, where $\br{f,g}$ is the induced morphism $A\to{}B\times{}B$ in $\B$, and $\top_A$ is the terminal object of $\fib{C}^A$. The 2-cells of the 2-categorical structure on $\B$ are given by the $\fib{C}$-homotopies. After this, the definitions of the remaining elements of the 2-categorical structure more or less suggest themselves.

The most natural source of Grothendieck fibrations are the ``codomain'' or ``family'' fibrations $\fibr{F(\C)}{C^\to}C$, in which $\C^\to$ is the category of morphisms in $\C$, and $\fib{F(\C)}$ is the functor sending each morphism to its codomain. This is a fibration whenever $\C$ has finite limits, and is in fact a $\wedgeq$-fibration. In these cases, the ``equality objects'' $\Eq_B$ are just the diagonal morphisms $\Delta\colon{}B\to{}B\times{}B$, and the resulting 2-categorical structure is trivial (i.e., the only 2-cells are identities).

The cases of interest are slight variations on the codomain fibrations. Here, we start with the category $\C$ of topological spaces (or Kan complexes, or more generally the category of fibrant objects of any Quillen model category), form the codomain fibration $\fib{F(\C)}$, and then take the homotopy category of each fiber of $\fib{F(\C)}$. The result, $\fib{HoF(\C)}$ is still a $\wedgeq$-fibration; but now, the equality objects $\Eq_B$ are the ``path-space fibrations'' $B^I\to{}B\times{}B$, and the $\fib{HoF(\C)}$-homotopies are the (homotopy classes of) homotopies, in the usual sense.

Let us say something about which aspects of our results are already known and which are (as far as we know) new. The notion of \emph{equality} in a fibration was introduced in \cite{lawvereequality}. The basic ingredients which go into the definition of our 2-category are well-known (see \cite[Lemma~3.4.5]{jacobscatlogic}) -- for example, the definition of the ``vertical'' composition simply amounts to the proof that this notion of equality is transitive. However, the fact that these ingredients can be used to define a 2-category has not, to our knowledge, been observed, though we should mention that \cite[p.~214]{jacobscatlogic} effectively constructs the ``homotopy category'' of our 2-category in the case of fibrations whose fibers are pre-orders (whence it follows that the $\Hom$-categories of the associated 2-category are also pre-orders.).

The idea that \emph{equality} is related to \emph{homotopy} is central to \emph{Homotopy Type Theory} (see \cite{awodeywarrenid,kapulkinlumsdaine,warrenthesis}), which was the direct inspiration for this work (and for \cite{voltron}). In particular, \cite{warrenthesis} observes that what we call $\fib{F_\fb(\C)}$ (see Definition~\ref{defn:mfib-subfibs}) is a $\wedge$-fibration, and that the path objects in it satisfy a ``weak'' version of the universal property of equality objects, though our ``fibration of homotopy categories'' $\fib{HoF(\C)}$ (see \S\ref{subsec:modelcats-hofib}), in which the path objects have the stronger universal property, is not considered (see the introduction to Part~\ref{sec:modelcats} for more on this).

At the time of first writing, the fibration $\fib{HoF(\C)}$ was, to our knowledge, new, but we since learned from Chaitanya Subramaniam that it has also been constructed by P. Cagne in \cite{cagnethesis}, where it is also shown to be equivalent to the fibration described in \S\ref{subsubsec:pseudo-functor-approach}.

The paper is organized as follows:

{\bf Part~\ref{sec:proof}:} We introduce the notion of homotopies in fibrations, and use it to construct a 2-categorical structure on the base of a $\wedgeq$-fibration. We prove some additional properties about this 2-categorical structure, namely its compatibility with the finite products on the base, and with the pseudo-functor associated to the fibration.

{\bf Part~\ref{sec:modelcats}:} We give examples of $\wedgeq$-fibrations. In particular, we associate to any model category $\C$ a fibration whose fibers are the homotopy categories of the slice categories of $\C$.

{\bf Part~\ref{sec:classical-2cat}:} We relate the 2-categorical structure on a model category arising from Parts~\ref{sec:proof}~and~\ref{sec:modelcats} to the standard 2-categorical structure on such a category.

{\bf Acknowledgments:}
We thank M. Makkai for reading and giving helpful comments on an early version of this paper, and Arpon Raksit for many helpful discussions.

We also thank the anonymous referee for many insightful recommendations and substantial improvements to the paper.

\subsection{Preliminaries on fibrations}\label{subsec:fibs}
We will now fix our notation and terminology regarding Grothendieck fibrations. This will be fairly cursory, and we will not give much in the way of motivation or proofs. For the latter things, we refer the reader to \cite{makkai-lauchli1,jacobscatlogic}.

\subsubsection{}\emph{Categories}.
We will use standard and hopefully familiar notation regarding categories. For basic notions of category theory, we refer to \cite{cwm}.

We note that we will always write $\id_A$, or just $\id$, for the identity morphism in a category, and we will avoid using the notation $\hid_A$, since that notation is reserved for ``identity $\fib{C}$-homotopies'' (see Definition~\ref{defn:id-c-htpy}).

When we say that some claim follows ``from a diagram chase'' in a given diagram, we mean that two particular morphisms, each represented by a composite of arrows in the diagram, are equal, and the proof is by a sequence of equalities, each coming from the commutativity of some sub-diagram. The commutativity of the sub-diagrams will always be evident -- often it follows from the definition of one of the morphisms involved -- and will usually be left to the reader.

We will be dealing with \emph{2-categories}. Given a 2-category $\C$ and objects $A,B\in\Ob\C$, we denote by $\HOM_\C(A,B)$ the corresponding Hom category. Each 2-category has an \emph{underlying} category, obtained by disregarding the 2-cells. Conversely, given a category $\C$, we can talk about a \emph{2-category structure} on $\C$, meaning a 2-category with underlying category $\C$ -- this will be our preoccupation in Part~\ref{sec:proof}.

\subsubsection{}\label{subsubsec:fib-basics}\emph{Fibrations}.
A \emph{prefibration} $\fibr{C}CB$ is just a functor $\fib{C}\colon\cat{C}\to\cat{B}$. $\C$ is the \emph{total category} of $\fib{C}$ and $\B$ is its \emph{base category} (and $\fib{C}$ is a prefibration \emph{over} its base category $\B$).
\vspace{-5pt}

For the rest of \S\ref{subsubsec:fib-basics}, fix a prefibration $\fibr{C}CB$.
\vspace{-5pt}

Given an object $B$ in $\B$, the \emph{fiber} of $\fib{C}$ over $B$, denoted $\fib{C}^{B}$, is the subcategory of $\C$ consisting of objects $P$ with $\fib{C}(P)=B$ and morphisms $p\colon{}P\to{}Q$ with $\fib{C}(p)=\id_B$. An object in the fiber over $B$ is said to \emph{lie over} $B$ or to be an \emph{object over} $B$. Similarly, a morphism $p$ in $\C$ with $\fib{C}(p)=f$ is said to \emph{lie over} $f$ or to be a \emph{morphism over} $f$, and is also said to be a \emph{lift} of $f$. We might also say $p$ lies over an object $B$ if $p$ lies over $\id_B$. Note that the property of, say, $p$ in $\C$ lying over $f$ in $\B$, depends on $\fib{C}$, so we should really say something like ``$p$ $\fib{C}$-lies-over $f$''. However, the prefibration $\fib{C}$ will always be clear from context in this and similar expressions.

When displaying diagrams in the total category $\C$ of a prefibration, we will usually display underneath it a diagram in the base category $\B$ so that each displayed object and morphism of $\C$ is positioned (approximately) above the object or morphism of $\B$ over which it lies. For example, in the diagram in \S\ref{subsubsec:wedge-fib}, $P$, $Q$, $P\wedge{}Q$, and the morphisms connecting them lie over $B$, and the horizontal morphisms lie over $f$.

Given a morphism $f\colon{}A\to{}B$ in $\B$ and objects $P,Q$ over $A,B$ in $\C$, we denote by $\Hom^{\fib{C}}_f(P,Q)$, or just $\Hom_f(P,Q)$, the set of morphisms $P\to{}Q$ lying over $f$. Note that for an object $A$ in $\B$ and objects $P$, $Q$ over $A$, $\Hom_{\id_A}(P,Q)$ is the same as $\Hom_{\fib{C}^A}(P,Q)$.

Recall that a morphism $q\colon{}Q\to{}R$ in $\C$ lying over $g\colon{}B\to{C}$ in $\B$ is \emph{cartesian} if for each morphism $f\colon{}A\to{}B$ in $\B$ and each $P\in\Ob\fib{C}^A$, the map $(g\ccmpf\mathord{\text{--}})\colon\Hom_{f}(P,Q)\to\Hom_{gf}(P,R)$ induced by composition with $q$ is a bijection. We will sometimes emphasize that a morphism in a diagram is cartesian by marking it like so $\cto$.

One (obvious) property of cartesian morphisms which we will frequently use is the following ``cancellation property'': if $r\colon{}Q\to{}R$ is a cartesian morphism in $\C$, and $p,q\colon{}P\to{}Q$ lie over the same morphism in $\B$, then $rp=rq$ implies $p=q$. Indeed, many of the propositions below state that two given morphisms are equal, and what will often be proven is that they are equal after composing with a certain cartesian morphism.

$\fib{C}$ is said to be a (Grothendieck) \emph{fibration} if for each morphism $f\colon{}A\to{}B$ in $\B$ and each object $P$ over $\B$, there is a cartesian lift of $f$ with codomain $P$.

\subsubsection{}\label{subsubsec:cleavage}\emph{Cleavages}. Let $\fibr{C}CB$ be a fibration. Recall that a \emph{cleavage} of $\fib{C}$ is a choice, for each morphism $f\colon{}A\to{}B$ in $\B$ and each $P\in\Ob\fib{C}^B$, of a cartesian lift of $f$ with codomain $P$. Assuming the axiom of choice, any fibration admits a cleavage.

We will use the following notational convention. Whenever we are dealing with a cloven fibration, we will, unless stated otherwise, denote the cartesian lifts in it by $\crt{f}Q\colon{}f^*Q\to{}Q$, or simply $\ct\colon{}f^*Q\to{}Q$.
Also, given a morphism $p\colon{}P\to{}Q$ over a composite $A\tox{f}B\tox{g}C$, we will denote by $\cind{p}$ the unique morphism $P\to{}g^*Q$ over $f$ such that $\ct\cind{p}=p$, i.e., such that
\[
  \begin{tikzcd}
    &g^*Q\car[dr, "\ct"]&\\[-10pt]
    P\ar[rr, "p"]\ar[ru, "\cind{p}", dashed]&&Q\\[-10pt]
    A\ar[r, "f"]&B\ar[r, "g"]&C
  \end{tikzcd}
\]
commutes (the notation hides the dependency on $f$, but this will always be clear from context).

In case we are considering a fibration and have \emph{not} chosen a cleavage, we may still use the above notation, but in this case it will be merely suggestive. For example, if we are considering a cartesian lift of $f$ with codomain $Q$, we may like to call it $\crt{f}Q$ or $\ct$, and call its domain $f^*Q$.

Recall that given a cleavage of $\fib{C}$, there is a associated to each morphism $f\colon{}A\to{}B$ in $\B$ a functor $f^*\colon\fib{C}^B\to\fib{C}^A$, the \emph{pullback functor along $f$}, taking $P\in\Ob\fib{C}^B$ to $f^*P$ and taking $p\colon{}P\to{}Q$ in $\fib{C}^B$ to $\cind{p\cdot\crt{f}P}$:
\[
  \begin{tikzcd}
    f^*P\ar[r, "\crt{f}P"]\ar[d, "f^*p=\cind{p\cdot\crt{f}P}"', dashed]&P\ar[d, "p"]\\
    f^*Q\car[r, "\crt{f}Q"]&Q.\\[-10pt]
    A\ar[r, "f"]&B
  \end{tikzcd}
\]

Note that whenever a morphism $p\colon{}P\to{}Q$ is cartesian, so is each induced morphism $\cind{p}\colon{}P\to{}f^*Q$. In particular, given a composite $h=gf\colon{}A\to{}C$ and an object $Q$ lying over $C$, the morphism $\ccind\colon{}h^*Q\to{}f^*Q$ is always cartesian.

Also note that the equation $\cind{q}p=\cind{qp}$ holds whenever it makes sense. We will often use these facts without explicit mention.

\subsubsection{}\label{subsubsec:fpnotations}\emph{Finite product categories}. By a \emph{finite product category} (or \emph{f.p.\ category}), we mean a category in which there exists a terminal object, and for each pair of objects $A$,$B$, there exists a product diagram \mbox{$A\ot{}P\to{}B$}.

A functor between f.p.\ categories is an \emph{f.p.\ functor} if it takes terminal objects to terminal objects and product diagrams to product diagrams.

When considering an f.p.\ category $\B$, we will often assume that $\B$ admits a choice of product diagram over each pair of objects, and will fix such a choice, as well as a choice of terminal object. Whenever we have fixed such choices, we will -- unless stated otherwise -- denote the chosen product of $A$ and $B$ by $A\times{}B$, the chosen product projections by $A\xot{\pi_1}A\times{}B\tox{\pi_2}B$, and the chosen terminal object by $\tm_{\C}$ (or just $\tm$).

However, when the category under consideration is a fiber $\fib{C}^A$ of some fibration $\fib{C}$, we will instead denote the chosen products by $P\wedge{}Q$ and the chosen terminal object by $\top_{A}$.

In either case, we denote by $\br{f,g}$ the morphism $X\to{}Y\times{}Z$ (or $X\to{}Y\wedge{}Z$) induced by $f\colon{}X\to{}Y$ and $g\colon{}X\to{}Z$, and by $!_{X}$ the unique morphism from $X$ to the terminal object. We write $f\times{}g$ for $\br{f\pi_1,g\pi_2}$, and $\Delta_X$ for $\br{\id_X,\id_X}$.

We assume the $\times$ and $\wedge$ are left-associative so, e.g., $A\times{B}\times{C}=(A\times{B})\times{C}$. We write $\br{f,g,h}$ for $\br{\br{f,g},h}$, $\Delta_B^3$ for $\br{\br{\id_B,\id_B},\id_B}$, as well as $\pi_1,\pi_2,\pi_3\colon{}A\times{B}\times{C}\to{A}$ for $\pi_1\pi_1,\pi_2\pi_1,\pi_2$, and similarly with $\br{f,g,h,k}$ and so on.

In case we have \emph{not} chosen distinguished binary products and terminal object, we may still use the above notation, but in this case it will be merely suggestive (this is similar to the convention in \S\ref{subsubsec:cleavage}). For example, if we wish to consider a product diagram over objects $A$ and $B$, we may like to call its vertex $A\times{}B$ and its projections $\pi_1$ and $\pi_2$.

\subsubsection{}\label{subsubsec:wedge-fib}\emph{$\wedge$-fibrations}.
A fibration $\fibr{C}CB$ is a \emph{$\wedge$-fibration} if it satisfies the following three conditions:
\vspace{-10pt}
\begin{enumerate}[(i)]
\item Each fiber $\fib{C}^A$ is an f.p.\ category.

\item\label{item:wedgefibtm}
Given a cartesian morphism $\crt{f}\top_A\colon{}f^*\top_B\to{}\top_B$ in $\C$ over a morphism $f\colon{}A\to{}B$ in $\B$, where $\top_B$ is terminal in $\fib{C}^B$, $f^*\top_B$ is terminal in $\fib{C}^A$.

\item\label{item:wedgefibprod}
Given a commutative diagram
\[
  \begin{tikzcd}
    &[-40pt]f^*(P\wedge{}Q)\ar[ldd]\ar[rd]\car[rrr]&[-20pt]&
    &[-30pt]P\wedge{}Q\ar[rd]&[-20pt]\\
    &&f^*Q\ar[rrr, "\carsym" {anchor=center, pos=0.4}]&&&Q\\
    f^*P\car[rrr]&&&P\ar[from=uur, crossing over]&&\\
    &A\ar[rrr, "f"]&&&B
  \end{tikzcd}
\]
in which $P\ot{}P\wedge{}Q\to{}Q$ is a product diagram in $\fib{C}^B$ and the horizontal arrows are cartesian over $f$, we have that $f^*P\ot{}f^*(P\wedge{}Q)\to{}f^*Q$ is a product diagram in $\fib{C}^A$.
\end{enumerate}

We will sometimes refer to the last two properties as ``stability'' of products and terminal objects (under pullbacks). Given a cleavage of $\fib{C}$, the conjunction of (ii) and (iii) is equivalent to the condition that each pullback functor $f^*$ is an f.p.\ functor.

A \emph{$\wedge$-cleavage} of $\fib{C}$ is a cleavage together with a choice of binary products and terminal object in each fiber.

Whenever we have fixed a $\wedge$-cleavage of $\fib{C}$, we will use the following notation.

Given a morphism $f\colon{}A\to{}B$ in $\B$ and a product $Q\xot{\pi_1}Q\wedge{}R\tox{\pi_2}R$ in the fiber $\fib{C}^B$, and given morphisms $q\colon{}P\to{}Q$ and $r\colon{}P\to{}R$ over $f$, we denote by $\brr{q,r}$ the unique morphism over $f$ satisfying $\pi_1\brr{q,r}=q$ and $\pi_2\brr{q,r}=r$ (we leave it to the reader to see that there is a unique such morphism).

Note that if $A=B$ and $f=\id_A$, then $\brr{q,r}=\br{q,r}$. If $P$ is itself a product, we also write $q\wwdge{}r$ for $\brr{q\pi_1,r\pi_2}$. Note that if $q$ and $r$ lie over an identity morphism then $q\wwdge{}r=q\wedge{}r$.

Note that we have the following familiar equations (whenever they make sense):
\[
  \brr{p,q}r=\brr{pr,qr}
  \quad
  \quad
  (s\wwdge{}t)\brr{p,q}=\brr{sp,tq}
  \quad
  \quad
  (s\wwdge{}t)(p\wwdge{}q)=sp\wwdge{}tq
\]

For a morphism $g\colon{}B\to{}C$ in $\B$ and an object $Q\in\Ob\C^B$ we denote by $\exx_g$ (or just $\exx$) the unique morphism $Q\to\top_C$ over $g$. Note that if $B=C$ and $g=\id_B$, then $\exx_g=!_Q$ and that, given a morphism $p\colon{}P\to{}Q$ over $f\colon{}A\to{}B$, we have $\exx_gp=\exx_{gf}\colon{}P\to\top_C$.

As usual (see \S\S\ref{subsubsec:cleavage}~and~\ref{subsubsec:fpnotations}), we may still use these notations even if we have not chosen a $\wedge$-cleavage, but in this case they will be merely suggestive.

\subsubsection{}\label{subsubsec:cocart}\emph{Cocartesian morphisms}. We recall the notion of a cocartesian morphism, dual to that of cartesian, i.e., a morphism in the total category of a prefibration $\fibr{C}CB$ is \emph{cocartesian} if it is cartesian in the prefibration $\fibr{C^\op}{\C^\op}{\B^\op}$.

For each statement involving cartesian morphisms, there is of course a dual statement for cocartesian morphisms.

If $\fib{C}$ is a fibration, then a morphism $p\colon{}P\to{}Q$ in $\C$ over $f\colon{}A\to{}B$ in $\B$ is cocartesian if and only if it is \emph{weakly cocartesian}, meaning that $(\mathord{\text{--}}\ccmpf{}p)\colon\Hom_{\fib{C}^B}(Q,R)\to\Hom_f(P,R)$ is a bijection for every $R\in\Ob\fib{C}^B$.

For the rest of \S\ref{subsec:fibs}, let $\fibr{C}CB$ be a fibration.

Suppose we have a cocartesian morphism $p\colon{}P\to{}Q$ in $\C$ over $g\colon{}C\to{}D$ in $\B$, and a morphism $k\colon{}B\to{}D$. We say that $p$ is \emph{stable along $k$} if the following condition holds. \[
  \begin{tikzcd}[row sep=8pt]
    &P'\ar[rr, "p'"]\car[ld, "s"']&&Q'\car[ld, "t"]\\
    P\rac[rr, "p"]&&Q&\\[-5pt]
    &A\ar[rr, "f"]\ar[ld, "h"']&&B\ar[ld, "k"]\\
    C\ar[rr, "g"]&&D&
  \end{tikzcd}
\]
Given any commutative square in $\C$ lying over a square in $\B$, as shown above, if the square in $B$ is a pullback square, and $s$ and $t$ are cartesian, then $p'$ is cocartesian (note that in this situation, $p'$ is uniquely determined by $p$, $s$, and $t$ since $t$ is cartesian). This is also known as the \emph{Beck-Chevalley condition}.

A morphism $\pi_2\colon{}X\times{}Y\to{}Y$ in any category $\D$ is a \emph{product projection} if there is a morphism $\pi_1\colon{}X\times{}Y\to{}X$ in $\D$ for which $X\xot{\pi_1}X\times{}Y\tox{\pi_2}Y$ is a product diagram. A morphism $\Delta_X\colon{}X\to{}X\times{}X$ is a \emph{diagonal morphism} if there exists a product diagram $X\xot{\pi_1}X\times{}X\tox{\pi_2}X$ with $\pi_1\Delta_X=\pi_2\Delta_X=\id_X$. A morphism $\id_X\times\Delta_Y\colon{}X\times{}Y\to{}X\times{}Y\times{}Y$ is a \emph{generalized diagonal morphism} if there exists a pullback diagram
\[
  \begin{tikzcd}
    X\times{}Y\ar[d, "\pi_2'"']\ar[r, "\id_X\times\Delta_Y"]&[30pt]
    X\times{}Y\times{}Y\ar[d, "\pi_2"]\\
    Y\ar[r, "\Delta_Y"]&Y\times{}Y
  \end{tikzcd}
\]
with $\Delta_Y$ a diagonal morphism and $\pi_2$ a product projection. Note that it follows that $\pi_2'$ is a product projection as well. Note also that in an f.p.\ category, given the diagonal morphism $\Delta_Y$ and the product projection $\pi_2$, such a pullback diagram always exists.

Suppose $\fib{C}$ is a $\wedge$-fibration, and let $q\colon{}Q\to{}Q'$ in $\C$ over $f\colon{}A\to{}B$ in $\B$ be cocartesian. We say that $q$ satisfies \emph{Frobenius reciprocity} if the following condition holds. Given any commutative diagram
\[
  \begin{tikzcd}[column sep=40pt]
    Q\rac[r, "q"]&Q'\\
    Q\wedge{}f^*P\ar[r, "q\wwdge{}\crt{f}p"]\ar[u]\ar[d]&Q'\wedge{}P\ar[u]\ar[d]\\
    f^*P\car[r, "\crt{f}p"]&P\\[-10pt]
    A\ar[r, "f"]&B
  \end{tikzcd}
\]
with $Q\ot{}Q\wedge{}f^*P\to{}f^*P$ and $Q'\ot{}Q'\wedge{}P\to{}P$ product diagrams, if $\crt{f}p$ is cartesian over $f$, then $q\wwdge{}\crt{f}p$ is cocartesian.

\subsubsection{}\emph{$\wedgeq$-fibrations}.
$\fib{C}$ is a \emph{$\wedgeq$-fibration} if it satisfies the following three conditions.
\vspace{-10pt}
\begin{itemize}
\item[(i)] $\fib{C}$ is a $\wedge$-fibration.
\item[(ii)] For every generalized diagonal morphism $\id_A\times\Delta_B\colon{}A\times{}B\to{}A\times{}B\times{}B$ in $\B$ and every terminal object $\top_{A\times{}B}$ of $\fib{C}^{A\times{}B}$, there is a cocartesian lift of $\id_A\times\Delta_B$ with domain $\top_{A\times{}B}$.
\item[(iii)] Every cocartesian lift as in (ii) satisfies Frobenius reciprocity and is stable with respect to all product projections $\pi_2\colon{}C\times{}A\times{}B\times{}B\to{}A\times{}B\times{}B$.
\end{itemize}

A \emph{$\wedgeq$-cleavage} of $\fib{C}$ is $\wedge$-cleavage together with a choice of binary products and terminal object in $\B$, and a choice, for each $B\in\Ob\B$, of a cocartesian lift of $\Delta_B\colon{}B\to{}B\times{}B$ with domain $\top_B$. For $B\in\Ob\B$, we will denote the codomain of this chosen lift by $\Eq_B$, and the lift itself by $\rho_B\colon\top_B\to\Eq_B$. ($\rho$ stands for ``reflexivity''.)

Note that by the stability demanded in (iii) and by the definition of generalized diagonal morphism, it is enough in (ii) to demand cocartesian lifts of (non-generalized) diagonal morphisms.

Also note that all the conditions in the definition of a $\wedgeq$-fibration are ``isomorphism invariant''. Hence, in (ii), for example, it suffices to check the condition for \emph{some} terminal object $\top_{A\times{}B}$ in $\fib{C}^{A\times{}B}$, and for \emph{some} pullback $\id_A\times{}\Delta_B$ of \emph{some} diagonal morphism $\Delta_B\colon{}B\to{}B\times{}B$ along \emph{some} product projection \mbox{$\pi_2\colon{}A\times{}B\times{}B\to{}B\times{}B$}. In particular, if we have chosen a $\wedgeq$-cleavage of $\fib{C}$, it suffices to check that above conditions for the specified structure in the $\wedgeq$-cleavage. This ``isomorphism invariance'' is not exactly trivial, but is straightforward to (formulate precisely and) check, and we leave it to the reader.

As a final remark on the definition, we note that each instance of Frobenius reciprocity produces a new cocartesian morphism; hence, a $\wedgeq$-fibration necessarily has many cocartesian morphisms not explicitly required in the definition. One might wonder if these other cocartesian morphisms automatically satisfy Frobenius reciprocity and stability along product projections as well. They do (this, too, we leave to the reader).

\section{2-categorical structure in fibrations}\label{sec:proof}
The goal of Part~\ref{sec:proof} is to define a 2-category structure on the base category of any (cleavable) $\wedgeq$-fibration $\fibr{C}CB$. The 2-cells between two morphisms $f,g\colon{}A\to{}B$ will be ``homotopies'' between $f$ and $g$. From the logical point of view, these are ``proofs according to $\fib{C}$'' that two morphisms $f,g\colon{}A\to{}B$ are equal.

According to this point of view, we think of the objects of $\B$ as denoting sets, and the objects in the fiber over the object $A$ as denoting predicates on $A$, the morphisms between them being implications. The terminal object $\top_A$ is then the trivial predicate ``true'', and the equality object $\Eq_B$ over $B\times{}B$ is the equality predicate ``$b=b'$''. Moreover, given a morphism $f\colon{}A\to{}B$, we think of the pullback $f^*$ as performing ``substitution'', i.e.\ $f^*$ takes the predicate $P(b)$ to $P(f(a))$. Hence, a ``proof that $f$ and $g$ are equal'' -- i.e.\ that the predicate ``$f(a)=g(a)$'' is always true -- should be an implication $\top_A\to\br{f,g}^*\Eq_B$, or equivalently, a morphisms $\top_A\to\Eq_B$ over $\br{f,g}$.

The idea that this (or any) notion of ``equality'' should have something to do with \emph{homotopies} -- which it does (see Part~\ref{sec:classical-2cat}) -- is familiar from homotopy theory; it often happens that each point of a space represents some (say, geometric) object, and that each path between two points gives rise to an identification between the corresponding objects. Thus one considers \emph{spaces} rather than \emph{sets} of objects, and ``path-connectedness'' rather than \emph{equality}. This idea is expressed most explicitly in Homotopy Type Theory which, as we mentioned in the introduction, is the inspiration for this work. In \cite{voltron}, we make the logical connection to the present work more explicit.

Having defined homotopies in this manner, it remains to define the 2-categorical structure; i.e., the composition operations. Here, the logical point of view is again helpful. For example, defining the composition of homotopies $f\to{}g$ and $g\to{}h$ amounts to producing a proof of $f(a)=h(a)$ from $f(a)=g(a)$ and $g(a)=h(a)$. However, this kind of thinking only goes so far; for example, proving associativity of this composition amounts to showing that ``two different proofs are equal'', which does not have a formal counterpart in the rules of predicate logic.\footnote{However, it does have a counterpart in Martin-L\"of type theory, and in fact, the constructions in this section are closely related to establishing a groupoid structure on identity types, as carried out, e.g., in \cite[p.7,~Proposition~4.1]{hofmannstreicher}. The precise relationship is somewhat subtle, since the ``type theory'' (i.e., notion of fibration) we use here is very minimal (for example, it has no function types), and it would be interesting to work it out in detail.}

Next, we will present an elegant alternative description, suggested by the anonymous referee, of the 2-categorical structure on $\B$, in which this structure is induced by embedding $\B$ into a larger 2-category -- namely the 2-category $\Cat(\C)$ of internal categories in $\C$.

Finally, we will prove two more properties of the 2-category $\B$ which are crucial to our application in \cite{voltron}. The first is that the finite products in $\B$ are also finite products in the 2-categorical sense.

The second is that the pseudo-functor $\B^\op\to\Cat$ associated to the cleavage of $\fib{C}$ extends to a pseudo-functor of 2-categories, or equivalently -- and it is actually this second statement that we prove\footnote{In an earlier version of this paper, we worked directly with pseudo-functors, but we changed to the cleaner approach using 1-discrete 2-fibrations upon a suggestion from the anonymous referee.} -- that the fibration $\fib{C}$ can be extended to a \emph{1-discrete 2-fibration} (see Definition~\ref{defn:1d2f}). The significance of this is roughly as follows. Given two parallel morphisms $f,g\colon{}A\to{}B$, we have the two pullback functors $f^*,g^*\colon\fib{C}^B\to\fib{C}^A$. If $f$ and $g$ are homotopic, we would expect these functors to be naturally isomorphic. The pseudo-functoriality says that this is so, and moreover that this association of natural isomorphisms to homotopies takes composition of homotopies to composition of natural transformations.

The extension of $\fib{C}$ to a 1-discrete 2-fibration with 2-categorical products in the base also allows us to prove a universal property of the 2-categorical structure, suggested by the referee, which characterizes it up to isomorphism.

Part~\ref{sec:proof} is organized as follows. In \S\ref{subsec:c-htpy-vert-comp}, we define the notion of $\fib{C}$-homotopy, and then define the ``vertical'' composition and show that it defines a category. In \S\ref{subsec:c-htpy-gpd}, we show that this category is in fact a groupoid -- i.e., that all the morphisms are invertible. In \S\ref{subsec:c-htpy-2cat}, we define the ``horizontal'' composition, and show that, together with the ``vertical'' composition, this forms a 2-category. In \S\ref{subsec:intcats}, we give the alternative presentation of the 2-categorical structure in terms of internal categories. In \S\ref{subsec:1d2fs}, we carry out the extension of the fibration to a 1-discrete 2-fibration, and in \S\ref{subsec:c-htpy-finprod}, we show that the 2-category has finite products. Finally, in \S\ref{subsec:c-htpy-canon}, we present the universal property characterizing the 2-categorical structure.

\subsection{Homotopies in fibrations}\label{subsec:c-htpy-vert-comp}
We now introduce the notion of homotopies in fibrations, and define the ``vertical'' composition operation. Throughout \S\ref{subsec:c-htpy-vert-comp}, let $\fibr{C}CB$ be a $\wedgeq$-cloven $\wedgeq$-fibration.

The definition of homotopy was sketched above. As we also mentioned there, given three morphisms $f,g,h\colon{}A\to{}B$, the definition of the composition of homotopies $f\to{}g$ and $g\to{}h$ can be viewed logically as a proof of transitivity $f(a)=g(a)\wedge{}g(a)=h(a)\To{}f(a)=h(a)$. This, in turn can be reduced to the general statement $b_1=b_2\wedge{}b_2=b_3\To{}b_1=b_3$ for $b_1,b_2,b_3\in{}B$. This is how we will proceed, noting that, in terms of the fibration, the predicates $b_i=b_j$ are represented by the pullback of $\Eq_B$ along $\br{\pi_i,\pi_j}\colon{}B^3\to{}B\times{}B$.

\subsubsection{}
\defn
Given two morphisms $f,g\colon{}A\to{B}$ in $\B$, a \emph{$\fib{C}$-homotopy} from $f$ to $g$ is a morphism $\top_A\to\Eq_B$ over $\br{f,g}\colon{}A\to{}B\times{}B$.\vspace{-15pt}

To denote a $\fib{C}$-homotopy $\alpha$ from $f$ to $g$, we use the notation $f\tox{\alpha}g$ or
$
\begin{tikzcd}[row sep=0pt, column sep=10pt]
  &\ar[dd, shorten >=2pt, shorten <=-2pt, Rightarrow, "\ \alpha", pos=0.4]&\\[-5pt]
  A\ar[rr, "f", bend left]\ar[rr, "g"', bend right]&&B.\\
  &{}&
\end{tikzcd}
$\vspace{-10pt}

\subsubsection{}\defn
Given $B\in\Ob\C$ and natural numbers $1\le{}i,j\le{}n$, we write $\Eq_B^{ij}\in\Ob\fib{C}^{B^n}$ for the pullback $\br{\pi_i,\pi_j}^*\Eq_B$ of $\Eq_B\in\Ob\fib{C}^{B\times{}B}$ along $\br{\pi_i,\pi_j}\colon{}B^n\to{}B\times{}B$. The notation hides the dependency on $n$, but it will always be clear from context.

\subsubsection{}\defn
For $B\in\Ob\B$, we define $\rho_B^{ij}$ to be the unique morphism $\top_B\to\Eq_B^{ij}$ over $\Delta_B^n\colon{}B\to{}B^n$ making the diagram
\[
  \begin{tikzcd}
    &\Eq_B^{ij}\car[dr, "\ct"]&\\[-5pt]
    \top_B\ar[rr, "\rho_B"]\ar[ru, "\rho_B^{ij}", dashed]&&\Eq_B\\[-10pt]
    B\ar[r, "\Delta_B^n"]&B^n\ar[r, "\br{\pi_i,\pi_j}"]&B\times{}B
  \end{tikzcd}
\]
commute, i.e., $\rho_B^{ij}=\cind{\rho_B}$. Again, the dependence of $\rho_B^{ij}$ on $n$ is concealed in the notation, but it will always be clear from context.

\subsubsection{}\label{prop:tripeq}\prop
For every $B\in\Ob\B$, the morphism $\brr{\rho_B^{12},\rho_B^{23}}\colon\top_B\to\Eq_B^{12}\wedge\Eq_B^{23}$ over $\Delta_B^3$ is cocartesian.
\pf
We will show that $\brr{\rho_B^{12},\rho_B^{23}}$ is the composite of two cocartesian morphisms and hence cocartesian. The first will just be $\rho_B\colon\top_B\to\Eq_B$ over $\Delta_B\colon{}B\to{}B\times{}B$. The second will be a morphism $\Eq_B\to\Eq_B^{12}\wedge\Eq_B^{23}$ over $\id_B\times\Delta_B\colon{}B\times{}B\to{}B\times{}B\times{}B$ (where, abusing notation, we write $\id_B\times\Delta_B$ for $\br{\pi_1,\pi_2,\pi_2}$), which we will now construct.

Let us denote by $\tilde\rho_B^{23}$ the unique morphism $\top_{B\times{}B}\to\Eq_B^{23}$ over $\id_B\times\Delta_B$ making the diagram
\[
  \begin{tikzcd}
    \top_{B\times{}B}\ar[r, "\tilde\rho_B^{23}", dashed]
    \car[d, "\exx_{\pi_2}"']&\Eq_B^{23}\car[d, "\ct"]\\
    \top_B\rac[r, "\rho_B" pos=0.6]&\Eq_B
  \end{tikzcd}
  \quad\text{lying over}\quad
  \begin{tikzcd}
    B\times{}B\ar[r, "\id_B\times\Delta_B"]\ar[d, "\pi_2"']
    \ar[rd, phantom, "\lrcorner", pos=0]&[20pt]
    B\times{}B\times{}B\ar[d, "\br{\pi_2,\pi_3}"]\\
    B\ar[r, "\Delta_B"]&B\times{}B
  \end{tikzcd}
\]
commute. By the stability of $\rho_B$ along the product projection $\br{\pi_2,\pi_3}$, $\tilde\rho_B^{23}$ is cocartesian.

Now, consider the following commutative diagram.
\[
  \begin{tikzcd}[column sep=50pt]
    \top_{B\times{}B}\rac[r, "\tilde\rho_B^{23}"]&\Eq_B^{23}\\
    \Eq_B\ar[r, "\brr{\cind{\id_{\Eq_B}}, \tilde\rho_B^{23}!}"]
    \ar[u, "!"]&\Eq_B^{12}\wedge\Eq_B^{23}\ar[u, "\pi_2"']\\
    \Eq_B\car[r, "\cind{\id_{\Eq_B}}"]\ar[from=u, "\id"']&[20pt]\Eq_B^{12}\ar[from=u, "\pi_1"]\\[-5pt]
    B\times{}B\ar[r, "\id_B\times\Delta_B"]&B\times{}B\times{}B
  \end{tikzcd}
\]
Since the left and right sides are product diagrams, $\cind{\id_{\Eq_B}}$ is cartesian, and $\tilde\rho_B^{23}$ is cocartesian, it follows by Frobenius reciprocity that $\brr{\cind{\id_{\Eq_B}},\tilde\rho_B^{23}!}$ is cocartesian.

It remains to show that $\brr{\rho_B^{12},\rho_B^{23}}$ is equal to the composite
$\brr{\cind{\id_{\Eq_B}}\rho_B,\tilde\rho_B^{23}!\rho_B}$ of $\rho_B$ and $\brr{\cind{\id_{\Eq_B}},\tilde\rho_B^{23}!}$. That $\cind{\id_{\Eq_B}}\rho_B$ is equal to $\rho_B^{12}$ follows from the fact that $\rho_B^{12}=\cind{\rho_B}$, and that $\tilde\rho_B^{23}!\rho_B$ is equal to $\rho_B^{23}$ follows from a diagram chase in the following diagram.
\[
  \begin{tikzcd}[baseline=(QED.base)]
    \top_B\ar[rrr, bend left=25pt, "\rho_B^{23}"]\ar[r, "\rho_B"]\ar[rrd, "\id"]&
    \Eq_B\ar[r, "!"]&
    \top_{B\times{}B}\ar[r, "\tilde\rho_B^{23}"']\ar[d, "\exx_{\pi_1}"]&
    \Eq_B^{23}\car[d, "\ct"]\\[10pt]
    &&\top_{B}\ar[r, "\rho_B"]& |[alias=QED]| \Eq_B
  \end{tikzcd}
  \tag*{\qed}
\]

\subsubsection{}\label{defn:trB}
\defn
For an object $B$ in $\B$, we define $\tr_B\colon\Eq_B^{12}\wedge\Eq_B^{23}\to\Eq_B^{13}$ to be the (by Proposition~\ref{prop:tripeq} unique) morphism in $\fib{C}^{B^3}$ making the following diagram commute. (``tr'' stands for ``transitivity'')
\[
  \begin{tikzcd}
    &\Eq_B^{12}\wedge\Eq_B^{23}\ar[d, "\tr_B", dashed]\\
    \top_B\ar[r, "\rho_B^{13}"']\rac[ru, "\brr{\rho_B^{12},\rho_B^{23}}"]&\Eq_B^{13}\\[-5pt]
    B\ar[r, "\Delta_B^3"]&B^3
  \end{tikzcd}
\]

\subsubsection{}\label{defn:comp-c-htpy}
\defn
Given morphisms $f,g,h\colon{}A\to{B}$ in $\B$ and $\fib{C}$-homotopies
$f\tox{\alpha}g\tox{\beta}h$
(i.e., morphisms $\alpha,\beta\colon\top_A\to\Eq_B$ over $\br{f,g}$ and $\br{g,h}$ respectively), we define the \emph{vertical composite} $\beta\hvcmp\alpha$ to be the $\fib{C}$-homotopy from $f$ to $h$ given by the morphism
$\ct\ccmp\tr_B\ccmp\brr{\cind\alpha,\cind\beta}\colon\top_A\to\Eq_B$
over $\br{f,h}$, as shown below.
\[
  \begin{tikzcd}
    &\Eq_B^{12}\wedge\Eq_B^{23}\ar[d, "\tr_B"]&[20pt]\\
    \top_A\ar[ru, "\brr{\cind{\alpha},\cind{\beta}}"]
    &\Eq_B^{13}\ar[r, "\ct"]&\Eq_B\\[-10pt]
    A\ar[r, "\br{f,g,h}"]&B^3\ar[r, "\br{\pi_1,\pi_3}"]&B\times{}B
  \end{tikzcd}
\]
In other words, $\beta\hvcmp\alpha$ is the unique morphism over $\br{f,h}$ such that the following diagram commutes.
\[
  \begin{tikzcd}
    &\Eq_B^{12}\wedge\Eq_B^{23}\ar[d, "\tr_B"]\\
    \top_A\ar[ur, "\brr{\cind{\alpha},\cind{\beta}}"]\ar[r, "\cind{\beta\hvcmp\alpha}"]
    &\Eq_B^{13}\\[-10pt]
    A\ar[r, "\br{f,g,h}"]&B^3
  \end{tikzcd}
\]

\subsubsection{}\label{defn:id-c-htpy}
\defn
Given a morphism $f\colon{}A\to{B}$ in $\B$, we define the \emph{identity $\fib{C}$-homotopy at $f$}, denoted $\hid_f$, to be the morphism $\rho_B\ccmp\exx_f\colon\top_A\to\Eq_B$ over $\br{f,f}$, as shown below.
\[
  \begin{tikzcd}
    \top_A\ar[r, "\exx_f"]&\top_B\ar[r, "\rho_B"]&\Eq_B\\[-5pt]
    A\ar[r, "f"]&B\ar[r, "\Delta_B"]&B\times{}B
  \end{tikzcd}
\]

\subsubsection{}
\defn
Given an object $B\in\Ob\B$ and natural numbers $1\le{}i,j,k\le{}n$, we define \mbox{$\tr_B^{ijk}\colon\Eq_B^{ij}\wedge\Eq_B^{jk}\to\Eq_B^{ik}$} to be the unique morphism in $\fib{C}^{B^n}$ making the following diagram commute.
\[
  \begin{tikzcd}
    \Eq_B^{ij}\wedge\Eq_B^{jk}\ar[d, "\tr_B^{ijk}"', dashed]\ar[r, "\ccind\wwdge\ccind"]&
    \Eq_B^{12}\wedge\Eq_B^{23}
    \ar[d, "\tr_B"]
    \\
    \Eq_B^{ik}\car[r, "\ccind"]&
    \Eq_B^{13}\\[-10pt]
    B^n\ar[r, "\br{\pi_i,\pi_j,\pi_k}"]&B^n
  \end{tikzcd}
\]

\subsubsection{}\label{lem:trgen}\lem
Given morphisms $f_1,\ldots,f_n\colon{}A\to{B}$ in $\B$ and $\fib{C}$-homotopies
$f_i\tox{\alpha}f_j\tox{\beta}f_k$, the following diagram commutes.
\[
  \begin{tikzcd}
    &\Eq_B^{ij}\wedge\Eq_B^{jk}\ar[d, "\tr_B^{ijk}"]\\
    \top_A\ar[r, "\cind{\beta\hvcmp\alpha}"]\ar[ru, "\brr{\cind{\alpha},\cind\beta}"]
    &\Eq_B^{ik}\\[-10pt]
    A\ar[r, "\br{f_1,\ldots,f_n}"]&B^n
  \end{tikzcd}
\]
\pf
This follows from a diagram chase in the following diagram.
\[
  \begin{tikzcd}[baseline=(QED.base)]
    \top_A\ar[r, "\brr{\cind{\alpha},\cind\beta}"']\areq[d]
    \ar[rr, bend left=15pt, "\brr{\cind{\alpha},\cind\beta}"]
    &[15pt]\Eq_B^{ij}\wedge\Eq_B^{jk}\ar[d, "\tr_B^{ijk}"]\ar[r, "\ccind\wwdge\ccind"']
    &\Eq_B^{12}\wedge\Eq_B^{23}\ar[d, "\tr_B"]\\
    \top_A\ar[r, "\cind{\beta\hvcmp\alpha}"]\ar[rr, bend right=15pt, "\cind{\beta\hvcmp\alpha}"']&
    \Eq_B^{ik}\car[r, "\ccind"]&\Eq_B^{13}\\
    A\ar[r, "\br{f_1,\ldots,f_n}"]&B^n\ar[r, "\br{\pi_i,\pi_j,\pi_k}"]&|[alias=QED]|B^3
  \end{tikzcd}
  \tag*{\qed}
\]

\subsubsection{}\label{lem:trgenrho}\lem
For any $B\in\Ob\B$ and any $1\le{}i,j,k\le{}n$, the following diagram commutes.
\[
  \begin{tikzcd}
    &\Eq_B^{ij}\wedge\Eq_B^{jk}\ar[d, "\tr_B^{ijk}"]\\
    \top_B\ar[r, "\rho_B^{ik}"]\ar[ru, "\brr{\rho_B^{ij},\rho_B^{jk}}"]
    &\Eq_B^{ik}\\[-10pt]
    B\ar[r, "\Delta_B^n"]&B^n
  \end{tikzcd}
\]
\pf
This is proven in the same way as Lemma~\ref{lem:trgen}.
\qed

\subsubsection{}\label{lem:quadeq}\lem
For every $B\in\Ob\B$, the morphism
\[\brr{\rho_B^{12},\rho_B^{23},\rho_B^{34}}\colon\top_B\to\Eq_B^{12}\wedge\Eq_B^{23}\wedge\Eq_B^{34}\]
over $\Delta_B^4\colon{}B\to{}B^4$ is cocartesian.
\pf
This claim is obviously analogous to Proposition~\ref{prop:tripeq} and the proof is essentially the same.

We write $\brr{\rho_B^{12},\rho_B^{23},\rho_B^{34}}$ as the composite of two cocartesian morphisms
\[
  \begin{tikzcd}
    \top_B\ar[r, "\brr{\rho_B^{12},\rho_B^{23}}"]&[20pt]
    \Eq_B^{12}\wedge\Eq_B^{23}\ar[r]&\Eq_B^{12}\wedge\Eq_B^{23}\wedge\Eq_B^{34}\\[-10pt]
    B\ar[r, "\Delta_B^3"]&B^3\ar[r, "\br{\pi_1,\pi_2,\pi_3,\pi_3}"]&B^4.
  \end{tikzcd}
\]
The first morphism is cocartesian by Proposition~\ref{prop:tripeq}. The second morphism is defined in the same way as the second morphism in the proof of Proposition~\ref{prop:tripeq}, treating the codomain as a product of $\Eq_B^{12}\wedge\Eq_B^{23}$ and $\Eq_B^{34}$. We leave the details to the reader.
\qed

\subsubsection{}\thm\label{thm:homcat-is-cat}
Each hom-set $\Hom_\B(A,B)$ is the object set of a category with morphisms the \mbox{$\fib{C}$-homotopies} and composition given by vertical composition of $\fib{C}$-homotopies. We denote this category by $\HOM_\B(A,B)$. Also, the identity morphisms of $\HOM_\B(A,B)$ are the identity $\fib{C}$-homotopies.
\pf
We first show that vertical composition is associative.

Let $f,g,h,k\colon{}A\to{B}$ be morphisms in $\B$ and let
$f\tox{\alpha}g\tox{\beta}h\tox{\gamma}k$
be $\fib{C}$-homotopies. To show that $\gamma\hvcmp(\beta\hvcmp\alpha)=(\gamma\hvcmp\beta)\hvcmp\alpha$, it suffices to show that the the two composites $\top_A\to\Eq_B^{14}$ in the diagram
\begin{equation}\label{eq:pent}\begin{gathered}
    \begin{tikzcd}
      \top_A\ar[rd, "\br{\br{\cind\alpha,\cind\beta},\cind\gamma}"]\\
      &\pl(){\Eq_B^{12}\wedge\Eq_B^{23}}\wedge\Eq_B^{34}
      \ar[rr, "\mathrm{assoc}", "\sim"']
      \ar[d, "\tr_B^{123}\wedge\id_{\Eq_B^{34}}"']&&
      \Eq_B^{12}\wedge\pl(){\Eq_B^{23}\wedge\Eq_B^{34}}
      \ar[d, "\id_{\Eq_B^{12}}\wedge \tr_B^{234}"]\\
      &\Eq_B^{13}\wedge\Eq_B^{34}\ar[dr, "\tr_B^{134}"']&&
      \Eq_B^{12}\wedge\Eq_B^{g4}\ar[dl, "\tr_B^{124}"]\\
      &&\Eq_B^{14}&
    \end{tikzcd}
  \end{gathered}\end{equation}
are equal since by Lemma~\ref{lem:trgen}, they are equal to $\cind{\gamma\hvcmp(\beta\hvcmp\alpha)}$ and $\cind{(\gamma\hvcmp\beta)\hvcmp\alpha}$, respectively. In fact, we will see that the two composites $(\Eq_B^{12}\wedge\Eq_B^{23})\wedge\Eq_B^{34}\to\Eq_B^{14}$ are equal. To see this, it suffices to see that their composites with the (by Lemma~\ref{lem:quadeq}) cocartesian morphism
\[
  \brr{\brr{\rho_B^{12},\rho_B^{23}},\rho_B^{34}}\colon\top_B\to(\Eq_B^{12}\wedge\Eq_B^{23})\wedge\Eq_B^{34}
\]
are equal, and by Lemma~\ref{lem:trgenrho}, these are both equal to $\rho_B^{14}$.

Next, we must show that for each $f\colon{}A\to{B}$, the $\fib{C}$-homotopy $\hid_f$ is an identity with respect to composition of $\fib{C}$-homotopies. We will only show that it is an identity on one of the two sides, since the proof for the other side is the same.

Given another morphism $g\colon{}A\to{B}$ and a $\fib{C}$-homotopy $f\tox{\alpha}g$, we must show that $\alpha\hvcmp\hid_f=\alpha$. Now, $\alpha\hvcmp\hid_f$ is by definition the composite
\[
  \top_A
  \tox{\brr{\cind{\hid_f},\cind{\alpha}}}
  \Eq_B^{12}\wedge\Eq_B^{23}
  \tox{\tr_B}
  \Eq_B^{13}
  \tox{\ct}
  \Eq_B.
\]
Note that $\cind{\hid_f}$ and $\cind{\alpha}$ both factor through $\alpha\colon\top_A\to\Eq_B$, namely as
\[
  \begin{tikzcd}
    \top_A\ar[r, "\alpha"]&\Eq_B\ar[r, "\exx_{\pi_1}"]&\top_B\ar[r, "\rho_B^{12}"]&\Eq_B^{12}\\[-10pt]
    A\ar[r, "\br{f,g}"]&B\times{}B\ar[r, "\pi_1"]&B\ar[r, "\Delta_B^3"]&B^3
  \end{tikzcd}
  \quad\text{and}\quad
  \begin{tikzcd}
    \top_A\ar[r, "\alpha"]&\Eq_B\ar[r, "\cind{\id_{\Eq_B}}"]&[15pt]\Eq_B^{23}\\[-10pt]
    A\ar[r, "\br{f,g}"]&B\times{}B\ar[r, "\br{\pi_1,\pi_1,\pi_2}"]&B^3
  \end{tikzcd}
\]
respectively. Hence, it suffices to see that the composite
\begin{equation}\label{eq:idpfcomp}
  \Eq_B\tox{\brr{\rho_B^{12}{\exx,\cind{\id_{\Eq_B}}}}}\Eq_B^{12}\wedge\Eq_B^{23}
  \tox{\tr_B}\Eq_B^{13}
  \tox{\ct}\Eq_B
\end{equation}
is equal to the identity.

This follows from a diagram chase in the following diagram.
\[
  \begin{tikzcd}[baseline=(QED.base)]
    \top_B\rac[r, "\rho_B"']\areq[d]
    \ar[rr, bend left=15pt, "\brr{\rho_B^{12},\rho_B^{23}}"]
    &\Eq_B\ar[r, "\brr{\rho_B^{12}{\exx,\cind{\id_{\Eq_B}}}}"']&[30pt]
    \Eq_B^{12}\wedge\Eq_B^{23}\ar[d, "\tr_B"]&\\
    \top_B\ar[rr, "\rho_B^{13}"]\ar[rrr, bend right=12pt, "\rho_B"']&&
    \Eq_B^{13}\ar[r, "\ct"]&\Eq_B\\
    B\ar[r, "\Delta_B"]&B\times{}B\ar[r, "\br{\pi_1,\pi_1,\pi_2}"]&
    B^3\ar[r, "\br{\pi_1,\pi_3}"]& |[alias=QED]| B\times{}B
  \end{tikzcd}
  \tag*{\qed}
\]

\subsection{Invertibility of  homotopies}\label{subsec:c-htpy-gpd}
Continuing with a $\wedgeq$-cloven $\wedgeq$-fibration $\fibr{C}CB$, we next show that each hom-category $\HOM_\B(A,B)$ is\vspace{-3pt} in fact a groupoid.

\subsubsection{}
\defn
For an object $B\in\B$, we define $\sym_B\colon\Eq_B\to\Eq_B$ to be the unique morphism over $\br{\pi_2,\pi_1}\colon{}B\times{}B\to{}B\times{}B$ making the following diagram commute.
\[
  \begin{tikzcd}
    &\Eq_B\ar[dr, "\sym_B", dashed]&[10pt]\\[-10pt]
    \top_B\ar[rr, "\rho_B"]\rac[ru, "\rho_B"]&&\Eq_B\\[-10pt]
    B\ar[r, "\Delta_B"]&B\times{}B\ar[r, "\br{\pi_2,\pi_1}"]&B\times{}B
  \end{tikzcd}
\]

\subsubsection{}
\defn
Given morphisms $f,g\colon{}A\to{B}$ and a $\fib{C}$-homotopy $f\tox{\alpha}g$, we define its \emph{inverse}, $\alpha^{-1}$ to be the $\fib{C}$-homotopy from $g$ to $f$ given by the morphism $\sym_B\ccmp\,\alpha$, as shown below.
\[
  \begin{tikzcd}
    \top_A\ar[r, "\alpha"]&\Eq_B\ar[r, "\sym_B"]&[10pt]\Eq_B\\[-10pt]
    A\ar[r, "\br{f,g}"]&B\times{}B\ar[r, "\br{\pi_2,\pi_1}"]&B\times{}B
  \end{tikzcd}
\]

\subsubsection{}\label{thm:groupoid}\thm
Given morphisms $f,g\colon{}A\to{B}$ in $\B$, every $\fib{C}$-homotopy $f\tox{\alpha}g$ is invertible with respect to vertical composition of $\fib{C}$-homotopies, with inverse $\alpha^{-1}$.
\pf
We will only show that $\alpha^{-1}$ is an inverse on one side of $\alpha$, since the proof is identical for the other side.

We must show that the composite $\ct\ccmp{\tr_B}\ccmp{\brr{\cind\alpha,\cind{\alpha\I}}}$ (shown below) is equal to $\hid_f$.
\[
  \begin{tikzcd}
  \top_A\ar[r, "\brr{\cind{\alpha},\cind{\alpha^{-1}}}"]&[15pt]
  \Eq_B^{12}\wedge\Eq_B^{23}\ar[d, "\tr_B"]&\\
  &\Eq_B^{13}\ar[r, "\ct"]&
  \Eq_B\\[-10pt]
  A\ar[r, "\br{f,g,f}"]&B^3\ar[r, "\br{\pi_1,\pi_3}"]&B\times{}B
\end{tikzcd}
\]

Note that both $\cind{\alpha}$ and $\cind{\alpha\I}$ factor through $\alpha\colon\top_A\to\Eq_B$, namely as
\[
  \begin{tikzcd}
    \top_A\ar[r, "\alpha"]&\Eq_B\ar[r, "\cind{\id_{\Eq_B}}"]&[15pt]\Eq_B^{12}\\[-10pt]
    A\ar[r, "\br{f,g}"]&B\times{}B\ar[r, "\br{\pi_1,\pi_2,\pi_1}"]&B^3
  \end{tikzcd}
  \quad\text{and}\quad
  \begin{tikzcd}
    \top_A\ar[r, "\alpha"]&\Eq_B\ar[r, "\cind{\sym_B}"]&[15pt]\Eq_B^{23}\\[-10pt]
    A\ar[r, "\br{f,g}"]&B\times{}B\ar[r, "\br{\pi_1,\pi_2,\pi_1}"]&B^3
  \end{tikzcd}
\]
respectively. Hence, using that $\ct\cind{\hid_{\pi_1}}\exx_{\br{f,g}}=\hid_{\pi_1}\exx_{\br{f,g}}=\hid_{f}$, it suffices to show that the square in the following diagram commutes.
\[
  \begin{tikzcd}[column sep=50pt]
    \top_A\ar[r, "\alpha"]\ar[rd, "\exx_{\br{f,g}}"']&
    \Eq_B\ar[r, "\brr{\cind{\id},\cind{\sym_B}}"]\ar[d, "!"']&
    \Eq_B^{12}\wedge\Eq_B^{23}\ar[d, "\tr_B"]\\
    &\top_{B\times{}B}\ar[r, "\cind{\hid_{\pi_1}}"]&\Eq_B^{13}\ar[r, "\ct"]&\Eq_B\\[-10pt]
    A\ar[r, "\br{f,g}"]&B\times{}B\ar[r, "\br{\pi_1,\pi_2,\pi_1}"]&
    B^3\ar[r, "\br{\pi_1,\pi_3}"]&B\times{}B
  \end{tikzcd}
\]

This follows from a diagram chase in the following diagram.
\[
  \begin{tikzcd}[column sep=50pt, baseline=(QED.base)]
    \top_B\areq[d]\rac[r, "\rho_B"']\ar[rr, "\brr{\rho_B^{12},\rho_B^{23}}", bend left=15pt]&
    \Eq_B\ar[r, "\brr{\cind{\id},\cind{\sym_B}}"']\ar[d, "!"']&
    \Eq_B^{12}\wedge\Eq_B^{23}\ar[d, "\tr_B"]\\
    \top_B\ar[r, "\exx_{\Delta_B}"]\ar[rr, "", bend right=15pt, "\rho_B^{13}"']
    &\top_{B\times{}B}\ar[r, "\cind{\hid_{\pi_1}}"]&\Eq_B^{13}\\[5pt]
    B\ar[r, "\Delta_B"]&B\times{}B\ar[r, "\br{\pi_1,\pi_2,\pi_1}"]& |[alias=QED]| B^3
  \end{tikzcd}
  \tag*{\qed}
\]

\subsection{The 2-categorical structure}\label{subsec:c-htpy-2cat}
We continue with a fixed $\wedgeq$-cloven $\wedgeq$-fibration $\fib{C}$.

\subsubsection{}
\defn
Given morphisms $h,k\colon{}B\to\C$ in $\C$ and a $\fib{C}$-homtopy
$h\tox{\beta}k$,
we denote by $\check{\beta}$ the unique morphism \mbox{$\Eq_B\to\Eq_C$}\vspace{-10pt} over $h\times{}k$ making the following diagram commute.

\[
  \begin{tikzcd}
    &\Eq_B\rac[dr, "\check\beta", dashed]&\\
    \top_B\ar[rr, "\beta"]\ar[ru, "\rho_B"]&&\Eq_C\\[-10pt]
    B\ar[r, "\Delta_B"]&B\times{}B\ar[r, "h\times{}k"]&C\times{}C
  \end{tikzcd}
\]

\subsubsection{}\label{defn:horiz-comp-c-htpy}
\defn
Given $\fib{C}$-homotopies
\[
  \begin{tikzcd}[row sep=0pt, column sep=10pt]
    &\ar[dd, shorten >=2pt, shorten <=-2pt, Rightarrow, "\ \alpha", pos=0.3]&
    &\ar[dd, shorten >=2pt, shorten <=-2pt, Rightarrow, "\ \beta", pos=0.4]&\\[-5pt]
    A\ar[rr, "f", bend left]\ar[rr, "g"', bend right]&&
    B\ar[rr, "h", bend left]\ar[rr, "k"', bend right]&&C,\\
    &{}&&{}&
  \end{tikzcd}
\]
we define the \emph{horizontal composite} of $\alpha$ and $\beta$, which we denote by $\beta\hhcmp\alpha$, to be the \mbox{$\fib{C}$-homotopy} from $hf$ to $kg$ given by the composite
\[
  \top_A\tox{\alpha}\Eq_B\tox{\check\beta}\Eq_C.
\]

\subsubsection{}\label{prop:2-assoc}\prop
Horizontal composition is associative. That is, given morphisms of $\B$ and $\fib{C}$-homotopies
\[
  \begin{tikzcd}[row sep=0pt, column sep=10pt]
    &\ar[dd, shorten >=2pt, shorten <=-2pt, Rightarrow, "\ \alpha", pos=0.4]&
    &\ar[dd, shorten >=2pt, shorten <=-2pt, Rightarrow, "\ \beta", pos=0.4]&
    &\ar[dd, shorten >=2pt, shorten <=-2pt, Rightarrow, "\ \gamma", pos=0.4]&\\[-5pt]
    A\ar[rr, "f", bend left]\ar[rr, "g"', bend right]&&
    B\ar[rr, "h", bend left]\ar[rr, "k"', bend right]&&
    C\ar[rr, "l", bend left]\ar[rr, "m"', bend right]&&D\\
    &{}&&{}&&{}&
  \end{tikzcd}
\]
we have $(\gamma\hhcmp\beta)\hhcmp\alpha=\gamma\hhcmp(\beta\hhcmp\alpha)$.
\pf
We need to show that $(\widecheck{\gamma\hhcmp\beta})\alpha=\check\gamma\check\beta\alpha$. Hence, it suffices to show that $\widecheck{\gamma\hhcmp\beta}=\check\gamma\check\beta$. This follows from a diagram chase in the following diagram.
\[
  \begin{tikzcd}[baseline=(QED.base)]
    &\Eq_B\ar[dr, "\check\beta"]\ar[rrdd, "\widecheck{\gamma\hhcmp\beta}", bend left]&&\\[-10pt]
    &&\Eq_C\ar[rd, "\check\gamma", near start]&\\[-10pt]
    \top_B\rac[uur, "\rho_B"]\ar[rru, "\beta"]\ar[rrr, "\gamma\hhcmp\beta"]&&&\Eq_D\\[-5pt]
    B\ar[r, "\Delta_B"]&B\times{}B\ar[r, "h\times{}k"]&C\times{}C\ar[r, "l\times{}m"]&
    |[alias=QED]| D\times{}D
  \end{tikzcd}
  \tag*{\qed}
\]

\subsubsection{}\label{prop:comp-bifun}\prop
For any objects $A,B,C\in\B$, horizontal composition extends the composition map $\Hom_\B(A,B)\times\Hom_\B(B,C)\to\Hom_\B(A,C)$ to a bifunctor $\HOM_\B(A,B)\times\HOM_\B(B,C)\to\HOM_\B(A,C)$.

This means that for morphisms $A\tox{f}B\tox{k}C$ in $\B$, we have $\hid_k\hhcmp\hid_f=\hid_{kf}$, and that given further morphisms and $\fib{C}$-homotopies
\[
  \begin{tikzcd}[column sep=10pt]
    &\ar[d, shorten >=3pt, shorten <=10pt, Rightarrow, "\ \alpha", pos=0.6]&
    &\ar[d, shorten >=3pt, shorten <=10pt, Rightarrow, "\ \gamma", pos=0.6]&\\[-5pt]
    A\ar[rr, "f", bend left=55pt]\ar[rr, "h"', bend right=55pt]\ar[rr, "g", pos=0.1]&{}
    \ar[d, shorten >=14pt, shorten <=-2pt, Rightarrow, "\ \beta", pos=0.1]&
    B\ar[rr, "k", bend left=55pt]\ar[rr, "m"', bend right=55pt]\ar[rr, "l", pos=0.2]&{}
    \ar[d, shorten >=14pt, shorten <=-2pt, Rightarrow, "\ \delta", pos=0.1]&C\\
    &{}&&{}&
  \end{tikzcd}
\]
we have $(\delta\hvcmp\gamma)\hhcmp(\beta\hvcmp\alpha)=
(\delta\hhcmp\beta)\hvcmp(\gamma\hhcmp\alpha)$.
\pf
The first claim follows from the commutativity of
\[
\begin{tikzcd}
  &&\Eq_B\ar[rrd, "\widecheck{\hid_k}"]&&\\
  &\top_B\ar[ru, "\rho_B"']\ar[rrd, "\exx_{k}", pos=0.7]\ar[rrr, "\hid_k"]&&&\Eq_C\\
  \top_A\ar[ru, "\exx_f"']\ar[rrr, "\exx_{kf}"]\ar[rruu, bend left, "\hid_f"]
  \ar[rrrru, "\hid_{kf}", rounded corners,
      to path={-- ([xshift=15pt,yshift=-20pt]\tikztostart.east)
               -- ([xshift=-20pt,yshift=-45pt]\tikztotarget.south)\tikztonodes
               -- (\tikztotarget)}]
  \ar[rrrru, "\hid_k\hhcmp\hid_f", rounded corners,
      to path={-- ([xshift=15pt,yshift=80pt]\tikztostart.north)
               -- ([xshift=-20pt,yshift=45pt]\tikztotarget.north)\tikztonodes
               -- (\tikztotarget)}]
  &&&\top_C\ar[ru, "\rho_C"]&\\[10pt]
  A\ar[r, "f"]&B\ar[rr, "k"]\ar[rd, "\Delta_B"']&&C\ar[r, "\Delta_C"]&C\times{}C.\\[-8pt]
  &&B\times{}B\ar[rru, "k\times{}k"']&&
\end{tikzcd}
\]

For the second claim, we must show that the following two composites are equal.
\[
  \top_A\tox{\brr{\cind{\alpha},\cind{\beta}}}
  \Eq_B^{12}\wedge\Eq_B^{23}\tox{\tr_B}
  \Eq_B^{13}\tox{\ct}
  \Eq_B\tox{\widecheck{\delta\hvcmp\gamma}}
  \Eq_C
\]
\[
  \top_A\tox{\brr{\cind{\gamma\hhcmp\alpha},\cind{\delta\hhcmp\beta}}}
  \Eq_C^{12}\wedge\Eq_C^{23}\tox{\tr_C}
  \Eq_C^{13}\tox{\ct}
  \Eq_C
\]

We first note that $\cind{\gamma\hhcmp\alpha}$ and $\cind{\delta\hhcmp\beta}$ can be factored respectively as
\[
  \top_A\tox{\cind{\alpha}}
  \Eq_B^{12}\tox{\cind{\check\gamma\ccmp\ct}}
  \Eq_C^{12}
  \quad\text{and}\quad
  \top_A\tox{\cind{\beta}}
  \Eq_B^{23}\tox{\cind{\check\delta\ccmp\ct}}
  \Eq_C^{23}.
\]

Hence, it suffices to see that the following diagram commutes.
\[
  \begin{tikzcd}
    \Eq_B^{12}\wedge\Eq_B^{23}\ar[d, "\tr_B"]
    \ar[rr, "\cind{\check\gamma\ccmp\ct}\ \wwdge\  \cind{\check\delta\ccmp\ct}"]&&
    \Eq_C^{12}\wedge\Eq_C^{23}\ar[d, "\tr_C"]&\\
    \Eq_B^{13}\ar[dr, "\ct"]&&\Eq_C^{13}\ar[dr, "\ct"]\\
    &\Eq_B\ar[rr, "\widecheck{\delta\hvcmp\gamma}"]&&\Eq_C\\
    B^3\ar[rr, "k\times{}l\times{}m"]\ar[dr, "\br{\pi_1,\pi_3}"']&&
    C^3\ar[dr, "\br{\pi_1,\pi_3}"]\\
    &B\times{}B\ar[rr, "k\times{}m"]&&C\times{}C
  \end{tikzcd}
\]

This follows from a diagram chase in the following diagram.
\[
  \begin{tikzcd}[baseline=(QED.base)]
    \top_B\rac[r, "\brr{\rho_B^{12},\rho_B^{23}}", pos=0.65]
    \ar[rrdd, "\rho_B", bend right]
    \ar[dr, "\rho_B^{13}"', bend right=5pt, near end]
    \ar[rrr, "\brr{\cind\gamma,\cind\delta}", bend left=25pt]&[30pt]
    \Eq_B^{12}\wedge\Eq_B^{12}\ar[d, "\tr_B"]
    \ar[rr, "\cind{\check\gamma\ccmp\ct}\ \wwdge\ \cind{\check\delta\ccmp\ct}"]&&
    \Eq_C^{12}\wedge\Eq_C^{23}\ar[d, "\tr_C"]&\\
    &\Eq_B^{13}\ar[dr, "\ct"]&&\Eq_C^{13}\ar[dr, "\ct"]\\
    &&\Eq_B\ar[rr, "\widecheck{\delta\hhcmp\gamma}"]&& |[alias=QED]| \Eq_C
  \end{tikzcd}
  \tag*{\qed}
\]

\subsubsection{}\label{thm:2cat}\thm
$\B$ can be extended to a 2-category, in which the hom-categories are the categories $\HOM_\B(A,B)$, and in which the composition of 2-cells is given by horizontal composition of $\fib{C}$-homotopies.
\pf
This is precisely the content of Propositions~\ref{prop:2-assoc}~and~\ref{prop:comp-bifun}, together with the following fact: given morphisms in $\B$ and a $\fib{C}$-homotopy as in
\[
  \begin{tikzcd}[row sep=0pt]
    &&[-5pt]\ar[dd, shorten >=2pt, shorten <=-2pt, Rightarrow, "\ \alpha", pos=0.4]&[-5pt]&\\[-5pt]
    A\ar[r, "\id_A"]&B,\ar[rr, "f", bend left]\ar[rr, "g"', bend right]&&B\ar[r, "\id_B"]&B\\
    &&{}&&
  \end{tikzcd}
\]
we have $\alpha\hhcmp{\hid_{\id_A}}=\alpha=\hid_{\id_B}\hhcmp\alpha$. To see this last fact, note that $\hid_{\id_A}$ and $\hid_{\id_B}$ are just $\rho_A$ and $\rho_B$, respectively. It follows that $\widecheck{\hid_{\id_B}}=\id_{\Eq_B}$. Hence
\[
  \alpha\hhcmp\hid_{\id_A}=\check\alpha\ccmp\rho_A=\alpha
  \quad\quad\text{and}\quad\quad
  \hid_{\id_B}\hhcmp\alpha=\widecheck{\hid_{\id_B}}\ccmp\alpha=\id_{\Eq_B}\ccmp\,\alpha=\alpha.
  \tag*{\qed}
\]

\subsection{Internal categories}\label{subsec:intcats}
In this section, we give an alternative presentation of the 2-categorical structure on $\B$ which was suggested by the referee.

Recall that for $A,B\in\Ob\B$ the morphisms in the category $\HOM_{\B}(A,B)$ were defined as all morphisms $\top_A\to\Eq_B$ in $\C$. Now, in general, when a hom-set in a category carries some algebraic structure, it is often an indication that the codomain object itself carries this structure (in the ``internal'', categorical sense). In this case, this suggests that $\Eq_B$ is the object of arrows of an internal category in $\C$. This is indeed the case, and pursuing this leads to the following characterization of the 2-categorical structure on $\B$ (Theorem~\ref{thm:intcat-descr}): there is a fully-faithful functor from $\B$ to $\Cat(\C)$, the category of internal categories in $\C$, and the 2-categorical structure on $\B$ is the one induced via this functor from the natural 2-categorical structure on $\Cat(\C)$. Below, we carry out this construction.

We note that Theorem~\ref{thm:groupoid} on the invertibility of homotopies can also be stated in this language -- namely, as saying that the internal categories considered here are in fact ``internal groupoids''.

The notion of internal category, and the associated notions such as internal functor, as well as the 2-categorical structure on $\Cat(\C)$, are defined, and their properties proved, in the usual manner: by simplemindedly imitating the usual definitions, replacing ``Set'' by ``object of $\C$'' (see, e.g, \cite[p.~267]{cwm}).

\subsubsection{}
Let us fix a few more notational conventions. Given a pair of morphisms $A\tox{f}B\xot{g}C$ in a category $\C$, we may use the notation $A\times_BC$ for a pullback (suppressing, as is common, the dependence on $f$ and $g$), and in this case, we denote by $\pi_1\colon{}A\times_BC\to{}A$ and $\pi_2\colon{}A\times_BC\to{}C$ the associated projections, and given morphisms $p\colon{}D\to{}A$ and $q\colon D\to{}C$ with $fp=qg$, we denote by $\br{p,q}\colon{}D\to{}A\times_BC$ the induced morphism.

\subsubsection{}\label{defn:intcat}
\defn
Let $\C$ be a category. An \emph{internal category} $C$ in $\C$ consists of (i) a reflexive graph object in $\C$ -- i.e., objects and morphisms $C_1\underset{t}{\overset{s}{\rightrightarrows}} C_0 \tox{\iota} C_1$ such that $s\iota=t\iota=\id_{C_0}$ -- together with (ii) a choice $C_2=C_1\times_{C_0}C_1$ of pullback of the morphisms $C_1\tox{t}C_0\xot{s}C_1$, and (iii) a morphism $C_2\tox{m}C_1$. These are moreover required to satisfy the usual unitality and associativity conditions (namely $m\br{\iota s,\id_{C_1}}=\id_{C_1}=m\br{\id_{C_1},\iota t}$ and $m\br{m\br{\pi_1,\pi_2},\pi_3}=m\br{\pi_1,m\br{\pi_2,\pi_3}}\colon{}C_3\to C_1$, respectively). The latter involves the choice of a triple pullback $C_3=C_1\times_{C_0}C_1\times_{C_0}C_1$ (i.e., a limit of the diagram $C_1\tox{t}C_0\xot{s}C_1\tox{t}C_0\xot{s}C_1$), but is of course independent of this choice.

Given internal categories $C$ and $D$, an \emph{internal functor} $F\colon{}C\to D$ consists of morphisms $F_i\colon{}C_i\to D_i$ ($i=1,2$) making the resulting squares involving $s,t,\iota,m$ commute (the last of these involves the induced morphism $F_2=\br{F_1\pi_1,F_1\pi_2}\colon{}C_2\to D_2$) -- i.e., preserving domain, codomain, identities, and composition.

Given two internal functors $F,G\colon{}C\to{}D$, an \emph{internal natural transformation} $\alpha\colon{}F\to{}G$ is a morphism $\alpha_0\colon{}C_0\to D_1$ with $s\alpha_0=F$ and $t\alpha_0=G$ satisfying the usual naturality condition $m\br{F_1,\alpha_0 t}=m\br{\alpha_0 s, G_1}\colon{}C_1\to D_1$.

Internal categories of $\C$ and internal functors form a category $\Cat(\C)$. In fact, $\Cat(\C)$ is naturally a 2-category in the following way, generalizing the usual 2-category of categories (see, e.g., \cite[Proposition~8.1.4]{borceux1}).

The 2-cells are given by the internal natural transformations. Given internal functors $F,G,H\colon{}C\to{}D$ and internal natural transformations $\alpha\colon{}F\to{}G$ and $\beta\colon{}G\to{}H$, the vertical composite of $\alpha$ and $\beta$ is given by the composite morphism $C_0\tox{\br{\alpha_0,\beta_0}}D_2\tox{m}D_1$.

Next, given internal functors and internal natural transformations as in
\[
  \begin{tikzcd}[row sep=0pt, column sep=10pt]
    &\ar[dd, shorten >=2pt, shorten <=-2pt, Rightarrow, "\ \alpha", pos=0.3]&
    &\ar[dd, shorten >=2pt, shorten <=-2pt, Rightarrow, "\ \beta", pos=0.4]&\\[-5pt]
    C\ar[rr, "F", bend left]\ar[rr, "G"', bend right]&&
    D\ar[rr, "H", bend left]\ar[rr, "K"', bend right]&&E,\\
    &{}&&{}&
  \end{tikzcd}
\]
the horizontal composite of $\alpha$ and $\beta$ is given by the composite $C_0\tox{\alpha_0}D_1\tox{\br{H_1,\beta_0 t}}E_2\tox{m}E_1$ (or equivalently $m\br{\beta_0 s,K_1}\alpha_0$).

The verification that this indeed defines a 2-category is a (somewhat lengthy but) straightforward exercise that we leave to the reader.

\subsubsection{}
\defn
Let us now fix a $\wedgeq$-cloven $\wedgeq$-fibration $\fibr{C}CB$ until the end of \S\ref{subsec:intcats}.

We define, for each object $B\in\Ob\B$, a category object $\intind{B}$ in $\C$ as follows. The underlying reflexive graph is $\Eq_B\underset{\exx_{\pi_2}}{\overset{\exx_{\pi_1}}{\rightrightarrows}} \top_B \tox{\rho_B} \Eq_B$. For the pullback $\intind{B}_2=\Eq_B\times_{\top_B}\Eq_B$, we take $\Eq_B^{12}\wedge\Eq_B^{23}$ with projections $\Eq_B\xot{\ct\pi_1}\Eq_B^{12}\wedge\Eq_B^{23}\tox{\ct\pi_2}\Eq_B$ (to see that this is indeed a pullback, note that given $\Eq_B\xot{p}P\tox{q}\Eq_B$ with $\exx_{\pi_1}p=\exx_{\pi_2}q$, the induced morphism $P\to\Eq_B^{12}\wedge\Eq_B^{23}$ is just $\brr{\cind{p},\cind{q}}$). We define the composition morphism $m\colon\Eq_B^{12}\wedge\Eq_B^{23}\to\Eq_B$ to be $\ct\tr_B$.

The unitality property now follows from \eqref{eq:idpfcomp} in the proof of Theorem~\ref{thm:homcat-is-cat}.

Similarly, for the associativity property, we can take as a triple pullback $C_1\times_{C_0}\times C_1\times_{C_0}\times C_1$ the object $\Eq_B^{12}\wedge\Eq_B^{23}\wedge\Eq_B^{34}$ with projections $\ct\pi_1,\ct\pi_2,\ct\pi_3$, and the associativity property can then be seen to follow from the commutativity of \eqref{eq:pent} in the proof of Theorem~\ref{thm:homcat-is-cat}.

Next, for each morphism $f\colon{}B\to{}C$ in $\B$, we define an internal functor $\intind{f}\colon\intind{B}\to\intind{C}$ by taking $\intind{f}_0$ to be $\exx_f\colon\top_B\to\top_C$ and $\intind{f}_1$ to be $\widecheck{\hid_f}\colon\Eq_B\to\Eq_C$. Preservation of domain and codomain is immediate, and preservation of identities follows from the definition of $\widecheck{\hid_f}$.

For the preservation of composites, we consider the following diagram, in which we would like to see that the rightmost parallelogram (in $\C$) commutes.
\[
  \begin{tikzcd}[column sep=10pt]
    \top_B\ar[rr, "\rho_B"]\ar[rd, "\exx_f"']&&
    \Eq_B\ar[rr, "\brr{\cind{\rho_B^{12}}\exx_{\pi_2},\cind{\id}}"]
    \ar[rd, "\widecheck{\id_f}"]
    &&
    \Eq_B^{12}\wedge\Eq_B^{23}\ar[rr, "\ct\tr_B"]
    \ar[rd, "\cind{\widecheck{\id_f}\ct}\wwdge\cind{\widecheck{\id_f}\ct}" near end]&[-30pt]&
    \Eq_B
    \ar[rd, "\widecheck{\id_f}"]\\
    &\top_C\ar[rr, "\rho_C"]&&
    \Eq_C\ar[rr, "\brr{\cind{\rho_C^{12}}\exx_{\pi_2},\cind{\id}}"]&&
    \Eq_C^{12}\wedge\Eq_C^{23}\ar[rr, "\ct\tr_C"]&&
    \Eq_C\\[-10pt]
    B\ar[rr, "\Delta_B"]\ar[rd, "f"]&&
    B^2\ar[rr, "\Delta_B\times\id_B"]\ar[rd, "f\times{}f"]&&
    B^3\ar[rr, "\br{\pi_1,\pi_2}"]\ar[rd, "f\times{}f\times{}f"]&&B^2\ar[rd, "f\times{}f"]\\
    &C\ar[rr, "\Delta_C"]&&C^2\ar[rr, "\Delta_C\times\id_C"]&&C^3\ar[rr, "\br{\pi_1,\pi_2}"]&&C^2
  \end{tikzcd}
\]
Since the first two parallelograms commute and the composite of the first two morphisms in the first row is cocartesian, it suffices to see that the outside of the diagram commutes, but this is so since the composites of the second and third horizontal morphisms in the first two rows are identity morphisms.

Next, given a $\fib{C}$-homotopy $\alpha\colon\top_B\to\Eq_C$ from $f$ to $g$, we wish to define an internal natural transformation $\intind{\alpha}\colon\intind{f}\to\intind{g}$. For this, we need a morphism $\intind{\alpha}_0\colon\top_B=\intind{B}_0\to\intind{C}_1=\Eq_C$, and we just take $\alpha$ itself. It is immediate that $s\intind\alpha_0=\intind{f}_0$ and $t\intind\alpha_0=\intind{g}_0$.

For the naturality condition, we need to verify the equality $\ct\tr_C\brr{\cind{\widecheck{\hid_f}},\cind{\alpha\exx_{\pi_2}}}=
\ct\tr_C\brr{\cind{\alpha\exx_{\pi_1}},\cind{\widecheck{\hid_g}}}\colon\Eq_B\to\Eq_C$. But after precomposing with the cocartesian morphism $\rho_B\colon\top_B\to\Eq_B$, both sides are equal to $\alpha$ by the part of Theorem~\ref{thm:homcat-is-cat} concerning identity $\fib{C}$-homotopies.

\subsubsection{}\label{thm:intcat-descr}
\thm
The assignments $B\mapsto\intind{B}$, $f\mapsto\intind{f}$, and $\alpha\mapsto\intind{\alpha}$ define a 2-functor $\B\to\Cat(\C)$ which is ``2-fully-faithful'' (i.e., it induces isomorphisms $\HOM_{\B}(B,C)\to\HOM_{\Cat(\C)}(\intind{B},\intind{C})$).

In other words, the 2-category structure from Theorem~\ref{thm:2cat} admits the following alternative characterization (up to isomorphism): it is obtained by considering $\B$ as a full subcategory of $\Cat(\C)$ via the fully faithful functor given by $B\mapsto\intind{B}$ and $f\mapsto\intind{f}$, and then passing to the full sub-2-category of (the 2-category) $\Cat(\C)$ on the objects and morphisms of $\B$.
\pf
That this defines a (1-)functor follows from the fact that $\widecheck{\hid_{\id_B}}=\id_{\Eq_B}$, as shown in the proof of Theorem~\ref{thm:2cat}, and the fact that $\widecheck{\id_g}\cdot\widecheck{\id_f}=\widecheck{\id_{gf}}$, which follows immediately from the definitions.

That this functor is faithful follows from the obvious fact that $\exx_f=\exx_g$ implies $f=g$.

To see that it is full, we need to show that any internal functor $F=(F_0,F_1)\colon\intind{B}\to\intind{C}$ is equal to $\intind{f}$ for some $f$. Clearly, we must have $F_0=\exx_f\colon\top_B\to\top_C$ where $f$ is the morphism over which $F_0$ lies. To see that $F_1=\widecheck{\id_f}$, we note first that $F_1$ must lie over $f\times{}f$, and then the asserted equality follows from the cocartesianness of $\rho_B$ and the fact that $F$ preserves identities.

That the assignment $\alpha\mapsto\intind{\alpha}$ defines a bijection from the set of $\fib{C}$-homotopies $f\to{}g$ to the set of internal natural transformations $\intind{f}\to\intind{g}$ is clear since $\intind{\alpha}_0=\alpha$.

It now remains to see that the operation $\alpha\mapsto\intind{\alpha}$ preserves horizontal and vertical composition.

Given $\fib{C}$-homotopies $\alpha\colon\top_A\to{}\Eq_B$ from $f$ to $g$ and $\beta\colon\top_B\to\Eq_C$ from $h$ to $k$, we have
\[
  \beta\hhcmp\alpha=\check\beta\cdot\alpha\quad\text{and}\quad
  (\intind\beta\circ\intind\alpha)_0
  =\ct\tr_B\brr{\cind{\widecheck{\id_h}},\cind{\beta\exx_{\pi_2}}}\alpha.
\]
To see that these are equal, it suffices to show that $\ct\tr_B\brr{\cind{\widecheck{\id_h}},\cind{\beta\exx_{\pi_2}}}=\widecheck{\beta}\colon\Eq_B\to\Eq_C$, which we can check after precomposing with the cocartesian $\rho_B\colon\top_B\to\Eq_B$. Hence, it remains to prove
$\ct\tr_B\brr{\cind{\widecheck{\id_h}},\cind{\beta\exx_{\pi_2}}}\rho_B=\beta$. But we have
\[
  \text{LHS}=
  \ct\tr_B\brr{\cind{\widecheck{\id_h}\rho_B},\cind{\beta\exx_{\pi_2}\rho_B}}=
  \ct\tr_B\brr{\cind{\id_h},\cind{\beta}}=\beta.
\]

Next, given $\fib{C}$-homotopies $\alpha\colon\top_A\to\Eq_B$ from $f$ to $g$ and $\beta\colon\top_B\to\Eq_C$ from $g$ to $h$, we have
\[
  \beta\hvcmp\alpha=\ct\tr_B\brr{\cind\alpha,\cind\beta},
\]
which, after unfolding the definitions, is seen to be precisely the vertical composite of $\intind{\alpha}$ and $\intind{\beta}$.
\qed

\subsection{1-discrete 2-fibrations}\label{subsec:1d2fs}
In this section, we will prove that not only the base of a $\wedgeq$-cloven fibration $\fibr{C}CB$, but also the total category, can be given a natural 2-categorical structure, so that $\fib{C}$ becomes a \emph{1-discrete 2-fibration} (Definition~\ref{defn:1d2f}) or \emph{1D2F}. Let us briefly comment on the significance of this notion.

We recall that a fibration is called \emph{discrete} if each fiber is a discrete category (every morphism is an identity morphism). As we mentioned in the introduction, the notion of fibration was introduced simultaneously with the essentially equivalent notion of \emph{pseudo-functor}, the idea behind which is that a fibration $\fib{C}$ over $\B$ can be described instead by some kind of functor $\B^\op\to\Cat$ taking each object $A\in\Ob\B$ to its fiber $\fib{C}^A$. Under this correspondence, the discrete fibrations correspond exactly to functors $\B^\op\to\Cat$ factoring through $\Set\hookrightarrow\Cat$, and the operation taking a pseudo-functor to its associated fibration (the ``\emph{Grothendieck construction}'') recovers in this case the well-known \emph{category of elements} of a presheaf.

In passing from the notion of discrete fibration to that of fibration, one replaces the category $\Set$ with the 2-category $\Cat$, and thus increases the ``categorical dimension'' by one. From this point of view, it is therefore most natural to consider morphisms $\B\to\Cat$ not only from a \emph{1-category} $\B^\op$, but from a general \emph{2-category}. Applying the appropriate generalization of the Grothendieck construction, we should then obtain a special kind of 2-functor with codomain $\B$ whose ``fibers'' are all 1-categories. This is the notion of 1D2F.

Extending the fibration $\fib{C}$ to a 1D2F is thus tantamount to extending the associated pseudo-functor $\B^\op\to\Cat$ to the 2-categorical structure on $\B$. Concretely, this means that, not only should each morphism in $\B$ induce a pullback functor between the appropriate fibers, but also each 2-cell in $\B$ should induce a \emph{natural transformation} between the associated pullback functors.

Thus, for each $\fib{C}$-homotopy $f\tox{\alpha}g$ between a pair of morphisms in $f,g\colon{}A\to{}B$ in $\B$ and each $P\in\Ob\fib{C}^B$, we need to produce a certain morphism $f^*P\to{}g^*P$. From the logical point of view, this amounts to producing a proof of $P(f(a))\To{}P(g(a))$ for each predicate $P(b)$, given a proof of $f(a)=g(a)$. This, in turn, can be reduced to proving $(b_1=b_2\wedge{}P(b_1))\To{}P(b_2)$, which in terms of the fibration means producing a morphism $\pi_1^*P\wedge\Eq_B\to\pi_2^*P$ in $\fib{C}^{B\times{}B}$, and this is precisely how we will proceed.

Finally, we note that, just as one does not need the general notion of fibration to define that of a discrete fibration, one can directly define 1-discrete 2-fibrations without defining 2-fibrations in general, and this is what we do. Similarly, constructing 1D2Fs is simpler than the general task of constructing 2-functors. In Lemma~\ref{lem:constructing-1d2fs}, we explain exactly what one needs in order to extend a fibration to a 1-discrete 2-fibration (given a 2-category structure on the base).

\subsubsection{}\label{defn:1d2f}
\defn
A \emph{pre-2-fibration} is simply a 2-functor $\fibr{C}CB$. We use similar terminology for pre-2-fibrations as we do for prefibrations: $\B$ is the \emph{base} 2-category; $\C$ is the \emph{total} 2-category; a 0-, 1-, or 2-cell in $\C$ \emph{lies over} its image in $\B$; and so on. The \emph{fiber} $\fib{C}^A$ of $\fib{C}$ over $A\in\Ob\B$ is the sub-2-category consisting of 0-cells, 1-cells, and 2-cells lying over $A$, $\id_A$, and $\id_{\id_A}$, respectively. The \emph{underlying prefibration} of a pre-2-fibration is just the induced functor on the underlying 1-categories.

The pre-2-fibration $\fib{C}$ is a \emph{1-discrete 2-fibration}\footnote{These have also considered by M. Lambert in \cite[Definition~2.2.15]{lambertthesis}, where they are called ``discrete 2-fibrations'', and were also known to our anonymous referee, who called them ``locally discrete fibrations''. We learned the concept from M. Makkai.} (or 1D2F) if (i) the underlying prefibration is a fibration, and (ii) for every 2-cell $\alpha\colon{}f\to{}g$ in $\B$ and every 1-cell $p$ over $f$, there is a unique 2-cell over $\alpha$ with domain $p$, as depicted below (this says that the functor $\HOM_\C(P,Q)\to\HOM_\B(A,B)$ is a \emph{discrete op-fibration}).
\[
  \begin{tikzcd}[column sep=30pt]
    P
    \ar[r, bend left, "p"{name=p}]
    \ar[r, bend right, ""'{name=q}, dashed]
    &Q
    \ar[Rightarrow, from=p, to=q, dashed, shorten <=3pt, shorten >=3pt]
    \\[10pt]
    A
    \ar[r, bend left, "f"{name=f}]
    \ar[r, bend right, "g"'{name=g}]
    &B
    \ar[Rightarrow, from=f, to=g, "\alpha", shorten <=3pt, shorten >=3pt]
  \end{tikzcd}
\]

Note that the fibers of a 1D2F are 1-categories and that, if $\B$ is a 1-category, then a 1D2F over $\B$ (seen as a 2-category with only identity 2-cells) is the same thing as a fibration over $\B$.

\subsubsection{}\label{lem:constructing-1d2fs}
\lem
Let $\fibr{C}CB$ be a fibration, with $\B$ a 2-category. Then any extension $\fibr{D}DB$ of $\fib{C}$ to a 1D2F is determined up to isomorphism by the function $F$ whose domain is the set of pairs $(\alpha,p)$ consisting of a 2-cell $\alpha\colon{}f\to{}g$ in $\B$ and a lift $p$ of $f$ in $\C$, and which assigns to $(\alpha,p)$ the codomain of the unique lift of $\alpha$ in $\D$ with domain $p$.

That is, given another extension $\fibr{D'}{D'}B$ with the same associated function $F$, there is a unique isomorphism of 2-categories $\D\to\D'$ extending $\id_\C$ and commuting with the 2-functors $\fib{D}$, $\fib{D'}$.

Moreover, an arbitrary function $F$ assigning to each pair $(\alpha,p)$ as above some lift of the codomain of $\alpha$ comes from a 1D2F if and only if it satisfies the following three conditions:
\begin{enumerate}[(i)]
\item If $\alpha\colon{}f\to{}g$ and $\beta\colon{}h\to{}k$ are horizontally composable 2-cells in $\B$ with composite $\gamma$, then for any lifts $p$, $q$ of $f$, $h$, we have $F(qp,\gamma)=F(q,\beta)\cdot{}F(p,\alpha)$.
\item If $\alpha\colon{}f\to{}g$ and $\beta\colon{}g\to{}h$ are vertically composable 2-cells in $\B$ with  composite $\gamma$, then for any lift $p$ of $f$, we have $F(F(p,\alpha),\beta)=F(p,\gamma)$.
\item If $\alpha\colon{}f\to{}f$ is an identity 2-cell in $\B$, then for any lift $p$ of $f$, we have $F(\alpha,p)=p$.
\end{enumerate}
\pf
Suppose we are given two extensions $\fib{D}$, $\fib{D'}$ with the same function $F$, and we want to show they are isomorphic. Note that a 2-functor $\D\to\D'$ extending the identity is determined by where it sends 2-cells. If such a 2-functor is to commmute with $\fib{D}$, $\fib{D'}$, it must send a 2-cell $p\to{}q$ in $\D$ lying over $\alpha\colon{}f\to{}g$ to the unique 2-cell in $\D'$ lying over $\alpha$ with domain $p$.

Let us see that the above prescription actually defines a 2-functor. By definition, this prescription preserves the domain of 2-cells, and it also preserves codomains, since $\D$ and $\D'$ have the same function $F$. Next, for 2-cells $\alpha$ and $\beta$ in $\D$ with (horizontal or vertical) composite $\gamma$, we must see that the image of $\gamma$ in $\D'$ is the composite of the images of $\alpha$ and $\beta$. However, these must both be the unique 2-cell of $\D'$ with the appropriate domain and codomain which lies over the image of $\gamma$ in $\B$.

For the ``moreover'' claim, it is clear that the function $F$ coming from a 1D2F satisfies the conditions (i)-(iii). Conversely, given any $F$ satisfying (i)-(iii), we can take the 2-cells of $\D$ to be the set of pairs $(\alpha,p)$ comprising the domain of $F$, where the domain and codomain of $(\alpha,p)$ are $p$ and $F(\alpha,p)$, respectively, and we can take the extension $\fib{D}$ to send $(\alpha,p)$ to $\alpha$.

The requirement that $\fib{D}$ be a 2-functor forces the definition of composition in $\D$; for example, given 2-cells $(\alpha,p)$ and $(\beta,q)$ of $\D$ with $q=F(\alpha,p)$, their composite is forced to be the 2-cell $(\gamma,p)$ where $\gamma$ is the composite of $\alpha$ and $\beta$. That this composite has the appropriate codomain is ensured by the condition (ii).

Finally, this prescription defines a 2-category, since for each equation which is required in the definition of a 2-category, the two sides are automatically equal, as there is a unique 2-cell with the appropriate domain and lying over the appropriate 2-cell in $\B$. For example, if $(\alpha,p)$, $(\beta,q)$, and $(\gamma,r)$ are vertically composable 2-cells, then in the associativity equation, both sides must be equal to the unique 2-cell with domain $p$ lying over the composite of $\alpha$, $\beta$, $\gamma$.

The existence of identity 2-cells is guaranteed by the condition (iii).
\qed

\subsubsection{}\label{defn:natb}
\defn
For the rest of \S\ref{subsec:1d2fs}, fix a $\wedgeq$-cloven $\wedgeq$-fibration $\fibr{C}CB$, where $\B$ is considered to have the 2-categorical structure given by Theorem~\ref{thm:2cat}.

Given an object $B\in\Ob\B$ and an object $P\in\Ob\fib{C}^B$, we define $\nat^B_P\colon\pi_1^*P\wedge\Eq_B\to\pi_2^*P$ to be the unique morphism over $B\times{}B$ making the diagram
\[
  \begin{tikzcd}
    &\pi_1^*P\wedge\Eq_B\ar[d, dashed, "\nat_B^P"]\\
    P\rac[ru, "\brr{\cind{\id_P},\rho_B!}"]\ar[r, "\cind{\id_P}"', pos=0.6]&\pi_2^*P\\[-10pt]
    B\ar[r, "\Delta_B"]&B\times{}B
  \end{tikzcd}
\]
commute, where $\brr{\cind{\id_P},\rho_B!}$ is cocartesian by Frobenius Reciprocity.

\subsubsection{}\label{defn:alpha-related}
\defn
Given morphisms $f,g\colon{}A\to{}B$ in $\B$ and a $\fib{C}$-homotopy $f\tox{\alpha}g$, as well as lifts $p\colon{}P\to{}Q$ of $f$ and $q\colon{}P\to{}Q$ of $g$, we say that $p$ and $q$ are \emph{$\alpha$-related} if $q$ is the unique morphism over $g$ such that the following diagram commutes (namely $q=\ct\cdot\nat^Q_B\cdot\brr{\bar{p},\alpha!}$).
\begin{equation}\label{eq:alpha-related-triangle}
  \begin{tikzcd}[column sep=60pt]
    &\pi_1^*Q\wedge\Eq_B\ar[d, "\nat_B^Q"]\\
    P\ar[r, "\cind{q}"]\ar[ru, "\brr{\cind{p},\alpha!}"]&\pi_2^*Q\\[-10pt]
    A\ar[r, "\br{f,g}"]&B\times{}B
  \end{tikzcd}
\end{equation}

\subsubsection{}\label{prop:1d2f-hcmp}
\prop
Given morphisms and $\fib{C}$-homotopies as in
\[
  \begin{tikzcd}[row sep=0pt, column sep=10pt]
    &&
    &&\\[-5pt]
    P\ar[rr, "p", bend left]\ar[rr, "q"', bend right]&&
    Q\ar[rr, "r", bend left]\ar[rr, "s"', bend right]&&R\\
    &{}&&{}&\\[13pt]
    &\ar[dd, shorten >=2pt, shorten <=-2pt, Rightarrow, "\ \alpha", pos=0.4]&
    &\ar[dd, shorten >=2pt, shorten <=-2pt, Rightarrow, "\ \beta", pos=0.4]&\\[-5pt]
    A\ar[rr, "f", bend left]\ar[rr, "g"', bend right]&&
    B\ar[rr, "h", bend left]\ar[rr, "k"', bend right]&&C,\\
    &{}&&{}&
  \end{tikzcd}
\]
if $p$ and $q$ are $\alpha$-related, and $r$ and $s$ are $\beta$-related, then $rp$ and $sq$ are ($\beta\hhcmp\alpha$)-related.
\pf
Assume the hypothesis. By the definition of $\beta\hhcmp\alpha$ and of ``($\beta\hhcmp\alpha$)-related'', we need to show that the outside of the following diagram commutes.
\[
\begin{tikzcd}
  &\pi_1^*Q\wedge\Eq_B\ar[r, "\cind{r\ct}\wwdge\check\beta"]\ar[d, "\nat_Q"]
  &\pi_1^*R\wedge\Eq_C\ar[d, "\nat_R"]\\
  P\ar[ru, "\brr{\cind{p},\alpha!}"]\ar[r, "\cind{q}"]&\pi_2^*Q\ar[r, "\cind{s\ct}"]&\pi_2^*R\\[-10pt]
  A\ar[r, "\br{f,g}"]&B\times{}B\ar[r, "h\times{}k"]&C\times{}C
\end{tikzcd}
\]
Now, the triangle commutes by the assumption that $p$ and $q$ are $\alpha$-related. That the square commutes follows from a diagram chase in
\[
\begin{tikzcd}
  &[40pt]\pi_1^*Q\wedge\Eq_B\ar[r, "\cind{r\ct}\wwdge\check\beta"]\ar[d, "\nat_Q"]
  &\pi_1^*R\wedge\Eq_C\ar[d, "\nat_R"]\\
  Q\rac[ru, "\brr{\cind{\id_Q},\rho_B!}"{pos=0.8, xshift=10pt}]\ar[r, "\cind{\id_Q}"' near end]\ar[rr, "\cind{s}"', bend right=20pt]
  \ar[rru, "\brr{r,\beta}", bend left=45pt]
  &\pi_2^*Q\ar[r, "\cind{s\ct}"]&\pi_2^*R\\[3pt]
  B\ar[r, "\Delta_B"]&B\times{}B\ar[r, "h\times{}k"]&C\times{}C
\end{tikzcd}
\]
using the definitions of $\nat_Q$ and $\check\beta$, and the fact that $r$ and $s$ are $\beta$-related.
\qed

\subsubsection{}\defn
Given objects $B\in\Ob\B$ and $P\in\Ob\fib{C}^B$ and natural numbers $1\le{}i,j\le{}n$ we define $\nat_B^{P,ij}$ to be the unique morphism $\Eq_B^{ij}\wedge\pi_i^*P\to\pi_j^*P$ in $\fib{C}^{B^n}$ making the following diagram commute (namely, $\nat_B^{P,ij}=\cind{\ct\ccmp\nat_B^P\ccmp(\ct\wwdge\,\ccind)}$).
\[
  \begin{tikzcd}
    \pi_i^*P\wedge\Eq_B^{ij}\ar[r, "\ccind\wwdge\ct"]\ar[d, "\nat_B^{P,ij}"', dashed]&
    \pi_1^*P\wedge\Eq_B\ar[d, "\nat_B^P"]\\
    \pi_j^*P\car[r, "\ccind"]&\pi_2^*P\\[-15pt]
    B^n\ar[r, "\br{\pi_i,\pi_j}"]&B\times{}B
  \end{tikzcd}
\]

\subsubsection{}\label{lem:natb-gen}
\lem
Given objects $B\in\Ob\B$ and $P\in\Ob\fib{C}^B$ and natural numbers $1\le{}i,j\le{}n$, the following diagram commutes.
\[
  \begin{tikzcd}[column sep=40pt]
    &\pi_i^*P\wedge\Eq_B^{ij}\ar[d, "\nat_B^{P,ij}"]\\
    P\ar[ru, "\brr{\cind{\id_P},\rho_B^{ij}!}"]\ar[r, "\cind{\id_P}"']&
    \pi_j^*P\\[-10pt]
    B\ar[r, "\Delta_B^n"]&B^n
  \end{tikzcd}
\]
\pf
This follows from a diagram chase in the following diagram.
\[
  \begin{tikzcd}[column sep=40pt, baseline=(QED.base)]
    P\areq[d]\ar[r, "\brr{\cind\id,\rho_B^{ij}!}"']
    \ar[rr, "\brr{\cind\id,\rho_B!}", bend left=15pt]&
    \pi_i^*P\wedge\Eq_B^{ij}\ar[d, "\nat_B^{P,ij}"]\ar[r, "\ccind\wwdge\ct"']&
    \pi_1^*P\wedge\Eq_B\ar[d, "\nat_B^P"]\\
    P\ar[r, "\cind\id"]\ar[rr, "\cind\id"', bend right=15pt]&\pi_j^*P\car[r, "\ccind"]&\pi_2^*P\\[5pt]
    A\ar[r, "\br{f_1,\ldots,f_n}"]&B^n\ar[r,"\br{\pi_i,\pi_j}"]& |[alias=QED]| B\times{}B
  \end{tikzcd}
  \tag*{\qed}
\]

\subsubsection{}\label{lem:alpha-related-gen}\lem
With $p$, $q$, and $\alpha$ as in Definition~\ref{defn:alpha-related}, if $p$ and $q$ are $\alpha$-related, then the following diagram commutes.
\[
  \begin{tikzcd}[column sep=60pt]
    &\pi_i^*Q\wedge\Eq_B^{ij}\ar[d, "\nat_B^{Q,ij}"]\\
    P\ar[r, "\cind{q}"]\ar[ru, "\brr{\cind{p},\alpha!}"]&\pi_j^*Q\\[-10pt]
    A\ar[r, "\br{f_1,\ldots,f_n}"]&B^n
  \end{tikzcd}
\]
\pf
This follows from a diagram chase in the following diagram.
\[
  \begin{tikzcd}[baseline=(QED.base)]
    &[40pt]\pi_i^*Q\wedge\Eq_B^{ij}\ar[d, "\nat_B^{Q,ij}"]\ar[r, "\ccind\wwdge\ct"]&
    \pi_1^*Q\wedge\Eq_B\ar[d, "\nat_B^Q"]\\
    P\ar[r, "\cind{q}"]\ar[ru, "\brr{\cind{p},\cind{\alpha}!}"]&\pi_j^*Q\car[r, "\ct"]&Q\\[-10pt]
    A\ar[r, "\br{f_1,\ldots,f_n}"]&B^n\ar[r, "\br{\pi_i,\pi_j}"]&|[alias=QED]| B\times{}B
  \end{tikzcd}
  \tag*{\qed}
\]

\subsubsection{}\label{prop:1d2f-vcmp}
\prop
Given morphisms and $\fib{C}$-homotopies as in
\[
  \begin{tikzcd}[column sep=10pt]
    P\ar[rr, "p", bend left=55pt]\ar[rr, "r", bend right=55pt, pos=0.43]\ar[rr, "q"]&{}
    &
    Q\\
    &\ar[d, shorten >=3pt, shorten <=10pt, Rightarrow, "\ \alpha", pos=0.6]&\\[-5pt]
    A\ar[rr, "f", bend left=55pt]\ar[rr, "h"', bend right=55pt]\ar[rr, "g", pos=0.2]&{}
    \ar[d, shorten >=14pt, shorten <=-2pt, Rightarrow, "\ \beta", pos=0.1]&
    B,\\
    &{}&&
  \end{tikzcd}
\]

if $p$ and $q$ are $\alpha$-related and $q$ and $r$ are $\beta$-related, then $p$ and $r$ are $\beta\hvcmp\alpha$-related.
\pf
Assume the hypothesis. By the definition of $\beta\hvcmp\alpha$, we need to show that the outside of the following diagram commutes.
\[
  \begin{tikzcd}
    &\pi_1^*Q\wedge(\Eq_B^{12}\wedge\Eq_B^{23})\ar[d, "\id\wedge\tr_B"]\\
    &\pi_1^*Q\wedge\Eq_B^{13}\ar[r, "\ccind\wwdge\ct"]\ar[d, "\nat_Q^{13}"]&
    \pi_1^*Q\wedge\Eq_B\ar[d, "\nat_Q"]\\
    P\ar[ruu, "\brr{\cind{p},\brr{\cind\alpha!,\cind\beta!}}", end anchor={[xshift=-20pt]}]
    \ar[rr, "\cind{r}"', bend right=13pt]\ar[r, "\cind{r}"]&
    \pi_3^*Q\ar[r, "\ccind"]&\pi_2^*Q\\
    A\ar[r, "\br{f,g,h}"]&B^3\ar[r, "\br{\pi_1,\pi_3}"]&B^2
  \end{tikzcd}
\]
The square and the triangle on the bottom commute, so it remains to see that the triangle on the left commutes, which is the same as the outside of the following diagram commuting.
\[
  \begin{tikzcd}
    &(\pi_1^*Q\wedge\Eq_B^{12})\wedge\Eq_B^{23}
    \ar[d, "\nat_B^{12}\wedge\id"]\ar[r, "\mathrm{assoc}", "\sim"']
    &\pi_1^*Q\wedge(\Eq_B^{12}\wedge\Eq_B^{23})\ar[d, "\id\wedge\tr_B"]\\
    &\pi_2^*Q\wedge\Eq_B^{23}\ar[rd, "\nat_Q^{23}" near start]&
    \pi_1^*Q\wedge\Eq_B^{13}\ar[d, "\nat_Q^{13}"]\\
    P\ar[ruu, "\brr{\brr{\cind{p},\cind\alpha!},\cind\beta!}", end anchor={[xshift=-20pt]}]
    \ar[rr, "\cind{r}"']\ar[ru, "\brr{\cind{q},\cind{\beta}!}"' {yshift=4pt}]&&\pi_2^*Q\\
  \end{tikzcd}
\]
Here, the two triangles commute by Lemma~\ref{lem:alpha-related-gen}, hence it remains to see that the trapezoid commutes.

Now, an application of Frobenius reciprocity to the (by Proposition~\ref{prop:tripeq}) cocartesian morphism $\brr{\rho_B^{12},\rho_B^{23}}$ shows that $\brr{\cind{\id_Q},\brr{\rho_B^{12},\rho_B^{23}}}\colon{}Q\to{}\pi_1^*Q\wedge(\Eq_B^{12}\wedge\Eq_B^{23})$ is cocartesian. Hence, it suffices to see that the above trapezoid commutes after precomposing with this cocartesian morphism. But by Lemma~\ref{lem:natb-gen} and the definition of $\tr_B$, both resulting composites are $\cind{\id_Q}\colon{}Q\to\pi_2^*Q$.
\qed

\subsubsection{}\label{prop:1d2f-ids}
\prop
Every morphism $p$ in $\C$ over a morphism $f$ in $\B$ is $\hid_f$-related to itself.
\pf
This follows from a diagram chase in
\[
  \begin{tikzcd}[column sep=50pt, baseline=(QED.base)]
    &&\pi_1^*Q\wedge\Eq_B\ar[d, "\nat_B^Q"]\\
    P\ar[rr, "\cind{p}"', bend right=15pt]
    \ar[rru, "\brr{\cind{p},\hid_f!}", bend left=25pt]\ar[r, "p"]&
    Q\ar[ru, "\brr{\cind{\id_Q},\rho_B!}"]\ar[r, "\cind{\id_Q}", pos=0.6]&\pi_2^*P\\
    A\ar[r, "f"]&B\ar[r, "\Delta_B"]& |[alias=QED]| B\times{}B
  \end{tikzcd}
  \tag*{\qed}
\]

\subsubsection{}\label{thm:2cat-psf}\thm
There is up to isomorphism a unique extension of the fibration $\fib{C}\colon\C\to\B$ to a 1D2F, where $\B$ is considered with its 2-category structure, such that, for given 1-cells $p$, $q$ in $\C$ lying over $f$, $g$ in $\B$ and a 2-cell $\alpha\colon{}f\to{}g$, there exists a lift $p\to{}q$ of $\alpha$ if and only if $p$ and $q$ are $\alpha$-related.
\pf
This follows immediately from Propositions~\ref{prop:1d2f-hcmp},~\ref{prop:1d2f-vcmp},~and~\ref{prop:1d2f-ids}, using Lemma~\ref{lem:constructing-1d2fs}.
\qed

\subsubsection{}
\cor
$\fib{C}$-homotopic morphisms induce isomorphic pullback functors.
\pf
We will only indicate the proof. The claim follows from an application of the (inverse of the) ``Grothendieck construction'' referred to above. In general, given any cloven 1D2F, each 2-cell $\alpha\colon{}f\to{}g$ induces a natural transformation $f^*\to{}g^*$, which is an isomorphism if $\alpha$ is. And by Theorem~\ref{thm:groupoid}, every $\fib{C}$-homotopy is an isomorphism 2-cell.
\qed

\subsection{2-categorical products}\label{subsec:c-htpy-finprod}
In this section, we will show that the 2-categorical structure on the base category of a $\wedgeq$-fibration has finite products in the 2-categorical sense. We will begin by recalling what this means.

The main effort in this section will be devoted to showing that, given objects $A,B$ in the base of a $\wedgeq$-fibration $\fib{C}$, a certain morphism $\Eq_{A\times{}B}\to\pi_1^*\Eq_A\wedge\pi_2^*\Eq_B$ in $\fib{C}^{A\times{}B}$ is an isomorphism. That these objects are isomorphic at all (which, logically speaking, says $\br{a_1,b_1}=\br{a_2,b_2}\Leftrightarrow{}a_1=a_2\wedge{}b_1=b_2$) is already proven in \cite[p.~10]{lawvereadj} (see also \cite[Exercise~3.4.7]{jacobscatlogic}).

\subsubsection{}\defn
Let $\C$ be a 2-category.

Given a pair of objects $A,B\in\Ob\C$, a \emph{product diagram based on $A$ and $B$} consists of an object $C\in\Ob\C$ and a pair of morphisms $A\xleftarrow{f}C\tox{g}B$ having the following universal property: for any object $D$, the functor $\br{f\circ\mathord{\text{--}},g\circ\mathord{\text{--}}}\colon\HOM_\C(D,C)\to\HOM_\C(D,A)\times\HOM_\C(D,B)$, induced by composition with $f$ and $g$, is an isomorphism of categories.

An object $A\in\Ob\C$ is a \emph{terminal} object if, for each object $X\in\Ob\C$, the category $\Hom_\C(X,A)$ has a single object and a single morphism.

$\C$ \emph{has finite products} if it has a terminal object and there is a product diagram based on each pair of objects.

We note that there are other (weaker) notions of 2-categorical products, but this is the only one we use. Note also that a product diagram in a 2-category is also a product diagram in the underlying category, and similarly for the terminal object.

\subsubsection{}\label{thm:2cat-prods}\thm
For any $\wedgeq$-cloven $\wedgeq$-fibration $\fibr{C}CB$, the 2-categorical structure on $\B$ given in Theorem~\ref{thm:2cat} has finite products.
\pf
Let $A,B\in\Ob\B$. We already know that, for any $C\in\Ob\B$, composition with $\pi_1\colon{}A\times{}B\to{}A$ and $\pi_2\colon{}A\times{}B\to{}B$ induces a bijection $\Hom_\B(C,A\times{B})\to\Hom_\B(C,A)\times\Hom_\B(C,B)$. What we still need to show is that, given morphisms
\[
  \begin{tikzcd}
    &C
    \ar[dl, "f"', bend right=15pt]\ar[dl, "g", bend left=15pt]
    \ar[dr, "h"', bend right=15pt]\ar[dr, "k", bend left=15pt]&\\
    A&&B,
  \end{tikzcd}
\]
composing horizontally with $\hid_{\pi_1}$ and $\hid_{\pi_2}$ induces a bijection
\[
  \HOM_{\B}(C,A\times{B})(\br{f,h},\br{g,k})\to
  \HOM_{\B}(C,A)(f,g)\times
  \HOM_{\B}(C,B)(h,k).
\]
Now, given a $\fib{C}$-homotopy $\top_C\tox{\alpha}\Eq_{A\times{}B}$ from $\br{f,h}$ to $\br{g,k}$, its image under the above morphism is given by composing with $\widecheck{\hid_{\pi_1}}\colon\Eq_{A\times{}B}\to\Eq_A$ and $\widecheck{\hid_{\pi_2}}\colon\Eq_{A\times{}B}\to\Eq_B$. Hence, it suffices to see that
\begin{equation}\label{eq:switcheroo-map}
  \begin{tikzcd}[column sep=70pt]
    \Eq_{A\times{}B}
    \ar[r, "\brr{\cind{\widecheck{\hid_{\pi_1}}},\cind{\widecheck{\hid_{\pi_2}}}}"]&
    \pi_1^*\Eq_A\wedge\pi_2^*\Eq_B\\[-10pt]
    (A\times{}B)\times(A\times{}B)\ar[r, "\br{\pi_1\times\pi_1,\pi_2\times\pi_2}"]&
    (A\times{}A)\times(B\times{}B)
  \end{tikzcd}
\end{equation}
is cartesian, since this would give us bijections
\[
  \begin{tikzcd}
    \Hom_{\br{\br{f,h},\br{g,k}}}(\top_C,\Eq_{A\times{}B})
    \ar[r, "\brr{\cind{\widecheck{\hid_{\pi_1}}},\cind{\widecheck{\hid_{\pi_2}}}}\circ\mathord{\text{--}}", "\sim"']
    \ar[ddr, bend right, pos=0.3, end anchor={[xshift=-50pt,yshift=10pt]},
        "\br{\widecheck{\hid_{\pi_1}}\circ\mathord{\text{--}},\widecheck{\hid_{\pi_2}}\circ\mathord{\text{--}}}"']
    &\Hom_{\br{\br{f,g},\br{h,k}}}(\top_C,\pi_1^*\Eq_{A}\wedge\pi_2^*\Eq_{B})
    \ar[d, "\br{\pi_1\circ\mathord{\text{--}},\pi_2\circ\mathord{\text{--}}}", "\vsim"']\\
    &\Hom_{\br{\br{f,g},\br{h,k}}}(\top_C,\pi_1^*\Eq_{A})\times
    \Hom_{\br{\br{f,g},\br{h,k}}}(\top_C,\pi_2^*\Eq_{B})
    \ar[d, "(\ct\circ\mathord{\text{--}})\times(\ct\circ\mathord{\text{--}})", "\vsim"']\\
    &\Hom_{\br{f,g}}(\top_C,\Eq_{A})\times
    \Hom_{\br{h,k}}(\top_C,\Eq_{B}).
  \end{tikzcd}
\]

Since the morphism \eqref{eq:switcheroo-map} lies over an isomorphism, it is an isomorphism if and only if it is cartesian, and also if and only if it is cocartesian. Let us see that it is cocartesian.

Now in general, given composable morphisms $p$ and $q$, if $p$ and $q\cdot{}p$ are both cocartesian, then so is $q$. Hence, it suffices to see that the composite
\[
  \begin{tikzcd}[column sep=70pt]
    \top_{A\times{}B}\rac[r, "\rho_{A\times{}B}"]
    \ar[rr, "\brr{\cind{\hid_{\pi_1}},\cind{\hid_{\pi_2}}}"', bend right=10pt]&
    \Eq_{A\times{}B}
    \ar[r, "\brr{\cind{\widecheck{\hid_{\pi_1}}},\cind{\widecheck{\hid_{\pi_2}}}}"]&
    \pi_1^*\Eq_A\wedge\pi_2^*\Eq_B\\[10pt]
    A\times{}B\ar[r, "\Delta_B"]
    &(A\times{}B)\times(A\times{}B)\ar[r, "\br{\pi_1\times\pi_1,\pi_2\times\pi_2}"]&
    (A\times{}A)\times(B\times{}B)
  \end{tikzcd}
\]
is cocartesian. We will show this by a similar argument to that used in Proposition~\ref{prop:tripeq}.

Namely, we will show that each morphism in the following factorization of
$\brr{\cind{\hid_{\pi_1}},\cind{\hid_{\pi_2}}}$ is cocartesian.
\begin{equation}
  \label{eq:prodeqmors}
  \begin{tikzcd}[column sep=70pt]
    \top_{A\times{}B}\ar[r, "\cind{\hid_{\pi_1}}"]&
    \pi_1^*\Eq_{A}
    \ar[r, "\brr{\ccind,\cind{\hid_{\pi_2}}!}"]&
    \pi_1^*\Eq_A\wedge\pi_2^*\Eq_B\\[-10pt]
    A\times{}B\ar[r, "\Delta_A\times{}\id_B"]
    &(A\times{}A)\times{}B\ar[r, "\id_{A\times{}A}\times\Delta_B"]&
    (A\times{}A)\times(B\times{}B)
  \end{tikzcd}
\end{equation}

The first of the morphisms in \eqref{eq:prodeqmors} is cocartesian by the stability of $\rho_A$ along the product projection $\pi_1\colon(A\times{}A)\times{}B\to{}A\times{}A$:
\[
  \begin{tikzcd}
    \top_{A\times{}B}\car[d, "\exx_{\pi_1}"']\ar[r, "\cind{\hid_{\pi_1}}"]&\pi_1^*\Eq_A\car[d, "\ct"]\\
    \top_{A}\rac[r, "\rho_A"]&\Eq_A
  \end{tikzcd}
  \quad\quad\quad
  \begin{tikzcd}
    A\times{}B\ar[r, "\Delta_A\times\id_B"]\ar[d ,"\pi_1"']
    \ar[rd, phantom, "\lrcorner", pos=0]
    &[15pt]
    (A\times{}A)\times{}B\ar[d, "\pi_1"]\\
    A\ar[r, "\Delta_A"]&A\times{}A.
  \end{tikzcd}
\]

Similarly, we have that $\cind{\hid_{\pi_2}}\colon\top_{(A\times{}A)\times{}B}\to\pi_2^*\Eq_B$ over $\id_{A\times{}A}\times\Delta_B$ is cocartesian, from which it follows by Frobenius reciprocity that the second morphism of \eqref{eq:prodeqmors} is cocartesian.

We have shown that $\B$ has 2-categorical binary products. It remains to check that the terminal object $\tm_\B$ is a 2-categorical terminal object, i.e., that for any $C\in\Ob\B$, there is a unique $\fib{C}$-homotopy $!_C\to!_C$. This is the case since $\Eq_{\tm_\B}$ is terminal in $\fib{C}^{\tm_\B\times\tm_\B}$ (to see this, note that $\exx_{\Delta_{\tm_\B}}\colon\top_{\tm_\B}\to\top_{\tm_\B\times\tm_\B}$ is (an isomorphism and hence) cocartesian).
\qed

\subsection{Universal property}\label{subsec:c-htpy-canon}
In this section, we show that, up to isomorphism, the extension of a $\wedgeq$-cloven $\wedgeq$-fibration to a 1D2F given by Theorem~\ref{thm:2cat-psf} does not depend on the chosen $\wedgeq$-cleavage.

Intuitively, this is rather clear: given any two $\wedgeq$-cleavages, one obtains canonical isomorphisms between the chosen objects of each $\wedgeq$-cleavage, and one could then try and check that under these isomorphisms, all of the constructions involved in the definition of the 1D2F-structure correspond.

However, as was pointed out to me by Arpon Raksit, a more satisfying solution to this problem would be to find a ``universal property'' that, given $\fib{C}$, characterizes the resulting 1D2F up to isomorphism. Such a universal property was suggested by the anonymous referee, and we carry out the details below.

The main point here is that a $\wedgeq$-fibration $\fib{C}$ is bound up with its induced 2-categorical structure in the following way: given \emph{any} 1D2F structure on $\fib{C}$, one can canonically associate to each of its 2-cells a $\fib{C}$-homotopy. Following the referee's suggestion, we call the original 1D2F \emph{univalent} if this assignment induces an isomorphism of 1D2Fs (importantly, this property can be stated without reference to a fixed cleavage -- see Definition~\ref{defn:univalent}). Then the 1D2F structure on $\fib{C}$ from Theorem~\ref{thm:2cat-psf} is the (up to isomorphism) unique univalent 1D2F structure on $\fib{C}$ (see Theorem~\ref{thm:univalent-main-thm} below).

In fact (still following the referee), we show something stronger (Theorem~\ref{thm:univalent-stronger-version}): that the canonical morphism from a $\wedgeq$-1D2F $\fib{C}$ to the 1D2F structure coming from Theorem~\ref{thm:2cat-psf} is universal among morphisms to univalent 1D2Fs. We note that this is similar to the universal property of the ``extensional collapse'' of a posetal fibration defined in \cite[p.~214]{jacobscatlogic} and \cite[Section~5]{maietti-rosolini-elementary-completion}, in which the notion of equality in the base category is forced to agree with that given by the fibration.

This leads to a second universal property of the 1D2F structure on a $\wedgeq$-fibration (Corollary~\ref{cor:univalent-other-universal-property}) which does not make any reference to univalent $\wedgeq$-1D2Fs.

Below, we will be dealing with $\wedgeq$-fibrations which aren't cloven, but we will, following our usual convention (as in \S\S\ref{subsubsec:cleavage},~\ref{subsubsec:fpnotations},~\ref{subsubsec:wedge-fib}) still use the same notation as if we had a cleavage. For example, given $B\in\Ob\B$, once we have chosen products $B\times{}B,B^3\in\Ob\B$, a terminal object $\top_B\in\Ob\fib{C}^B$, a cocartesian morphism $\rho_B\colon\top_B\to\Eq_B$ over $\Delta_B$, pullbacks $\Eq_B^{12},\Eq_B^{23},\Eq_B^{13}\in\Ob\B^3$, and a product $\Eq_B^{12}\wedge\Eq_B^{23}$, we will write $\tr_B$ for the morphism $\Eq_B^{12}\wedge\Eq_B^{23}\to\Eq_B^{13}$ given by the construction in Definition~\ref{defn:trB} with respect to these choices.

We will often need to fix several such choices, and for brevity, we may simply write ``fix a choice of $\top_A,B\times{}B,\ldots$'' rather than ``fix a terminal object $\top_A\in\Ob\fib{C}^A$, a product $B\xot{\pi_1}B\times{}B\tox{\pi_2}B$, \ldots''.

\subsubsection{}\defn
A 1D2F $\fibr{C}CB$ is a \emph{$\wedgeq$-1D2F} if (i) its underlying fibration is a $\wedgeq$-fibration, (ii) each 2-cell in $\B$ (and hence in $\C$) is invertible, and (iii) the products in $\B$ are products in the 2-categorical sense.

Note that each $\wedgeq$-fibration is a $\wedgeq$-1D2F when considered with the trivial 2-categorical structure and that, by Theorems~\ref{thm:groupoid}~and~\ref{thm:2cat-prods}, the extension of a $\wedgeq$-cloven $\wedgeq$-fibration to a 1D2F given by Theorem~\ref{thm:2cat-psf} is a $\wedgeq$-1D2F.

Given $\wedgeq$-1D2Fs $\fibr{C}CB$ and $\fibr{C}{C'}{B'}$, a \emph{morphism of $\wedgeq$-1D2Fs} from $\fib{C}$ to $\fib{C'}$ is a pair $(\Phi,\phi)$ of 2-functors $\Phi\colon\C\to\C'$ and $\phi\colon\B\to\B'$ such that (i) $\fib{C'}\Phi=\phi\fib{C}$, (ii) $\phi$ preserves products, and (iii) $\Phi$ preserves cartesian morphisms, cocartesian lifts of diagonal morphisms, and fiberwise products.

$\wedgeq$-1D2Fs and morphisms thereof form a category in an obvious manner.

\subsubsection{}\label{defn:univalent}\defn
Let $\fib{C}$ be a $\wedgeq$-1D2F. Given morphisms $f,g,h\colon{}A\to{}B$, a 2-cell $g\tox{\beta}h$ in $\B$, a choice of terminal objects $\top_A\in\Ob\fib{C}^A$ and $\top_B\in\Ob\fib{C}^B$, a product $B\xot{\pi_1}B\times{}B\tox{\pi_2}B$, and an equality object $\rho_B\colon\top_B\to\Eq_B$, as well as a morphism $p\colon\top_A\to\Eq_B$ lying over $\br{f,g}$, we define $\beta_!p\colon\top_A\to\Eq_B$ to be the unique morphism for which there exists a 2-cell $p\to\beta_!p$ lying over $\br{\id_f,\beta}$:
\[
  \begin{tikzcd}[column sep=50pt, row sep=30pt]
    \top_A
    \ar[r, bend left=20pt, "p"{name=p}]
    \ar[r, bend right=20pt, "\beta_!p"'{name=q}, dashed]
    &\Eq_B
    \ar[Rightarrow, from=p, to=q, dashed, shorten=3pt]
    \\
    A
    \ar[r, bend left=20pt, "\br{f,g}"{name=f}, pos=0.55]
    \ar[r, bend right=20pt, "\br{f,h}"'{name=g}, pos=0.55]
    &B\times{}B.
    \ar[Rightarrow, from=f, to=g, "\br{\id_f,\beta}", shorten=3pt]
  \end{tikzcd}
\]
The most important case is when $f=g$ and $p=\hid_f$, and in this case we write $\lft{\beta}$ for $\beta_!\hid_f$.

$\fib{C}$ is \emph{univalent} if for each pair $f,g\colon{}A\to{}B$ of morphisms in $\B$, and for any (and hence every) choice of $\top_A,B\times{}B,\top_B,\Eq_B$, the operation $\beta\mapsto\lft{\beta}$ establishes a bijection
\begin{equation}\label{eq:univalent-morphism}
  \Hom_{\HOM_{\B(A,B)}}(f,g)\tox\sim\Hom^{\fib{C}}_{\br{f,g}}(\top_A,\Eq_B).
\end{equation}

Let $\fibr{C}CB$ be a $\wedgeq$-1D2F until further notice.

\subsubsection{}\label{prop:construction-is-univalent}\prop
If $\fib{C}$ is the 1D2F induced by Theorem~\ref{thm:2cat-psf} from a $\wedgeq$-cleavage of the underlying $\wedgeq$-fibration $\abs{\fib{C}}$, then $\fib{C}$ is univalent.
\pf
The $\wedgeq$-cleavage of $\abs{\fib{C}}$ provides us, for each $A,B\in\Ob\fib{C}$, with the choices necessary to carry out the construction $\alpha\mapsto\lft{\alpha}$ from Definition~\ref{defn:univalent}, which in this case gives us, for each 2-cell $\alpha\colon{}f\to{}g$, a $\fib{C}$-homotopy $\top_A\to\Eq_B$ from $f$ to $g$. But $\alpha$ is by the definition of the 2-categorical structure on $\B$ also a $\fib{C}$-homotopy, hence the map (\ref{eq:univalent-morphism}) is in this case a map
\[
  \Hom_{\br{f,g}}(\top_A,\Eq_B)\to\Hom_{\br{f,g}}(\top_A,\Eq_B).
\]
To see that this is a bijection, let us show that it is in fact the identity.

Inspecting the definition of the above map, and of the 2-categorical structure on $\C$, we see that this amounts to showing that for each $\fib{C}$-homotopy $\alpha\colon\top_A\to\Eq_B$ from $f$ to $g$, the two morphisms $\hid_f,\alpha\colon\top_A\to\Eq_B$ are $\br{\hid_f,\alpha}$-related, i.e., that the triangle in the following diagram commutes.
\[
  \begin{tikzcd}
    &\pi_1^*\Eq_B\wedge\Eq_{B\times{}B}\ar[d, "\nat_{B\times{}B}^{\Eq_B}"]
    \ar[r, "\sim"',
    "\brr{\cind{\widecheck{\hid_{\pi_1}}}\pi_2,\ccind\pi_1,\cind{\widecheck{\hid_{\pi_2}}}\pi_2}"]
    &[70pt]
    \Eq_B^{12}\wedge\Eq_B^{23}\wedge\Eq_B^{34}
    \ar[d, dashed]\\
    \top_A\ar[ru, "\brr{\cind{\hid_f},\br{\hid_f,\alpha}}"]\ar[r, "\cind\alpha"]&
    \pi_2^*\Eq_B\ar[r, "\ccind", "\sim"']
    &\Eq_B^{14}\\
    A\ar[r, "\br{\br{f,f},\br{f,g}}"]&
    (B\times{}B)^2\ar[r, "\br{\pi_1\pi_2,\pi_1\pi_1,\pi_2\pi_1,\pi_2\pi_2}", "\sim"']&
    B^4
  \end{tikzcd}
\]
Here, the dashed morphism is the unique one making the rectangle commute. That the topmost horizontal arrow is an isomorphism follows from the fact that the morphism (\ref{eq:switcheroo-map}) from Theorem~\ref{thm:2cat-prods} is an isomorphism, as was established there. Hence, it suffices to see that the outside of the above diagram, displayed below, commutes.
\begin{equation}\label{eq:construction-univalent-proof-triangle}
  \begin{tikzcd}
    &\Eq_B^{12}\wedge\Eq_B^{23}\wedge\Eq_B^{34}\ar[d, dashed]\\
    \top_A\ar[ru, "\brr{\cind\hid_f,\cind\hid_f,\cind\alpha}"]
    \ar[r, "\cind\alpha"]&\Eq_B^{14}\\[-10pt]
    A\ar[r, "\br{f,f,f,g}"]&B^4.
  \end{tikzcd}
\end{equation}

Let us now identify the dashed morphism. The outside of the commutative diagram
\[
  \begin{tikzcd}
    &&\pi_1^*\Eq_B\wedge\Eq_{B\times{}B}\ar[d, "\nat_{B\times{}B}^{\Eq_B}"]
    \ar[r, "\sim"',
    "\brr{\cind{\widecheck{\hid_{\pi_1}}}\pi_2,\ccind\pi_1,\cind{\widecheck{\hid_{\pi_2}}}\pi_2}"]
    &[70pt]
    \Eq_B^{12}\wedge\Eq_B^{23}\wedge\Eq_B^{34}
    \ar[d, dashed]\\
    \top_B\rac[r, "\rho_B"]&
    \Eq_B\rac[ru, "\brr{\cind{\id_{\Eq_B}},\rho_{B\times{}B}!}"]\ar[r, "\cind{\id_{\Eq_B}}"]&
    \pi_2^*\Eq_B\ar[r, "\ccind", "\sim"']
    &\Eq_B^{14}\\[-5pt]
    B\ar[r, "\Delta"]&B\times{}B\ar[r, "\Delta_{B\times{}B}"]&
    (B\times{}B)^2\ar[r, "\br{\pi_1\pi_2,\pi_1\pi_1,\pi_2\pi_1,\pi_2\pi_2}", "\sim"']&
    B^4
  \end{tikzcd}
\]
gives a diagram
\[
  \begin{tikzcd}
    &\Eq_B^{12}\wedge\Eq_B^{23}\wedge\Eq_B^{34}\ar[d, dashed]\\
    \top_B\rac[ru]
    \ar[r, "\cind{\rho_B^{14}}"']&\Eq_B^{14}\\[-10pt]
    A\ar[r, "\Delta_B^4"]&B^4.
  \end{tikzcd}
\]
in which the displayed morphism is cocartesian. Hence, we see that the dashed morphism is precisely the morphism shown in (\ref{eq:pent}) in the proof of Theorem~\ref{thm:homcat-is-cat}. It follows that the triangle (\ref{eq:construction-univalent-proof-triangle}) commutes if and only if $\alpha$ is equal to $\alpha\hvcmp\hid_f\hvcmp\hid_f$, which of course it is.
\qed

\subsubsection{}\label{lem:total-cat-induced-2cells}\lem
Given objects and morphisms in $\B$ and $\C$ and a 2-cell $\alpha$ as shown below, with $\ct$ cartesian, if there is a 2-cell $\ct{}p\to\ct{}q$ over $h\circ{}\alpha$, then there is a 2-cell $p\to{}q$ over $\alpha$.
\[
  \begin{tikzcd}
    P\ar[rrd, bend left=30pt, shorten <=10pt, "\ct{}p" name=pp]
    \ar[rrd, shorten <=10pt, "\ct{}q"' {name=qq, inner sep=0pt}, bend left=5pt]
    \ar[from=pp, to=qq, Rightarrow, shorten=2pt]
    \ar[rd, "p" {name=p, inner sep=1pt}, bend right=5pt]\ar[rd, bend right=40pt, "q"' name=q]
    \ar[from=p, to=q, Rightarrow, dashed, shorten=2pt]
    &&\\[15pt]
    &g^*Q\car[r, "\ct"']&Q\\
    A\ar[r, bend left, "f" name=f]\ar[r, bend right, "g"' name=g]
    \ar[Rightarrow, from=f, to=g, shorten=3pt, "\ \alpha"]&B\ar[r, "h"]&C
  \end{tikzcd}
\]
The dual statement for cocartesian morphisms also holds.

Similarly, given a product $R\xot{\pi_1}R\wedge{}S\tox{\pi_2}S$ in $\fib{C}^B$ and morphisms $p,q\colon{}P\to{}R\wedge{}S$ over $f$ and $g$, if there are 2-cells $\pi_1p\to\pi_1q$ and $\pi_2p\to\pi_2q$ over $\alpha$, then there is a 2-cell $p\to{}q$ over $\alpha$.
\pf
Let us prove the first claim. The other two are proved similarly.

There is a (unique) 2-cell $\beta\colon{}p\to{}q'$ over $\alpha$ with domain $p$. Composing this with $\ct$ gives a 2-cell $\ct{}p\to\ct{}q'$ over $h\circ\alpha$. Since there is a unique lift of $h\circ\alpha$ with domain $\ct{}p$, it follows that $\ct{}q'=\ct{}q$, and hence $q'=q$ since $\ct$ is cartesian.
\qed

\subsubsection{}\label{lem:univalent-hvcmp-commute}\lem
Let $A,B\in\Ob\C$ and fix a choice of $\top_A$, $\top_B$, $B\times{}B$, $B^3$, $\Eq_B$, $\Eq_B^{12}$, $\Eq_B^{23}$, $\Eq_B^{12}\wedge\Eq_B^{23}$.

Given morphisms $f,g,h,k\colon{}A\to{}B$ and morphisms $p,q\colon\top_A\to\Eq_B$ lying over $\br{f,g}$ and $\br{g,h}$, respectively, as well as a 2-cell $\beta\colon{}h\to k$, we have $\beta_!(q\hvcmp p)=(\beta_!q)\hvcmp p$.
\pf
We must show that there is a 2-cell $q\hvcmp p\to(\beta_!q)\hvcmp p$ over $\br{\id_f,\beta}$.

Since there are 2-cells $p\to p$ and $q\to \beta_!q$ over $\br{\pi_1,\pi_2}\circ\br{\id_f,\id_g,\beta}$ and $\br{\pi_2,\pi_3}\circ\br{\id_f,\id_g,\beta}$, respectively, we have by (the first and last claims in) Lemma~\ref{lem:total-cat-induced-2cells} a 2-cell $\brr{\cind{p},\cind{q}}\to\brr{\cind{p},\cind{\beta_!q}}$ over $\br{\id_f,\id_g,\beta}$.
\[
  \begin{tikzcd}
    \top_A\ar[r, "\brr{\cind{p},\cind{q}}" name=f, pos=0.6, bend left=15pt]
    \ar[r, "\brr{\cind{p},\cind{\beta_!q}}"' name=g, pos=0.6, bend right=15pt]
    \ar[from=f, to=g, Rightarrow, shorten=3pt, dashed]
    &[30pt]\Eq_B^{12}\wedge\Eq_B^{23}
    \ar[r, "\ct\tr_B"]
    &
    \Eq_B
    \\[15pt]
    A\ar[r, "\br{f,g,h}" name=ff, bend left=15pt]
    \ar[r, "\br{f,g,k}"' name=fg, bend right=15pt]
    \ar[from=ff, to=fg, Rightarrow, "\br{\id_f,\id_g,\beta}", shorten=2pt, shift right=5pt]
    &
    B^3\ar[r, "\br{\pi_1,\pi_3}"]&
    B\times{}B
  \end{tikzcd}
\]
Composing this with $\ct\tr_B\colon\Eq_B^{12}\wedge\Eq_B^{23}\to\Eq_B$ then gives the desired 2-cell $q\hvcmp p\to(\beta_!q)\hvcmp p$ over $\br{\id_f,\beta}$.
\qed

\subsubsection{}\label{prop:univalent-hvcmp}\prop
Let $A,B\in\Ob\B$ and fix choices as in Lemma~\ref{lem:univalent-hvcmp-commute}. Given morphisms $f,g,h\colon{}A\to{}B$ in $\B$ and 2-cells $f\tox{\alpha}g\tox{\beta}h$, we have $\lft{\beta\cdot\alpha}=\lft\beta\hvcmp\lft\alpha$.
\pf
We have
\[
  \lft{\beta\cdot\alpha}=_{\mathrm{def}}(\beta\cdot\alpha)_!\hid_f=\beta_!(\alpha_!\hid_f)=
  \beta_!(\hid_g\hvcmp(\alpha_!\hid_f))=
  (\beta_!\hid_g)\hvcmp(\alpha_!\hid_f)=_{\mathrm{def}}
  \lft{\beta}\hvcmp\lft{\alpha},
\]
where the second equality is immediate from the definitions, the third equality comes from $\hid_g$ being a unit for $\hvcmp$, and the fourth equality is by Lemma~\ref{lem:univalent-hvcmp-commute}.
\qed

\subsubsection{}\label{prop:univalent-hhcmp}\prop
Let $A,B,C\in\Ob\B$ and fix a choice of $\top_A$, $\top_B$, $\top_C$, $B\times{}B$, $C\times{}C$, $\Eq_B$, $\Eq_C$. Given morphisms and 2-cells
\[
  \begin{tikzcd}[row sep=0pt, column sep=10pt]
    &\ar[dd, shorten >=2pt, shorten <=-2pt, Rightarrow, "\ \alpha", pos=0.3]&
    &\ar[dd, shorten >=2pt, shorten <=-2pt, Rightarrow, "\ \beta", pos=0.4]&\\[-5pt]
    A\ar[rr, "f", bend left]\ar[rr, "g"', bend right]&&
    B\ar[rr, "h", bend left]\ar[rr, "k"', bend right]&&C\\
    &{}&&{}&
  \end{tikzcd}
\]
in $\B$, we have $\lft{\beta\circ\alpha}=\lft\beta\hhcmp\lft\alpha$.
\pf
We need to show that there is a 2-cell $\hid_{hf}\to\lft{\beta}\hhcmp\lft{\alpha}$ lying over $\br{\id_{hf},\beta\circ\alpha}$.

Since there is 2-cell $\hid_h\to\lft\beta$ lying over $\br{\id_h,\beta}=(\id_h\times\beta)\circ{}\Delta_B$, we have by Lemma~\ref{lem:total-cat-induced-2cells} a 2-cell $\widecheck{\id_h}\to\widecheck{\lft\beta}$ lying over $\id_h\times\beta\colon{}h\times{}h\to{}h\times{}k$.
\[
  \begin{tikzcd}
    &&[8pt]\Eq_C\\
    \top_B\rac[r, "\rho_B"']
    \ar[rru, bend left=30pt, "\hid_h" name=pp]
    \ar[rru, "\lft{\beta}"' {name=qq, inner sep=0pt}, bend left=10pt]
    \ar[from=pp, to=qq, Rightarrow, shorten=3pt]
    &\Eq_B
    \ar[ru, "\widecheck\hid_h" {name=p, inner sep=1pt}, bend right=5pt]
    \ar[ru, bend right=40pt, "\widecheck{\lft{\beta}}"' name=q]
    \ar[from=p, to=q, Rightarrow, dashed, shorten=3pt]
    &&\\
    B\ar[r, "\Delta_B"]&B\times{}B
    \ar[r, bend left=20pt, "h\times{}h" name=f]\ar[r, bend right=20pt, "h\times{}k"' name=g]
    \ar[Rightarrow, from=f, to=g, shorten=3pt, "\id_h\times\beta", shift right=8pt]&C\times{}C
  \end{tikzcd}
\]
Composing this with the 2-cell $\hid_f\to\lft\alpha$ lying over $\br{\id_f,\alpha}$ then gives the desired 2-cell.
\qed

\subsubsection{}\label{prop:univalent-identities}\prop
Let $A,B\in\Ob\B$ and fix a choice of $\top_A$, $\top_B$, $B\times B$, $\Eq_B$. Then $\lft{\id_f}=\hid_f$ for any morphism $f\colon{}A\to B$.
\pf
Immediate from the definitions.
\qed

\subsubsection{}\label{prop:univalent-2cell-related}\prop
Let $A,B\in\Ob\B$ and $P\in\Ob\fib{C}^A$ and $Q\in\Ob\fib{C}^B$, and fix a choice of $\top_B$, $B\times{}B$, $\Eq_B$, $\pi_1^*Q$, $\pi_2^*Q$, $\pi_1^*Q\wedge\Eq_B$. Given a 2-cell
$
\begin{tikzcd}
  A\ar[r, bend left, "f" name=f]\ar[r, bend right, "g"' name=g]
  \ar[Rightarrow, "\alpha", from=f, to=g, shorten=2pt]&B
\end{tikzcd}
$
in $\B$ and lifts $p,q\colon{}P\to{}Q$ of $f$ and $g$, there exists a lift $p\to{}q$ of $\alpha$ if and only if the triangle (\ref{eq:alpha-related-triangle}) from Definition~\ref{defn:alpha-related}, with $\alpha$ replaced by $\{\alpha\}$, commutes -- i.e., if and only if $p$ and $q$ are $\lft{\alpha}$-related.
\pf
Given that there is a unique lift of $\alpha$ with domain $p$ and a unique morphism to which $p$ is $\lft{\alpha}$-related, it suffices to show that if $p$ and $q$ are $\lft{\alpha}$-related, then there is a lift $p\to{}q$ of $\alpha$.

Since there are 2-cells $\hid_f\to\lft\alpha$ and $\cind{p}\to\cind{p}$ over $\br{\id_f,\alpha}$, we have by Lemma~\ref{lem:total-cat-induced-2cells} a 2-cell $\brr{\cind{p},\id_f!}\to\brr{\cind{p},\lft\alpha!}$ over $\br{\id_f,\alpha}$.
\[
  \begin{tikzcd}
    P\ar[r, "\brr{\cind{p},\hid_f!}" name=f, pos=0.6, bend left=15pt]
    \ar[r, "\brr{\cind{p},\lft\alpha!}"' name=g, pos=0.6, bend right=15pt]
    \ar[from=f, to=g, Rightarrow, shorten=3pt]
    &[30pt]\pi_1^*Q\wedge\Eq_B\ar[r]
    \ar[r, "\ct\nat_B^Q"]
    &[10pt]\Eq_B\\[10pt]
    A\ar[r, "\br{f,f}" name=ff, bend left=15pt]
    \ar[r, "\br{f,g}"' name=fg, bend right=15pt]
    \ar[from=ff, to=fg, Rightarrow, "\br{\id_f,\alpha}", shorten=2pt]
    &
    B\times{}B
    \ar[r, "\pi_2"]&
    B
  \end{tikzcd}
\]
Composing this with $\ct\nat_B^Q\colon\pi_1^*Q\wedge\Eq_B\to{}Q$ gives the desired lift $p\to{}q$ of $\alpha$.
\qed

\subsubsection{}\label{prop:univalent-unique-extension}\prop
Given $\wedgeq$-1D2Fs $\fib{C}$ and $\fib{C'}$, with $\fib{C'}$ univalent, any morphism $\abs{\fib{C}}\to\fib{C'}$ of $\wedgeq$-1D2Fs, where $\abs{\fib{C}}$ is the underlying $\wedgeq$-fibration of $\fib{C}$, extends uniquely to a morphism $\fib{C}\to\fib{C'}$.
\pf
Let $(\Phi,\phi)\colon\abs{\fib{C}}\to\fib{C'}$ be a morphism of $\wedgeq$-1D2Fs that we wish to extend.

Now fix a $\wedgeq$-cleavage of $\abs{\fib{C}}$. In particular, this fixes all the data needed in Definition~\ref{defn:univalent} and the above Propositions.

Let us now introduce the following notation. Given a pair of objects $A,B\in\Ob\B$, we write $\phi{}A\times\phi{}B$ for $\phi(A\times{}B)$, and similarly with all the other choices made in the $\wedgeq$-cleavage of $\abs{\fib{C}}$. For example we write $\Eq_{\phi{}B}$ for $\Phi\Eq_B$ and $\rho_{\phi{}B}$ for $\Phi\rho_B$, and so on.

Note that the objects $\phi(A)\times\phi(B)$, $\Eq_{\phi{}B}$, and so on, do not constitute a $\wedgeq$-cleavage of (the underlying $\wedgeq$-fibration of) $\fib{C'}$, as we have not made the necessary choices for \emph{all} the objects and morphisms in $\fib{C'}$, but rather only for those in the image of $(\Phi,\phi)$ -- and moreover, our choices can ``conflict'', in the sense that if $\phi(A)=\phi(B)$, we needn't have (for example) $\top_{\phi{}A}=\top_{\phi{}B}$.

However, the notation is still ``correct'' in the sense that, for example, $\phi(B)\times\phi(B)$ is indeed a product of $\phi(B)$ with itself, and $\Eq_{\phi{}B}$ is indeed an equality object for $\phi{}B$.

Also, note that there is a potential ambiguity of the following kind. When we write $\hid_{\phi{}f}\colon\top_{\phi{}A}\to\Eq_{\phi{}B}$, do we mean the composite $\top_{\phi{}A}\tox{\exx_{\phi{}f}}\top_{\phi{}B}\tox{\rho_{\phi{}B}}\Eq_{\phi{}B}$, or the image under $\Phi$ of $\hid_f\colon\top_A\to\Eq_B$? The point, of course, is that the two possibilities always agree, as must be checked separately in each case which arises below (and which we leave to the reader).

Now, suppose we are given an extension $\fib{C}\to\fib{C'}$ of $(\Phi,\phi)$ (which we also denote by $(\Phi,\phi)$).

We then have for any 2-cell $\alpha$ in $\B$ that $\Phi\lft\alpha=\lft{\phi{\alpha}}$. Since $\fib{C'}$ is univalent, this uniquely determines the 2-cell $\phi\alpha$. This shows uniqueness of the extension $\phi$.

Next, given a 2-cell
$
\begin{tikzcd}
  P\ar[r, bend left, "p" name=f]\ar[r, bend right, "q"' name=g]
  \ar[Rightarrow, "\beta", from=f, to=g, shorten=2pt]&Q
\end{tikzcd}
$
in $\C$ over
$
\begin{tikzcd}
  A\ar[r, bend left, "f" name=f]\ar[r, bend right, "g"' name=g]
  \ar[Rightarrow, "\alpha", from=f, to=g, shorten=2pt]&B
\end{tikzcd}
$
in $\B$, we have by Proposition~\ref{prop:univalent-2cell-related} that $p$ and $q$ are $\lft\alpha$ related, i.e., that $\nat_B^Q\cdot\brr{\cind{p},\lft\alpha!_P}=\cind{q}$. It then follows that $\nat_{\phi{}B}^{\Phi{}Q}\cdot\brr{\cind{\Phi{}p},\lft{\phi\alpha}!_{\Phi{}P}}=\cind{\Phi{}q}$ -- i.e., that $\Phi{}p$ and $\Phi{}q$ are $\lft{\phi\alpha}$-related. Thus, by Proposition~\ref{prop:univalent-2cell-related} again there is an (of course, unique) lift $\Phi{}p\to\Phi{}q$ of $\phi\alpha$, and this must be $\Phi\beta$. This shows uniqueness of the extension $\Phi$.

It remains to see that the above prescriptions really define a morphism of $\wedgeq$-1D2Fs. It is clear from the above description that $\fib{C'}\Phi=\phi\fib{C}$, and since we are already assuming that $(\Phi,\phi)$ is a morphism of $\wedgeq$-1D2Fs when restricted to $\abs{\fib{C}}$, it only remains to see that $\phi$ and $\Phi$ are 2-functors, i.e.\ that they preserve horizontal and vertical composition and identity 2-cells.

In each case, this is a matter of inspection, the point being that the constructions in $\fib{C}$ from Propositions~\ref{prop:univalent-hvcmp},~\ref{prop:univalent-hhcmp},~and~\ref{prop:univalent-identities} are taken by the $\wedgeq$-morphism $(\Phi,\phi)$ to the corresponding constructions in $\fib{C'}$.

For example, given morphisms $f,g,h\colon{}A\to{}B$ and 2-cells $f\tox{\alpha}g\tox{\beta}h$, we have by Proposition~\ref{prop:univalent-hvcmp} that
$\lft{\beta\cdot\alpha}=\lft\beta\hvcmp\lft\alpha=
\ \ct\!\cdot\tr_B\cdot\brr{\cind{\lft\beta},\cind{\lft\alpha}}$,
the image of which under $\Phi$ is precisely
$\ct\cdot\tr_{\phi{}B}\cdot\brr{\cind{\lft{\phi\beta}},\cind{\lft{\phi\alpha}}}=
\lft{\phi\beta}\hvcmp\lft{\phi\alpha}$,
which by Proposition~\ref{prop:univalent-hvcmp} again is $\lft{\phi\beta\cdot\phi\alpha}$. On the other hand, we have that $\Phi\lft{\beta\cdot\sigma}=\lft{\phi(\beta\cdot\alpha)}$. Hence, by univalence, $\phi\beta\cdot\phi\alpha=\phi(\beta\cdot\alpha)$ as desired.
\qed

\subsubsection{}\thm\label{thm:univalent-main-thm}
Given a $\wedgeq$-cloven $\wedgeq$-fibration $\abs{\fib{C}}$, the inclusion $\abs{\fib{C}}\hookrightarrow\fib{C}$ into the (by Proposition~\ref{prop:construction-is-univalent} univalent) $\wedgeq$-1D2F of Theorem~\ref{thm:2cat-psf} is universal among morphisms to univalent 1D2Fs, i.e., any other such morphism $\abs{\fib{C}}\to\fib{D}$ factors uniquely through $\fib{C}$.

In particular, $\fib{C}$ is, up to isomorphism, the unique extension of $\abs{\fib{C}}$ to a univalent 1D2F.
\pf
Immediate from Proposition~\ref{prop:univalent-unique-extension}.
\qed

\subsubsection{}\thm\label{thm:univalent-stronger-version}
\begin{enumerate}[(i)]
\item The forgetful functor from the category of univalent $\wedgeq$-1D2Fs to the category of $\wedgeq$-fibrations is an equivalence.
\item The category of univalent $\wedgeq$-1D2Fs is a reflexive subcategory of the category of $\wedgeq$-1D2Fs (i.e., the inclusion has a left adjoint).
\item The forgetful functor from the category of $\wedgeq$-1D2Fs to the category of $\wedgeq$-fibrations has a left adjoint with object function given by (choosing a $\wedgeq$-cleavage for each $\wedgeq$-fibration and applying) Theorem~\ref{thm:2cat-psf}.
\end{enumerate}
\pf
That the forgetful functor in (i) is fully faithful follows from Proposition~\ref{prop:univalent-unique-extension}, and that it is surjective on objects follows from Proposition~\ref{prop:construction-is-univalent}, which also shows that there is a right inverse sending each $\wedgeq$-fibration to the $\wedgeq$-1D2F of Theorem~\ref{thm:2cat-psf} (with respect to some $\wedgeq$-cleavage).

To show that the univalent $\wedgeq$-1D2Fs form a reflexive subcategory, we need to find, for each $\wedgeq$-1D2F $\fibr{C}CB$, a morphism into this subcategory which is universal among such morphisms. Choose a $\wedgeq$-cleavage on the underlying $\wedgeq$-fibration $\abs{\fib{C}}$ of $\fib{C}$, and let $\fib{C'}$ be the 1D2F with underlying $\wedgeq$-fibration $\abs{\fib{C}}$ given by Theorem~\ref{thm:2cat-psf}. It follows from Propositions~\ref{prop:univalent-hvcmp},~\ref{prop:univalent-hhcmp},~and~\ref{prop:univalent-identities} that the construction $\alpha\mapsto\lft\alpha$ extends $\id_{\abs{\B}}$ to a morphism of 2-categories from $\B$ to the base category of $\fib{C'}$, and by Proposition~\ref{prop:univalent-2cell-related}, this extends uniquely to a morphism $r\colon\fib{C}\to\fib{C'}$ of $\wedgeq$-1D2Fs extending $\id_{\abs{\fib{C}}}$.

We claim that this $r$ is universal. Indeed, given any other univalent $\wedgeq$-1D2F $\fib{D}$, we have a commutative triangle
\[
  \begin{tikzcd}[column sep=-10pt]
    \Hom(\fib{C}',\fib{D})\ar[rr, "r\circ\mathord{\text{--}}"]\ar[rd]&&
    \Hom(\fib{C},\fib{D})\ar[ld]\\
    &\Hom(\abs{\fib{C}},\fib{D})&
  \end{tikzcd}
\]
of hom sets in the category of $\wedgeq$-1D2Fs, where the two diagonal (restriction) maps are bijections by Proposition~\ref{prop:univalent-unique-extension}, and hence the horizontal map is also a bijection, as desired.

Claim (iii) now follows from (i) and (ii) since the forgetful functor in question is the composite of the right-adjoint from (ii) (whose left-adjoint is the inclusion) and the equivalence from (i) (whose inverse is given on objects by Theorem~\ref{thm:2cat-psf}).
\qed

\subsubsection{}\cor\label{cor:univalent-other-universal-property}
Let $\abs{\fib{C}}$ be a $\wedgeq$-cloven $\wedgeq$-fibration, and $\fib{C}$ the $\wedgeq$-1D2F of Theorem~\ref{thm:2cat-psf}.

Let us say that a \emph{co-extension} of $\abs{\fib{C}}$ is a pair $(\fib{D},f)$ with $\fib{D}$ a $\wedgeq$-1D2F and $f\colon\abs{\fib{D}}\to\abs{\fib{C}}$ a morphism from the underlying $\wedgeq$-fibration of $\fib{D}$.

Then $(\fib{C},\id_{\abs{\fib{C}}})$ is a universal co-extension of $\fib{C}$, i.e., given any other co-extension $(\fib{D},f)$, there is a unique morphism $g\colon\fib{D}\to\fib{C}$ with $f=\id_{\abs{\fib{C}}}\cdot\abs{g}$, where $\abs{g}\colon\abs{\fib{D}}\to\abs{\fib{C}}$ is the restriction of $g$.
\pf
This is just a restatement of the adjunction from Theorem~\ref{thm:univalent-stronger-version}~(iii).
\qed

\section{Examples of $\wedgeq$-fibrations}\label{sec:modelcats}
In Part~\ref{sec:modelcats}, we give some examples of $\wedgeq$-fibrations, to which the results of Part~\ref{sec:proof} can be applied.\vspace{-3pt}

The main examples of fibrations are the ``codomain fibrations'' $\fibr{F(\C)}{C^\to}{C}$ (see Definition~\ref{defn:cod-prefib}), and variations thereof. $\fib{F(\C)}$ is a $\wedgeq$-fibration precisely when $\C$ has finite limits, so this gives a large class of examples. However, in this case the resulting 2-categorical structures is trivial (see \S\ref{subsubsec:cod-fib-triv}).

The main examples of interest to us (of which the above example is, in fact, a special case) come from a \emph{Quillen model category} $\C$ (see Definition~\ref{defn:model-cat}). From $\C$, we form a variation of $\fib{F}(\C)$, namely by constructing a fibration $\fib{HoF(\C)}$ over $\C$ whose fiber over $A$ is the \emph{homotopy category} (see Proposition~\ref{prop:hoc-props}) of $\C/A$ (i.e., of $\fib{F}(\C)^A$).

In \cite{warrenthesis}, it is already observed that the category of fibrations in a model category gives rise to a Grothendieck fibration (our $\fib{F_\fb(\C)}$ -- see Definition~\ref{defn:mfib-subfibs}) and that -- under certain conditions -- the path objects (see Definition~\ref{defn:modc-htpy}) satisfy a ``weak'' analogue of the defining property of equality objects, which is relevant to dependent type theory -- roughly, they satisfy the ``existence'' but not the ``uniqueness'' part of the universal property. The reason for this ``weakness'' is (again roughly speaking) that the universal property is really a ``homotopical one'' -- it involves a \emph{homotopy-equivalence} rather than a \emph{bijection} on $\Hom$-sets. Hence, one might hope to achieve the stronger universal property by passing to the quotient by the homotopy relation -- which is precisely what we do.

Of course, the name ``homotopies'' for the 2-cells defined in Part~\ref{sec:proof} is motivated by these examples -- indeed, two morphisms in $\C_\cfb$ are homotopic with respect to the associated fibration if and only if they are homotopic in the sense of the model structure on $\C$. In Part~\ref{sec:classical-2cat}, we will spell this out in detail.

We should mention an important caveat. Though $\fib{HoF(\C)}$ is always a $\wedge$-fibration, and always has equality objects, it seems that in order for the equality objects to satisfy Frobenius reciprocity and stability along product projections, we need to restrict to the fibrant objects of $\C$. In an earlier version of this paper, we had also required that $\C$ be \emph{right-proper} -- i.e., that weak equivalence are closed under pullbacks along fibrations -- but the anonymous referee pointed out that this is unnecessary since we are restricting to the fibrant objects (see Proposition~\ref{prop:right-proper-between-fibs}).

Part~\ref{sec:modelcats} is organized as follows. In \S\ref{subsec:modelcats-slices}, we recall the definition of the codomain fibration $\fib{F(\C)}$. In \S\ref{subsec:modelcats-modelcats}, we recall the definition of and some basic facts about model categories. In \S\ref{subsec:modelcats-hofib}, we define the fibration $\fib{HoF(\C)}$, in \S\ref{subsec:modelfib-hofib-wedge}, we show that it is a $\wedge$-fibration, and in \S\ref{subsec:modelfib-hofib-wedgeq}, we show that its restriction to $\C_\fb$ is a $\wedgeq$-fibration.

\subsection{Codomain fibrations}\label{subsec:modelcats-slices}
We recall the simplest examples of fibrations, namely the ``codomain'' or ``family'' fibration $\fibr{F(\C)}{\C^\to}{\C}$, for any category $\C$ with pullbacks. The name ``family fibration'' comes from the fact that a morphism $X\to{}A$ (i.e., an object in $\C^\to$ over $A$) can also be thought of as a family of sets indexed by $A$.

These will serve as a fairly uninteresting example of $\wedgeq$-fibrations, but more importantly will serve as the basis for the more interesting examples below.

For the rest of \S\ref{subsec:modelcats-slices}, let $\C$ be a category.

\subsubsection{}\label{defn:cod-prefib}
\defn
The \emph{arrow category} $\C^\to$ of $\C$ has objects triples $(X,A,x)$, with $A,X\in\Ob\C$ and $x\colon{}X\to{}A$ a morphism of $\C$, and the morphisms $(X,A,x)\to(Y,B,y)$ are pairs $(p\colon{}X\to{}Y,f\colon{}A\to{}B)$ such that $yp=fx$. We denote by $\fib{F(\C)}\colon\C^\to\to\C$ the ``codomain functor'', which takes $(X,A,x)$ to $A$ and $(p,f)$ to $f$.

For an object $A\in\Ob\C$, the \emph{slice category} $\C/A$ of $\C$ over $A$ is the fiber $\fib{F(\C)}^A$ over $A$ of the prefibration $\fibr{F(\C)}{\C^\to}{\C}$ (i.e., the subcategory of $\C^\to$ with objects $(X,A,x)$ and morphisms $(p,\id_A)$).

We will sometimes write $(X,x)$ instead of $(X,A,x)$ and $p$ instead of $(p,\id_A)$.

It is easy to see that a morphism in $\C^\to$ is cartesian if and only if it (seen as a square in $\C$) is a pullback square. It follows that $\fib{F(\C)}$ is a fibration if and only if it has pullbacks, and is in fact in this case a $\wedge$-fibration (the only non-trivial thing being stability of fiberwise products -- this amounts to a certain face of a commutative cube being a pullback square, which follows from certain other faces being pullback squares).

It is also easy to see that a morphism $(p,f)\colon(X,A,x)\to(Y,B,y)$ in $\C^\to$ is cocartesian if and only if $p\colon{}X\to{}Y$ is an isomorphism. Hence every morphism $f$ in $\B$ has a cocartesian lift $(\id_P,f)$ with domain any given $P$.

It follows that if $\C$ has, in addition to pullbacks, a terminal object (and is hence an f.p.\ category), then $\fib{F(\C)}$ is a $\wedgeq$-fibration (see \cite[pp.~81,193]{jacobscatlogic}). The stability of the cocartesian lifts and Frobenius reciprocity both amount to a certain edge of a certain commutative cube being an isomorphism, which follows from certain other edges being isomorphisms and certain faces being pullback squares.

\subsubsection{}\label{subsubsec:cod-fib-triv}
We now observe that all the $\fib{F(\C)}$-homotopies are identity $\fib{F(\C)}$-homotopies and hence, the 2-categorical structure induced on $\C$ by Theorem~\ref{thm:2cat} is trivial.

Indeed, given any equality object $\rho_B=(p,\Delta_B)\colon(B',B,b)=\top_B\to\Eq_B=(Y,B\times{}B,y)$, we have that $b$ and $p$ are isomorphisms and hence that $y\colon{}Y\to{}B\times{}B$ is a diagonal morphism. It follows that for $f,g\colon{}A\to{}B$, there can be at most one $\fib{F(\C)}$-homotopy $(q,\br{f,g})\colon(A',A,a)\to(Y,B\times{}B,y)$, and that it exists if only if $f=g$.

\subsection{Model categories}\label{subsec:modelcats-modelcats}
We now review some elements of the theory of model categories. These were introduced in \cite{quillen-ha} as an abstract framework for homotopy theory. This will be fairly brief, and we refer to \cite{may-ponto,hovey} for more background.

\subsubsection{}\label{defn:model-cat}
\defn
Given morphisms $i\colon{}A\to{}B$ and $p\colon{}X\to{}Y$ in a category $\C$, we say that $i$ satisfies the \emph{left lifting property} with respect to $p$, and $p$ satisfies the \emph{right lifting property} with respect to $i$ if for every commutative solid diagram
\[
  \begin{tikzcd}
    A\ar[d, "i"']\ar[r]&X\ar[d, "p"]\\
    B\ar[r]\ar[ru, dashed]&Y,
  \end{tikzcd}
\]
there exists a dashed morphism making the whole diagram commute.

A \emph{weak factorization system} in a category $\C$ consists of two sets $\mathcal{L},\mathcal{R}$ of morphisms of $\C$ such that (i) any morphism $f$ of $\C$ admits a factorization $f=pi$ with $i\in\mathcal{L}$ and $p\in\mathcal{R}$, and (ii) a morphism of $\C$ is in $\mathcal{L}$ (resp.\ in $\mathcal{R}$) if and only if it satisfies the left (resp.\ right) lifting property with respect to every morphism in $\mathcal{R}$ (resp.\ in $\mathcal{L}$).

A \emph{model structure} on a category $\C$ consists of three sets $\mathcal{C},\mathcal{F},\mathcal{W}$ of morphisms of $\C$, called the \emph{cofibrations}, \emph{fibrations},\footnote{Of course, this means we are now considering two different notions called ``fibration'' -- Grothendieck fibrations, and fibrations in a model category. However, this shouldn't cause any confusion.} and \emph{weak equivalences} of the model structure, such that (i) both $(\mathcal{C}\cap\mathcal{W},\mathcal{F})$ and $(\mathcal{C},\mathcal{F}\cap\mathcal{W})$ are weak factorization systems, and (ii) given a commutative diagram
\[
  \begin{tikzcd}[row sep=5pt]
    &C\ar[rd, "g"]&\\
    A\ar[ru, "f"]\ar[rr, "h"']&&B,
  \end{tikzcd}
\]
in which two of the morphisms are weak equivalences, the third is as well.

We refer to property (ii) as the ``two-of-three axiom''.

A \emph{model category} is a category $\C$ having finite limits and colimits, together with a model structure.

We note that this is what Quillen originally called a \emph{closed model category} \cite{quillen-ha} but is now normally just called a model category. The above is a slight reformulation of the definition from \cite{quillen-ha}, and can be found, e.g., in \cite[p.~427]{joyal-thy-of-quasi}. %
We note that sometimes (for example in \cite{hovey}), it is demanded that $\C$ be (not just finitely) complete and cocomplete, and that the cofibration-fibration factorizations are given by specified functors $\fib{C}^\to\to\fib{C}^{\cdot\to\cdot\to\cdot}$.

We will make the usual abuse of notation of identifying a model category with its underlying category.

\subsubsection{}
Let us fix some notational conventions concerning categories with finite coproducts.

For a category $\C$ having specified finite coproducts, we will use the notation $\init_\C$ to denote the chosen initial object of $\C$, and $\ex_A\colon\init_\C\to{}A$ to denote the unique morphism from $\init_\C$ to $A$. Given $A,B\in\Ob\C$, we will denote the chosen coproduct of $A$ and $B$ by $A+B$ and, for morphisms $f\colon{}A\to{}C$ and $g\colon{}B\to{}C$, we will denote by $[f,g]\colon{}A+B\to{}C$ the map induced by $f$ and $g$. We denote by $\nabla_A$ the \emph{codiagonal morphism} $[\id_A,\id_A]\colon{}A+A\to{}A$.

As usual (see \S\S\ref{subsubsec:cleavage},~\ref{subsubsec:fpnotations}~and~\ref{subsubsec:wedge-fib}), we may still use this notation even if finite coproducts are only assumed to exist but have not been specified, but in this case it will be merely suggestive.

\subsubsection{}\defn
An object $A$ in a model category $\C$ is \emph{fibrant} if the unique morphism $!_A\colon{}A\to\tm_\C$ to some (and hence -- since the fibrations include isomorphisms and are closed under composition -- any) terminal object is a fibration, and is \emph{cofibrant} if the unique morphism $\ex_A\colon\init_\C\to{}A$ from some (and hence any) initial object is a cofibration.

We denote by $\C_\cf$ (resp.\ $\C_\fb$ and $\C_\cfb$) the full subcategory on the fibrant (resp.\ cofibrant and cofibrant-fibrant) objects.

A morphism $f\colon{}A\to{}B$ in $\C$ is called a \emph{trivial fibration} if it is both a fibration and a weak equivalence, and a \emph{trivial cofibration} if it is both a cofibration and a weak equivalence.

By using the factorization axiom on the morphisms $\ex_A\colon\init_\C\to{}A$ and $!_A\colon{}A\to\tm_\C$, we can always find a trivial fibration $q\colon{}QA\to{}A$ with $QA$ cofibrant and a trivial cofibration $r\colon{}A\to{}RA$ with $RA$ fibrant. These are called \emph{cofibrant} and \emph{fibrant replacements} for $A$.

\subsubsection{}\label{defn:modc-htpy}
\defn
Given an object $A$ in a model category $\C$, a \emph{cylinder object} for $A$ is a factorization $A+A\tox{[\partial_1,\partial_2]}A\times{}I\tox{\sigma}A$ of a codiagonal map $\nabla\colon{}A+A\to{}A$, in which $[\partial_1,\partial_2]$ is a cofibration and $\sigma$ is a weak equivalence.

Note that by the factorization axiom, every $A\in\Ob\C$ has a cylinder object, and we can even assume that $\sigma$ is a trivial fibration.

We follow \cite{quillen-ha} in using the suggestive notation $A\times{}I$, but this does not mean that the object $A\times{}I$ is really a product.

Similarly, a \emph{path object} for $A$ is a factorization $A\tox{s}A^I\tox{\br{d_1,d_2}}A\times{}A$ of a diagonal map $A\to{}A\times{}A$ with $s$ a weak equivalence and $\br{d_1,d_2}$ a fibration.

Again, there exists a path object for every object, in which $s$ is a cofibration, and again, the notation $A^I$ is merely suggestive.

Given two morphisms $f,g\colon{}A\to{}B$ in $\C$, a \emph{left-homotopy} from $f$ to $g$ is a factorization of the induced map \mbox{$[f,g]\colon{}A+A\to{}B$} through some cylinder object $[\partial_1,\partial_2]\colon{}A+A\to{}A\times{}I$, and we say that $f$ and $g$ are \emph{left-homotopic}, and write $f\siml{}g$, if there exists a left-homotopy between them. Similarly, a \emph{right-homotopy} from $f$ to $g$ is a factorization of $\br{f,g}:A\to{}B\times{}B$ through some path object $\br{d_1,d_2}:B^I\to{}B\times B$, and we say that $f$ and $g$ are \emph{right-homotopic}, and write $f\simr{}g$, if there exists a right-homotopy between them.

By Proposition~\ref{prop:modc-facts}~\ref{item:modc-cyl-choose} below and its dual, if $A$ is cofibrant and $B$ is fibrant, then the relations $\siml$ and $\simr$ on $\Hom(A,B)$ agree and are an equivalence relation. In this case, we write $\sim$ for the relation $\siml=\simr$, and $\pi(A,B)$ for the quotient $\Hom_\C(A,B)/\sim$.

\subsubsection{}\label{prop:modc-facts}\prop
Let $\C$ be a model category, $A,B,C\in\Ob\C$, and let $f,g\colon{}A\to{}B$ and $h\colon{}B\to{}C$ be morphisms in $\C$. Each statement below also comes with a dual statement (in which the direction of morphisms are reversed, the words ``fibration'' and ``cofibration'' are interchanged, and so on).
\claim
\begin{enumerate}[(i)]
\item\label{item:modc-eqrel} If $A$ is cofibrant, then $\siml$ is an equivalence relation on $\Hom(A,B)$.

\item\label{item:modc-cyl-choose} If $B$ is fibrant and $A\times{}I$ is a cylinder object for $A$, then $f\simr{}g$ implies that there is a left-homotopy $A\times{}I\to{}B$ from $f$ to $g$. (In particular, if $B$ is fibrant, then $f\simr{}g$ implies $f\siml{}g$.)
\end{enumerate}
\pf
See \cite[p.~9,~Proposition~1.2.5]{hovey}.
\qed

\subsubsection{}\defn
Given categories $\C$ and $\D$ and a set $W\subset{}\Ar\C$ of morphisms in $\C$, we say that a functor $F\colon\C\to\D$ is a \emph{localization of $\C$ at $W$} (and by abuse of notation, also that $\D$ is the localization of $\C$ at $W$) if (i) $F$ takes morphisms in $W$ to isomorphisms and (ii) given any functor $F\colon\C\to\D'$ satisfying (i), there is a unique\footnote{There is a natural, weaker notion of localization, in which this uniqueness is guaranteed only up to isomorphism.} functor $G\colon\D\to\D'$ with $GF=F'$.

We now recall the definition and main properties of the homotopy category of a model category.

\subsubsection{}\label{prop:hoc-props}\prop
Given a model category $\C$, there exists an essentially unique category $\Ho(\C)$ with $\Ob\Ho(\C)=\Ob\C$ and functor $\gamma\colon\C\to\Ho(\C)$ which is the identity on objects and which is a localization of $\C$ at the weak equivalences.

``Essentially unique'' means: given another such $\gamma'\colon\C\to\Ho'(\C)$, there is a unique functor $F\colon\Ho(\C)\to\Ho'(\C)$ such that $F\gamma=\gamma'$, and moreover $F$ is an isomorphism which is the identity on objects.

Moreover, $\Ho(\C)$ has the following properties:
\begin{enumerate}[(i)]
\item\label{item:hoc-props-weq-iso} For a morphism $f$ in $\C$, $\gamma{}f$ is an isomorphism if and only if $f$ is a weak equivalence.
\item\label{item:hoc-props-cof-fib-htpy} If $f,g\colon{}A\to{}B$ are morphisms in $\C$, with $A$ cofibrant and $B$ fibrant, then $\gamma{}f=\gamma{}g$ if and only if $f\sim{}g$; i.e., $\gamma$ induces a bijection $\pi(A,B)\to\Hom_{\Ho(\C)}(A,B)$.
\item\label{item:hoc-props-cf-equiv} Denoting by $\Ho(\C_*)$, for $*\in\{\cf,\fb\,\cfb\}$, the full subcategory of $\Ho(\C)$ on the objects in $\Ob\C_*\subseteq\Ob\C=\Ob\Ho(\C)$, we have that the restriction $\gamma\colon\C_*\to\Ho(\C_*)$ is a localization of $\C_*$ at the weak equivalences, and that the inclusion $\Ho(\C_*)\hookrightarrow\Ho(\C)$ is an equivalence (and hence also that the inclusions $\Ho(\C_{\cfb})\hookrightarrow\Ho(\C_\cf),\Ho(\C_\fb)$ are equivalences).
\end{enumerate}
\pf
The essential uniqueness is immediate from the definition of localization. It is easy to see that a localization is always a bijection on objects, and being the identity on objects can of course be arranged.

As for existence, it is well-known that a localization at any set of morphisms always exists (though it can lead from a category that is locally small not to one that is not; in this case, this is precluded by \ref{item:hoc-props-cof-fib-htpy} and \ref{item:hoc-props-cf-equiv}).

One direction of \ref{item:hoc-props-weq-iso} is trivial. For the other direction, see \cite[p.~11,~Proposition~1.2.8]{hovey}. For \ref{item:hoc-props-cof-fib-htpy}~and~\ref{item:hoc-props-cf-equiv}, see \cite[I~p.~1.13]{quillen-ha}.
\qed

\subsection{The fibration $\fib{HoF(\C)}$}\label{subsec:modelcats-hofib}
In this section, we define the prefibration $\fibr{HoF(\C)}{\Ho(\C)}{\C}$ by introducing a model structure on $\C^\to$, originally due to A. Roig (in a more general setting, see \cite{roig-bifibrations}) and passing to its homotopy category.

We note that in an earlier version of this paper, we had sought in vain for such a model structure, and were forced instead to take a much more circuitous route to the definition of $\fib{HoF(\C)}$. This model structure was brought to our attention by P. Cagne, and was also known to the anonymous referee.

In \S\ref{subsubsec:pseudo-functor-approach}, we will describe a possible alternative construction of $\fib{HoF(\C)}$.

For the rest of \S\ref{subsec:modelcats-hofib}, let $\C$ be a model category.

\subsubsection{}\defn
We define a model structure on $\C^\to$ (which has all finite limits and colimits since $\C$ does) as follows. A morphism $(p,f)\colon(X,A,x)\to(Y,B,y)$ in $\C^\to$ is
\begin{itemize}
\item a fibration if ``the'' induced map $X\to{}A\times_B{}Y$ is a fibration i.e., if $(p,f)$ factors as $(\ct,f)(p',\id_A)$, with $\ct$ cartesian and $p'$ a fibration.
\item a cofibration if $p$ is a cofibration.
\item a weak equivalence if $f$ is an isomorphism and $p$ is a weak equivalence.
\end{itemize}
The verification that this is a model structure is straightforward, though somewhat lengthy, and in any case it follows from the more general theorems of \cite{roig-bifibrations,stanculescu-bifibrations} on model structures in bifibrations.

Whenever we refer to $\C^\to$ as a model category, we mean with respect to this model structure.

\subsubsection{}\prop\label{prop:slice-total-mcat-compat}
A morphism $(f,\id_A)\colon(X,A,x)\to{}(Y,A,x)$ in $\C^\to$ is a fibration, cofibration, or weak equivalence if and only if $f$ is. These make $\C/A$ into a model category. Moreover, an object in $\C/A$ is fibrant or cofibrant if and only if it is as an object of $\C^\to$. Explicitly, $(X,x)$ is fibrant if and only if $x$ is a fibration and is cofibrant if and only if $X$ is cofibrant.

Whenever we refer to $\C/A$ as a model category, we mean with respect to this model structure.
\pf
The proof of the first claim is by inspection. That this defines a model structure on $\C/A$ is well-known (see, e.g., \cite[p.~5]{hovey}) and easy to verify. The proofs of the last two claims are also by inspection.
\qed

\subsubsection{}\prop\label{prop:fiberwise-homotopy}
If two morphisms $(p,f),(q,g)\colon(X,A,x)\to(Y,B,y)$ in $\C^\to$ are left-homotopic, then $f=g$. Moreover, $(p,f),(q,f)\colon(X,A,x)\to(Y,B,y)$ are left-homotopic if and only if there exists a homotopy $X\times{}I\to{}Y$ from $p$ to $q$ making the following diagram commute.
\[
  \begin{tikzcd}
    X\times{}I\ar[d, "\sigma"']\ar[rd, "h"]&\\
    X\ar[d, "x"']&Y\ar[d, "y"]\\
    A\ar[r, "f"]&B
  \end{tikzcd}
\]

In particular, in the case that $(X,A,x)$ is cofibrant and $(Y,A,y)$ is fibrant, the set $\Hom^{\Ho(\C^\to)}_f((X,A,x),(Y,B,y))$ is the quotient of $\Hom^{\C^\to}_f((X,A,x),(Y,B,y))$ by the above relation.
\pf
Let us set $P=(X,A,x)$ and $Q=(Y,B,y)$.

We know that the coproducts $P+P$ in $\C^\to$ are exactly the objects of the form $(X+X,A+A,x+x)$. Next, it follows from the definitions that $P+P\tox{[\partial_1',\partial_2']}P\times{}I\tox{\sigma'}P$ is a cylinder object if and only if $P\times{}I$ is of the form $(X\times{}I,A',e\I{}x\sigma)$ with $\sigma'=(\sigma,e)$ and $[\partial_1',\partial_2']=([\partial_1,\partial_2],\nabla_Ae\I)$, for some cylinder object $X+X\tox{[\partial_1,\partial_2]}X\times{}I\tox{\sigma}X$ and some isomorphism $e\colon{}A'\tox{\sim}A$.

Given $P\times{}I$ of this form, a left-homotopy $P\times{}I\to{}Q$ from $(p,f)$ to $(g,q)$ is then a morphism $(h,h')\colon(X\times{}I,A',e\I{}x\sigma)\to(Y,B,y)$ such that the triangles
\[
  \begin{tikzcd}
    &X\times{}I\ar[dr, "h"]&\\
    X+X\ar[rr, "{[p,q]}"]\ar[ru, "{[\partial_1,\partial_2]}"]&&Y
  \end{tikzcd}
  \hspace{45pt}
  \begin{tikzcd}
    &A'\ar[dr, "h'"]&\\
    A+A\ar[ru, "e\I\nabla_A"]\ar[rr, "{[f,g]}"]&&Q
  \end{tikzcd}
\]
commute (in particular, $h'=fe=ge$). This proves both claims. \qed

\subsubsection{}\prop\label{prop:slice-total-htpy-compat}
Two morphisms $(p,\id_A),(q,\id_A)\colon(X,A,x)\to(Y,A,y)$ in $\C^\to$ are left-homotopic in $\C^\to$ if and only if they are in $\C/A$.
\pf
We show that $p$ and $q$ are homotopic in $\C/A$ if and only if they satisfy the condition described in Proposition~\ref{prop:fiberwise-homotopy}. The argument is similar to the one there, but simpler: one first verifies that the cylinder objects in $\C/A$ on $(X,x)$ are exactly of the form
\[
  (X+X,[x,x])\tox{[\partial_1,\partial_2]}(X\times{}I,x\sigma)\tox{\sigma}(X,x)
\]
with $X+X\tox{[\partial_1,\partial_2]}X\times{}I\to{}\tox{\sigma}{}X$ a cylinder object in $\C$, and the claim follows by inspecting the definition of left-homotopy.
\qed

\subsubsection{}\defn
Given two prefibrations $\fibr{D}DB$ and $\fibr{D'}{D'}B$ over a category $\B$, a \emph{morphism of prefibrations} $F\colon\fib{D}\to\fib{D'}$ is a functor $F\colon\D\to{}\D'$ for which
\[
  \begin{tikzcd}[column sep=7pt]
    \D\ar[rd]\ar[rr, "F"]&&\D'\ar[ld]\\
    &\B&
  \end{tikzcd}
\]
commutes. Note that by restriction, $F$ induces functors $\fib{D}^A\to\fib{D'}^A$ for each $A\in\Ob\B$.

If $\fib{D}$ and $\fib{D'}$ are fibrations, then $F$ is a \emph{morphism of fibrations} if it takes cartesian morphisms to cartesian morphisms.

If $\fib{D}$ and $\fib{D'}$ are $\wedge$-fibrations, then $F$ is a \emph{morphism of $\wedge$-fibrations} if, in addition, the induced functors $\fib{D}^A\to\fib{D'}^A$ are all f.p.\ functors.

\subsubsection{}\label{defn:mfib-subfibs}
\defn
For $*\in\{\cf,\fb,\cfb\}$, we define the pre-fibration $\fibr{F_*(\C)}{(C^\to)_*}{C}$ to be the restriction to $(\C^\to)_*$ of the functor $\fib{F(\C)}\colon\C^\to\to\C$. Note that by Proposition~\ref{prop:slice-total-mcat-compat}, the fiber $\fib{F_*(\C)}^A$ of $\fib{F_*(\C)}$ over $A$ is precisely the category $(\C/A)_*$.

\subsubsection{}\defn
Since the functor $\fib{F(\C)}\colon\C^\to\to\C$ clearly preserves weak equivalences it induces a functor $\Ho(\C^\to)\to\C$, which we denote by $\fibr{HoF(\C)}{\Ho(\C^\to)}{\C}$. Similarly, for $*\in\{\cf,\fb,\cfb\}$, we have a prefibration $\fibr{HoF_*(\C)}{\Ho(\C^\to)_*}{\C}$ induced from $\fib{F_*(\C)}$.

The functor $\gamma\colon\C^\to\to\Ho(\C^\to)$ induces a morphism of prefibrations $\gamma\colon\fib{F(\C)}\to\fib{HoF(\C)}$, and we similarly have morphisms $\gamma\colon\fib{F_*(\C)}\to\fib{HoF_*(\C)}$.

\subsubsection{}\prop\label{prop:slice-total-hoc-compat}
The functors $\C/A=\fib{F(\C)}^A\tox{\gamma}\fib{HoF(\C)}^A$ induced from $\gamma\colon\fib{F(\C)}\to\fib{HoF(\C)}$ factor through $\gamma\colon\C/A\to\Ho(\C/A)$ and induce isomorphisms $\Ho(\C/A)\to\fib{HoF(\C)}^A$.

Similarly, the restrictions $\fib{F_*(\C)}^A\to\fib{HoF_*(\C)}^A$ induce isomorphisms $\Ho(\C/A)_*\to\fib{HoF_*(\C)}^A$ for $*\in\{\cf,\fb,\cfb\}$.
\pf
That the functor $\C/A\to\fib{HoF(\C)}^A$ factors through $\Ho(\C/A)$ means that it takes weak equivalences to isomorphisms, which follows from Proposition~\ref{prop:slice-total-mcat-compat}. The same goes for $(\C/A)_*$.

Next, we claim that the induced functor $\Ho(\C/A)_\cfb\to\fib{HoF_\cfb(\C)}^A$ is an equivalence. It is the identity on objects, so it remains to see that it is fully faithful. Fix $P,Q\in\Ob(\C/A)_\cfb$ and consider the diagram
\[
\begin{tikzcd}
  \Hom_{\C/A}(P,Q)\ar[r, "\gamma"]\ar[d, hook']&\Hom_{\fib{HoF(\C)}^A}(P,Q)\ar[d, hook']\\
  \Hom_{\C^\to}(P,Q)\ar[r, "\gamma"]&\Hom_{\Ho(\C^\to)}(P,Q).
\end{tikzcd}
\]
where we would like to show that the top map induces a bijection from $\Hom_{\Ho(\C/A)}(P,Q)$. Let us first see that the top map is surjective. We will refer, for brevity, to the set in the top-left of the diagram by (TL), the set in the bottom-right by (BR), and so on.

By the definition of $\fib{HoF(\C)}$, (TR) consists of those elements of (BR) which are sent to $\id_A$ under the functor $\Ho(\C^\to)\to\C$ induced from $\fib{F(\C)}\colon\C^\to\to\C$. Since $P$ is cofibrant $Q$ is fibrant, it follows from Proposition~\ref{prop:hoc-props}~\ref{item:hoc-props-cof-fib-htpy} that the bottom map $\gamma$ is surjective, so that (TR) consists of the images of those elements of (BL) sent by $\fib{F(\C)}$ to $\id_A$ -- i.e., the image of (TL). This shows surjectivity.

By Proposition~\ref{prop:slice-total-mcat-compat}, $P$ and $Q$ are cofibrant and fibrant also in $\C^\to$, hence by Proposition~\ref{prop:hoc-props}~\ref{item:hoc-props-cof-fib-htpy} again, the bottom map (and hence the top map) identifies two morphisms if and only if they are left-homotopic in $\C^\to$. But by Proposition~\ref{prop:slice-total-htpy-compat}, this is the same as being left-homotopic in $\C/A$. Hence, by Proposition~\ref{prop:hoc-props}~\ref{item:hoc-props-cof-fib-htpy} once again, the top map induces a bijection $\Hom_{\Ho(\C/A)}(P,Q)\to$(TR).

We now have a commutative diagram of categories and functors
\[
  \begin{tikzcd}
    \Ho(\C/A)_\cfb\ar[r, "\gamma"]\ar[d, hook']&\fib{HoF_\cfb(\C)}^A\ar[d, hook']\\
    \Ho(\C/A)\ar[r, "\gamma"]&\fib{HoF(\C)}^A,
  \end{tikzcd}
\]
and we want to show the bottom arrow is an equivalence. We just showed that the top arrow is an equivalence, and the left arrow is an equivalence by
Proposition~\ref{prop:hoc-props}~\ref{item:hoc-props-cf-equiv}, so it remains to see that the right arrow is an equivalence. Since it is an inclusion of a full subcategory, we only need to see it is essentially surjective. But for every $P\in\Ob(\C/A)$, there is a weak equivalence in $\C/A$ (which is by Proposition~\ref{prop:slice-total-mcat-compat} also a weak equivalence in $\C^\to$) between $P$ and an object in $(\C/A)_\cfb$, the image of which under $\gamma$ gives the desired isomorphism.

The same argument if we replace $(\C/A)$ and $\fib{HoF(\C)}^A$ and $(\C/A)_*$ by $\fib{HoF_*(\C)}^A$ for $*\in\{\cf,\fb\}$.
\qed

\subsubsection{}\label{subsubsec:pseudo-functor-approach}
We now sketch a possible alternative approach to the construction of the fibration $\fib{HoF(\C)}$ which more directly implements the idea ``pass to the homotopy category of each fiber of $\fib{F(\C)}$'', and which has been implemented by P. Cagne in \cite{cagnethesis}.

In \cite[p.26]{hovey}, it is shown that the passage from a model category to its homotopy category is described by a pseudo-functor $\cat{Ho}\colon\cat{Mod}\to\Cat$ from the 2-category $\cat{Mod}$ of model categories and ``Quillen adjunctions'' to the 2-category of categories. Given a model category $\C$, the pseudo-functor $\fpsf{F(\C)}\colon\C^\op\to\Cat$ associated to any cleavage of $\fib{F(\C)}$ factors through the forgetful 2-functor $\cat{Mod}\hookrightarrow\Cat$, since the left-adjoint $\sum_f$ to each pullback functor $f^*$ is a ``left Quillen functor''. According to \cite[Proposition~4.3.3]{cagnethesis}, the fibration $\fib{HoF(\C)}$ is then the one associated to the composition $\cat{Ho}\circ\fpsf{F(\C)}$ of the factored pseudo-functor $\fpsf{F(\C)}\colon\C^\op\to\cat{Mod}$ with the pseudo-functor $\cat{Ho}\colon\cat{Mod}\to\Cat$.

\subsection{$\fib{HoF_\fb(\C)}$ is a $\wedge$-fibration}\label{subsec:modelfib-hofib-wedge}
We will now show that $\fib{HoF(\C)}$ is a $\wedge$-fibration. For this purpose, it will be more convenient to work with the equivalent fibration $\fib{HoF_\fb(C)}$, since $\fib{F_\fb(\C)}$ is itself a $\wedge$-fibration and, as we will show, the morphism $\gamma\colon\fib{F_\fb(\C)}\to\fib{HoF_\fb(\C)}$ is a morphism of $\wedge$-fibrations, which gives us a very explicit description of the $\wedge$-fibration structure of $\fib{HoF_\fb(\C)}$.

\subsubsection{}\label{fact:modc-fib-clos}
We recall that the fibrations in any model category are stable under pullbacks: if
\[
  \begin{tikzcd}
    A\ar[r]\ar[d, "p'"']\ar[rd, phantom, "\lrcorner", pos=0]&B\ar[d, "p"]\\
    C\ar[r]&D
  \end{tikzcd}
\]
is a pullback square in and $p$ is a fibration, so is $p'$. In fact, the elements of $\mathcal{R}$ in any weak factorization system $(\mathcal{L},\mathcal{R})$ are stable under pullbacks (and dually, the elements of $\mathcal{L}$ are stable under pushouts).

In particular, if $A\xot{\pi_1}A\times{}B\tox{\pi_2}B$ is a product diagram and $A$ (resp.\ $B$) is fibrant, then $\pi_2$ (resp.\ $\pi_1$) is a fibration.

\subsubsection{}\label{prop:modc-fb-prods}\prop
For any model category $\C$, the category $\C_{\fb}$ is an f.p.\ category, and the inclusion $\C_\fb\hookrightarrow\C$ is an f.p.\ functor.
\pf
Since $\C_\fb$ is a full subcategory of $\C$, it suffices to see that a finite product of fibrant objects is fibrant. That terminal objects are always fibrant is immediate, and that the binary product of fibrant objects is fibrant follows from \S\ref{fact:modc-fib-clos}.
\qed

\subsubsection{}\label{prop:mfib-fb-wedge}\prop
For any model category $\C$, the prefibration $\fibr{F_\fb(\C)}{(\C^\to)_\fb}{\C}$ is a $\wedge$-fibration, and the inclusion $\fib{F_\fb(\C)}\hookrightarrow\fib{F(\C)}$ is a morphism of $\wedge$-fibrations.
\pf
We first prove the claim with ``$\wedge$-'' removed. Since $(\C^\to)_\fb$ is a full subcategory of $\C^\to$, it suffices to see that any cartesian morphism in $\C^\to$ with codomain in $(\C^\to)_\fb$ has its domain in $(\C^\to)_\fb$. But this follows from \S\ref{fact:modc-fib-clos} since the cartesian morphisms in $\C^\to$ are precisely the pullback squares.

That the fibers of $\fib{F_\fb(\C)}$ have finite products which are preserved by the inclusion $\fib{F_\fb(\C)}^A\hookrightarrow\fib{F(\C)}^A$ follows from Proposition~\ref{prop:modc-fb-prods}. That the finite products in $\fib{F_\fb(\C)}$ are stable under pullbacks is immediate from the corresponding property in $\fib{F(\C)}$.
\qed

\subsubsection{}\label{prop:modc-hofb-prods}\prop
(cf. \cite[Example~1.3.11]{hovey})
For any model category $\C$, the category $\Ho(\C_{\fb})$ (hence also $\Ho(\C)$) is an f.p.\ category, and $\gamma\colon\C_\fb\to\Ho(\C_\fb)$ is an f.p.\ functor.
\pf
Let $\tm$ be terminal in $\C_\fb$. We need to see that for each $A\in\Ob\C_\fb$, there is a unique morphism $A\to\tm$ in $\Ho(\C_\fb)$. It suffices to see this for cofibrant $A$, since every object in $\Ho(\C_\fb)$ is isomorphic to such an $A$. But in this case, by
Proposition~\ref{prop:hoc-props}~\ref{item:hoc-props-cof-fib-htpy}. the morphisms $A\to\tm$ in $\Ho(\C_\fb)$ are just homotopy classes of morphisms $A\to\tm$ in $\C$, of which there is of course just one.

Next, let $B\xot{\pi_1}B\times{}C\tox{\pi_2}C$ be a product in $\C_\fb$. We need to see that for each $A\in\Ob\C_\fb$, composition with $\pi_1$ and $\pi_2$ induces a bijection $\Hom_{\Ho(\C_\fb)}(A,B\times{}C)\to\Hom_{\Ho(\C_\fb)}(A,B)\times\Hom_{\Ho(\C_\fb)}(A,C)$. Again, it suffices to consider $A$ cofibrant, so that we need to show that
\[
  \br{\pi_1\ccmpf\mathord{\text{--}},\pi_2\ccmpf\mathord{\text{--}}}\colon\pi(A,B\times{}C)\to\pi(A,B)\times\pi(A,C)
\]
is a bijection. That it is surjective is immediate, since
\[
  \br{\pi_1\ccmpf\mathord{\text{--}},\pi_2\ccmpf\mathord{\text{--}}}\colon\Hom_{\C_\fb}(A,B\times{}C)\to
  \Hom_{\C_\fb}(A,B)\times\Hom_{\C_\fb}(A,C)
\]
is already surjective pointwise, and not just on homotopy classes.

To see that it is injective, we need to check that given homotopic maps $f_1,f_2\colon{}A\to{}B$ and homotopic maps $g_1,g_2\colon{}A\to{}C$, the induced maps $\br{f_1,g_1},\br{f_2,g_2}\colon{}A\to{}B\times{}C$ are homotopic. By Proposition~\ref{prop:modc-facts}~\ref{item:modc-cyl-choose} we can choose left-homotopies $h_f\colon{}A\times{}I\to{}B$ from $f_1$ to $f_2$ and $h_g\colon{}A\times{}I\to{}C$ from $g_1$ to $g_2$ with a common cylinder object $A+A\tox{[\partial_1,\partial_2]}A\times{}I$. We then have an induced left-homotopy $\br{h_f,h_g}\colon{}A\times{}I\to{}B\times{}C$. To show that this is a left-homotopy between $\br{f_1,g_1}$ and $\br{f_2,g_2}$, we need to show that $\br{h_f,h_g}[\partial_1,\partial_2]=[\br{f_1,g_1},\br{f_2,g_2}]$. But
\[
  \mathrm{LHS}=
  [\br{h_f,h_g}\partial_1,\br{h_f,h_g}\partial_2]=
  [\br{h_f\partial_1,h_g\partial_1},\br{h_f\partial_2,h_g\partial_2}]=
  \mathrm{RHS}.
  \tag*{\qed}
\]

\subsubsection{}\label{prop:mfib-hofb-wedge}\prop
For any model category $\C$, $\fib{HoF_\fb(\C)}$ is a $\wedge$-fibration and $\gamma\colon\fib{F_\fb(\C)}\to\fib{HoF_\fb(\C)}$ is a morphism of \mbox{$\wedge$-fibrations}.
\pf
That the fibers of $\fib{HoF_\fb(\C)}$ are f.p.\ categories, and that functors $\fib{F_\fb(\C)}^A\to\fib{HoF_\fb(\C)}^A$ induced by $\gamma$ are f.p.\ functors, follows from Propositions~\ref{prop:modc-hofb-prods}~and~\ref{prop:slice-total-hoc-compat}.

Let $p\colon{}Q\to{}R$ in $(\C^\to)_\fb$ be a cartesian morphism over $g\colon{}B\to{}C$ in $\C$. We need to see that the image $\gamma{}p$ in $\Ho(\C^\to)_\fb$ is still cartesian; i.e., that for $f\colon{}A\to{}B$ in $\C$ and $P\in\Ob\fib{HoF_\fb(\C)}^A$, the map $(\gamma{}p)\ccmpf\mathord{\text{--}}\colon\Hom_{f}(P,Q)\to\Hom_{gf}(P,R)$ is a bijection. We argue as in Proposition~\ref{prop:modc-hofb-prods}. First, we can assume that $P$ is cofibrant, and hence, by Proposition~\ref{prop:fiberwise-homotopy}, we must show that $(p\ccmpf\mathord{\text{--}})\colon\Hom_{f}(P,Q)\to\Hom_{gf}(P,R)$ induces a bijection on left-homotopy classes. Surjectivity is clear, and injectivity follows by a similar -- but simpler -- argument to the one in Proposition~\ref{prop:modc-hofb-prods}.

It remains to see that the products in the fibers of $\fib{HoF_\fb(\C)}$ are stable under pullbacks. Let $f\colon{}A\to{}B$ be a morphism in $\C$ and $P,Q\in\Ob\fib{HoF_\fb(\C)}^B$. It suffices to see that for \emph{some} product $P\wedge{}Q$, and \emph{some} pullbacks $f^*P$, $f^*Q$, and $f^*(P\wedge{}Q)$ as in
\[
  \begin{tikzcd}
    &[-40pt]f^*(P\wedge{}Q)\ar[ldd, "f^*\pi_1"']\ar[rd, "f^*\pi_2"]\car[rrr, "\ct"]&[-20pt]&
    &[-30pt]P\wedge{}Q\ar[rd, "\pi_2"]&[-20pt]\\
    &&f^*Q\ar[rrr, "\carsym" anchor=center, near start, "\ct" outer sep=1pt]&&&Q\\
    f^*P\car[rrr, "\ct"]&&&P\ar[from=uur, crossing over, "\pi_1", near end]&&\\
    &A\ar[rrr, "f"]&&&B,
  \end{tikzcd}
\]
$f^*P\xot{f^*\pi_1}P\wedge{}Q\tox{f^*\pi_2}f^*Q$ is also a product diagram. Now, by Proposition~\ref{prop:modc-hofb-prods} and what we have just shown, we can obtain such a product $P\wedge{}Q$ and such cartesian morphisms by first choosing a product diagram and cartesian morphisms in $(\C^\to)_\fb$, and then taking their images in $\Ho(\C^\to)_\fb$. But now the diagram $f^*P\xot{f^*\pi_1}P\wedge{}Q\tox{f^*\pi_2}f^*Q$ in $\fib{F_\fb(\C)}^A$ is a product diagram since $\fib{F_\fb(\C)}$ (by Proposition~\ref{prop:mfib-fb-wedge}) is a $\wedge$-fibration. Hence, using Proposition~\ref{prop:modc-hofb-prods} again, $f^*P\xot{f^*\pi_1}P\wedge{}Q\tox{f^*\pi_2}f^*Q$ is a product diagram in $\fib{HoF_\fb(\C)}^A$ as desired.

The proof that terminal objects in $\fib{HoF_\fb(\C)}^B$ are stable under pullback is similar, but simpler.\qed

\subsection{$\fib{HoF_\fb(\C_\fb)}$ is a $\wedgeq$-fibration}\label{subsec:modelfib-hofib-wedgeq}
We now want to show that $\fib{HoF_\fb(\C)}$ is a $\wedgeq$-fibration. We will show that the necessary cocartesian lifts always exist; in fact, as with the fibration $\fib{F(\C)}$, we will show that \emph{any} $f\colon{}A\to{}B$ in $\C$ has a cocartesian lift with domain \emph{any} $P\in\Ob\fib{HoF_\fb(\C)}^A$. However, in order to show that these satisfy Frobenius reciprocity and that they are stable along product projections, it seems that we need to restrict to the fibrant objects of $\C$.

For the rest of \S\ref{subsec:modelfib-hofib-wedgeq}, let $\C$ be a model category.

\subsubsection{}\label{prop:mfib-hofb-cocart-crit}\prop
Given a morphism $(\hat{f},f)\colon(X,A,x)\to(Y,B,y)$ in $(\C^\to)_\fb$, the image $\gamma(\hat{f},f)$ in $\Ho(\C^\to)_\fb$ of $(\hat{f},f)$ is cocartesian if and only if $\hat{f}$ is a weak equivalence.
\pf
It suffices to prove this for $\Ho(\C^\to)$ (rather than $\Ho(\C^\to)_\fb$); since the inclusion $\fib{HoF_\fb(\C)}\to\fib{HoF(\C)}$ is an equivalence on total categories and on each fiber, it follows that it preserves cocartesian morphisms, and also that $\fib{HoF(\C)}$ is a fibration (since $\fib{HoF_\fb(\C)}$ is).

We now factor $(\hat{f},f)$ as $(X,A,x)\tox{(\id_X,f)}(X,B,fx)\tox{(\hat{f},\id_B)}(Y,B,y)$. Let us next see that $\gamma(\id_X,f)$ is cocartesian.

We can assume $X$ is cofibrant by passing to a cofibrant replacement $y\colon{}X'\to{}X$. Indeed, we then have a commutative square
\[
  \begin{tikzcd}[column sep=50pt]
    (X',A,xy)\ar[r, "{\gamma(\id_{X'},f)}"]\ar[d, "{\gamma(y,\id_A)}"']&
    (X',B,fxy)\ar[d, "{\gamma(y,\id_B)}"]\\
    (X,A,x)\ar[r, "{\gamma(\id_X,f)}"]&(X,B,fx)
  \end{tikzcd}
\]
in $\Ho(\C^\to)$ with the vertical arrows isomorphisms, and hence that $\gamma(\id_X,f)$ is cocartesian if and only if $\gamma(\id_{X'},f)$ is.

We need to check that for each $(Z,B,z)\in\Ob\C/B$ (which we can assume fibrant), the map $(\mathord{\text{--}}\circ{}\gamma(\id_X,f))\colon\Hom^{\Ho(\C^\to)}_{\id_B}((X,B,fx),(Z,B,z))\to\Hom^{\Ho(\C^\to)}_f(X,A,x),(Z,B,z))$ is a bijection.

That the corresponding map $(\mathord{\text{--}}\circ(\id_X,f))$, with $\Hom^{\Ho(\C^\to)}$ replaced by $\Hom^{\C^\to}$, is a bijection is obvious, and since $(X,A,x)$ is cofibrant and $(Z,B,z)$ is fibrant, it only remains to see that this bijection preserves the relation $\sim$. This follows from the explicit description of $\siml$ in $\C^\to$ from Proposition~\ref{prop:fiberwise-homotopy}.

Since $\gamma(\id_X,f)$ is cocartesian, $\gamma(\hat{f},\id_B)$ is cocartesian if and only if the composite $\gamma(\hat{f},f)$ is, so it remains to see that $\hat{f}$ is a weak equivalence if and only if $\gamma(\hat{f},\id_B)$ is cocartesian -- i.e., (since it lies over an isomorphism) if and only if $\gamma(\hat{f},\id_B)$ is an isomorphism. But this holds by Proposition~\ref{prop:hoc-props}~\ref{item:hoc-props-weq-iso} and the definition of weak equivalence in $\C^\to$.
\qed

\subsubsection{}\label{prop:mfib-hofb-has-cocarts}\prop
For each $f\colon{}A\to{}B$ in $\C$ and each $P\in\Ob\fib{HoF_\fb(\C)}^A$, there is a cocartesian lift of $f$ in $\Ho(\C^\to)_\fb$ with domain $P$.
\pf
Suppose $P=(X,A,x)$, and factor $fx\colon{}X\to{}B$ in $\C$ as a trivial cofibration $\hat{f}\colon{}X\to{}X'$ followed by a fibration $x'\colon{}X'\to{}B$. Then $(X',B,x')\in\Ob(\C^\to)_\fb$, and by Proposition~\ref{prop:mfib-hofb-cocart-crit}, the image of $(\hat{f},f)\colon(X,A,x)\to(X',B,x')$ in $\Ho(\C^\to)_\fb$ is cocartesian.
\qed

\subsubsection{}\label{prop:right-proper-between-fibs}\prop
Given a pullback square
\[
  \begin{tikzcd}
    A\ar[r]\ar[d, "p'"']\ar[rd, phantom, "\lrcorner", pos=0]&B\ar[d, "p"]\\
    C\ar[r, "g"]&D
  \end{tikzcd}
\]
in $\C$ in which $B$ and $D$ are fibrant, $g$ is a fibration and $p$ is a weak equivalence, $p'$ is also a weak equivalence.
\pf
See \cite[Proposition~13.1.2]{hirschhorn}
\qed

\subsubsection{}\label{prop:mfib-hofb-cocart-stab}\prop
Given a morphism $g\colon{}C\to{}D$ in $\C_\fb$, the cocartesian morphisms in $\Ho(\C^\to)_\fb$ lying over $g$ are stable along every fibration $k\colon{}B\to{}D$.
\pf
It suffices, for each $P\in\Ob\fib{HoF_\fb(\C)}^C$, to see that \emph{some} cocartesian lift of $g$ with domain $P$ is stable along every $k$.

By Proposition~\ref{prop:mfib-hofb-cocart-crit}, we can take as our cocartesian morphism the image $\gamma{}p\colon{}P\to{}Q$ of some morphism $p=(\hat{g},g)\colon{}P\to{}Q$ in $(\C^\to)_\fb$ with $\hat{g}$ a weak equivalence. Again, to see that $\gamma{}p$ is stable along the fibration $k$, it suffices to see that for each pullback square
\[
  \begin{tikzcd}
    A\ar[r, "f"]\ar[d, "h"']\ar[rd, phantom, "\lrcorner", pos=0]&B\ar[d, "k"]\\
    C\ar[r, "g"]&D
  \end{tikzcd}
\]
in $\C$, there exist \emph{some} cartesian morphisms $\ct\colon{}h^*P\to{}P$ over $h$ and $\ct\colon{}k^*Q\to{}Q$ over $k$, such that the unique morphism $p'$ over $f$ making the following diagram commute is cocartesian.
\[
  \begin{tikzcd}
    h^*{}P\ar[r, "p'", dashed]\ar[d, "\ct"']&
    k^*Q\ar[d, "\ct"]\\
    P\ar[r, "p"]&Q
  \end{tikzcd}
\]

Now, by Proposition~\ref{prop:mfib-hofb-wedge}, we can take our cartesian lifts $h^*P\to{}P$ and $k^*Q\to{}Q$ to be the image under $\gamma$ of cartesian lifts of $h$ and $k$ in $(\C^\to)_\fb$ -- where we recall that ``cartesian'' in $(\C^\to)_\fb$ means ``pullback square''.

We thus have a commutative cube
\[
  \begin{tikzcd}
    &\cdot\ar[rr, "\hat{f}"]
    \ar[dd]&&
    \cdot\ar[dd]\\
    \cdot\ar[rr, "\hat{g}", pos=0.7, crossing over]
    \ar[dd]\ar[from=ru]&&
    \cdot\ar[from=ru, "\hat{k}"]\\
    &A\ar[rr, "f", pos=0.3]&&B\\
    C\ar[rr, "g"]\ar[from=ru, "h"']&&D\ar[from=ru, "k"]
    \ar[from=uu, crossing over]
  \end{tikzcd}
\]
in $\C$ in which the right, left, and bottom faces (and hence -- by three application of the ``2-of-3'' rule -- also the top face) are pullback squares, and $\hat{g}$ is a weak equivalence. We want to show that $\hat{f}$ is a weak equivalence (since this would imply, by Proposition~\ref{prop:mfib-hofb-cocart-crit}, that $p'$ is cocartesian). But by \S\ref{fact:modc-fib-clos}, $\hat{k}$ is a fibration, and the domain and codomain of $\hat{g}$ are fibrant (since they admit fibrations to the fibrant objects $C$ and $D$), and hence $\hat{f}$ is a weak equivalence by Proposition~\ref{prop:right-proper-between-fibs}.
\qed

\subsubsection{}\label{prop:mfib-hofb-cocart-frob}\prop
Given a morphism $f\colon{}A\to{}B$ in $\C_\fb$, the cocartesian morphisms in $\Ho(\C^\to)_\fb$ lying over $f$ satisfy Frobenius reciprocity.
\pf
The argument is similar to the one in Proposition~\ref{prop:mfib-hofb-cocart-stab}.

Let $f\colon{}A\to{}B$ be a morphism in $\C$. It suffices to check that for each $Q\in\Ob\fib{HoF_\fb(\C)}^A$ and each $P\in\Ob\fib{HoF_\fb(\C)}^B$, there is \emph{some} cocartesian lift $q\colon{}Q\to{}Q'$ of $f$, \emph{some} cartesian lift $\ct\colon{}f^*P\to{}P$ and \emph{some} product diagrams $Q\xot{\pi_1}Q\wedge{}f^*P\tox{\pi_2}f^*P$ and $Q'\xot{\pi_1}Q'\wedge{}P\tox{\pi_2}P$ for which the induced morphism $q\wwdge\!\ct$ -- i.e., the unique morphism over $f$ making the following diagram commute -- is cocartesian.
\[
  \begin{tikzcd}[column sep=40pt]
    Q\rac[r, "q"]&Q'\\
    Q\wedge{}f^*P\ar[r, "q\wwdge{}\ct", dashed]\ar[u, "\pi_1"]\ar[d, "\pi_2"']&
    Q'\wedge{}P\ar[u, "\pi_1"']\ar[d, "\pi_2"]\\
    f^*P\car[r, "\ct"]&P\\[-10pt]
    A\ar[r, "f"]&B
  \end{tikzcd}
\]

Now, as in Proposition~\ref{prop:mfib-hofb-cocart-stab}, we choose the cartesian morphism $\ct\colon{}f^*P\to{}P$ to be the image of a cartesian lift of $f$ in $(\C^\to)_\fb$, and the cocartesian lift $q$ to be the image a morphism $(\hat{f},f)\colon{}Q\to{}Q'$ in $(\C^\to_\fb)$ with $\hat{f}$ a weak equivalence. Similarly, we choose (using Proposition~\ref{prop:modc-hofb-prods}) the product diagrams to be the images of product diagrams in $(\C/A)_{\fb}$ and $(\C/B)_\fb$ -- which, we recall, are pullback diagrams in $\C$.

We then, as in Proposition~\ref{prop:mfib-hofb-cocart-stab}, end up with a cube in $\C$ in which the same three faces are pullbacks and the same edge is a weak equivalence, and we have to show that the same edge is a weak equivalence. And of course, we do this using the same argument we used in Proposition~\ref{prop:mfib-hofb-cocart-stab}.
\qed

\subsubsection{}\defn
If $\D\subseteq\C$ is any full subcategory and $*$ is one of $\cf$, $\fb$, or $\cfb$ (the case of interest being $\D=\C_\fb\subseteq\C$ and $*=\fb$), we define the fibration $\fib{HoF_*(\D)}$ to be the \emph{restriction} of $\fib{HoF_*(\C)}$ to $\D$ -- i.e., the total category of $\fib{HoF_*(\D)}$ is the subcategory of $\Ho(\C^\to)_*$ consisting of those objects and morphisms lying over $\D$.

In general, the restriction of any $\wedge$-fibration to a full subcategory is still a $\wedge$-fibration; in particular, $\fib{HoF_\fb(\C_\fb)}$ is a $\wedge$-fibration.

\subsubsection{}\label{thm:mfib-hofb-wedgeq}\thm
The $\wedge$-fibration $\fib{HoF_\fb(\C_\fb)}$ is a $\wedgeq$-fibration.
\pf
By Proposition~\ref{prop:mfib-hofb-has-cocarts}, every morphism in $\C_\fb$ admits a cocartesian lift and by Proposition~\ref{prop:mfib-hofb-cocart-frob}, these satisfy Frobenius reciprocity.

Since $\C_\fb$ has only fibrant objects, every product projection in $\C_\fb$ is a fibration. Hence, by Proposition~\ref{prop:mfib-hofb-cocart-stab}, the cocartesian morphisms are stable along all product projections.
\qed

\section{The homotopy 2-category}\label{sec:classical-2cat}
By the results of Parts~\ref{sec:proof}~and~\ref{sec:modelcats}, we now have a 2-categorical structure on $\C_\fb$, and in particular on $\C_\cfb$, for any model category $\C$. There is another, more familiar (especially in the case $\C=\Top$) 2-categorical structure that one can put on $\C_\cfb$, and in Part~\ref{sec:classical-2cat}, we compare the two.

In the case of the category $\C=\Top$ of topological spaces, the 2-categorical structure on $\C_\cfb$ is roughly as follows: 0-cells are spaces, 1-cells are continuous maps, and 2-cells are homotopy-classes (with ``fixed endpoints'') of homotopies -- the reason one needs to take \emph{homotopy-classes} of homotopies is so that the composition is strictly associative, as is familiar from the definition of the fundamental group.

In the case of a general model category $\C$, the 2-categorical structure $\C_\cfb$ was more or less completely described already in \cite[I~pp.~2.1-2.8]{quillen-ha}. The idea, of course, is that the above description of the 2-categorical structure on $\Top$ uses only notions that are already available in a general model category. One subtlety is that there are now the two notions of left and right homotopy, and in fact in \emph{loc. cit.} this is exploited in a clever way to obtain a very clean notion of ``homotopy of homotopies'' (see \S\ref{subsec:htpies-btw-htpies}).

In an earlier version of this paper, we explicitly defined the 2-categorical structure on $\Top$, and showed that it was the same as ours. Below, we will instead show (Theorem~\ref{thm:model-structure-agrees-with-quillen}) that our 2-categorical structure -- i.e., the 2-cells, and the notion of horizontal and vertical composition -- agree with those from \cite{quillen-ha}. That this structure agrees with the expected one on $\Top$ is then a straightforward exercise involving judicious choices of appropriate path or cylinder objects.

It was also pointed out to me by the anonymous referee that the 2-categorical structure on $\C_\cfb$ exhibits a nice universal property, namely it a \emph{localization} in an obvious 2-categorical sense: any functor from (the 1-category) $\C_\cfb$ to a 2-category $\D$ taking every weak equivalence to an equivalence in $\D$ extends uniquely to the 2-categorical structure on $\C_\cfb$, though we will not prove this.

\subsection{Homotopies between homotopies}\label{subsec:htpies-btw-htpies}
We now introduce the notion of homotopy of homotopies from \cite[\S{}I-2]{quillen-ha}.

There are in fact three such notions: \emph{left homotopy of left homotopies}, \emph{right homotopy of right homotopies}, and \emph{correspondence} between left homotopies and right homotopies, always defined for homotopies between morphisms $A\to{}B$ with $A$ cofibrant and $B$ fibrant.

These turn out all to agree, in the sense that two left (or right) homotopies are left (or right) homotopic if and only they correspond to a common right (or left) homotopy (see Proposition~\ref{prop:right-left-correspondence}).

We will need an additional notion, that of \emph{strong left homotopy of right-homotopies} (see Definition~\ref{defn:right-right-htpy-htpy}) as this is the one that naturally occurs in our 2-category. It is also equivalent to the other notions, but apparently only when dealing with path objects $B^I$ for which $\sigma\colon{}B\to{}B^I$ is a cofibration. In general, it is only equivalent to a variation of the notion of right homotopy, which we call \emph{strong right homotopy} (see Proposition~\ref{prop:right-strong-left-right}).

We will have occasion below to deal with two path objects on the same object $B$ of a model category. We will use the usual notation $B^I$ for the first, and $B\tox{s'}B^{I'}\tox{\br{d_1',d_2'}}B\times{}B$ for the second.

For the rest of \S\ref{subsec:htpies-btw-htpies}, let $\C$ be a model category.

\subsubsection{}\label{defn:right-right-htpy-htpy}\defn
Let $f,g\colon{}A\to{}B$ be morphisms in $\C$ with $A$ cofibrant and $B$ fibrant, and let $k\colon{}A\to{}B^I$ and $k'\colon{}A\to{}B^{I'}$ be right-homotopies from $f$ to $g$. A \emph{right homotopy} from $k$ to $k'$ is a commutative diagram
\begin{equation}\label{eq:weak-right-homotopy-diag}
  \begin{tikzcd}[column sep=30, row sep=30]
    A\ar[r, "K"]\ar[d, "\br{k,k'}"']&B^J\ar[from=d, "t"']\ar[dl, "\br{e_1,e_2}"']\\
    B^I\times_{B\times{}B}B^{I'}\ar[from=r, "\br{s,s'}"]&B
  \end{tikzcd}
\end{equation}
with $\br{e_1,e_2}$ a fibration and $t$ a weak equivalence, where $B^I\times_{B\times{}B}B^{I'}$ is a pullback along the maps $\br{d_1,d_2}\colon{}B^I\to{}B\times{}B$ and $\br{d_1',d_2'}\colon{}B^{I'}\to{}B\times{}B$. Note that, despite the notation (which we borrow from \cite{quillen-ha}), $B^J$ is not a path object for $B$.

If $B^I=B^{I'}$ are the same path object, then a \emph{strong right homotopy} from $k$ to $k'$ is a commutative diagram as follows, with $\br{e_1,e_2}$ a fibration and $t$ a weak equivalence.
\begin{equation}\label{eq:strong-right-homotopy-diag}
  \begin{tikzcd}[column sep=40, row sep=30]
    A\ar[r, "K"]\ar[d, "\br{k,k'}"']&B^J\ar[from=d, "t"']\ar[dl, "\br{e_1,e_2}"']\\
    B^I\times_{B\times{}B}B^{I}\ar[from=r, "\br{\id_{B^I},\id_{B^I}}"]&B^I
  \end{tikzcd}
\end{equation}

Still assuming $B^I=B^{I'}$, a \emph{strong left homotopy} from $k$ to $k'$ is a left homotopy $H\colon{}A\times{}I\to{}B$ from $k$ to $k'$ for which the following diagram commutes.
\[
  \begin{tikzcd}
    A\times{}I\ar[r, "H"]\ar[d, "\sigma"']&B^I\ar[d, "\br{d_1,d_2}"]\\
    A\ar[r, "\br{f,g}"]&B\times{}B
  \end{tikzcd}
\]

We say that $k$ and $k'$ are \emph{right homotopic} if there is a right homotopy between them, and similarly with \emph{strongly right homotopic} and \emph{strongly left homotopic}.

There are dual notions of \emph{left homotopy} between left homotopies, and so on.

\subsubsection{}\label{prop:right-strong-left-right}\prop
With $A,B,k,k'$ as in Definition~\ref{defn:right-right-htpy-htpy}, and assuming $B^I=B^{I'}$, $k$ and $k'$ are strongly left homotopic if and only if they are strongly right homotopic, and in this case, they are right homotopic.

Moreover, if $s\colon{}B\to{}B^I$ is a cofibration and $k$ and $k'$ are right-homotopic, then they are strongly right homotopic.
\pf
For the first equivalence, note that the strong left- and right-homotopies from $k$ to $k'$ are just certain left- and right-homotopies between $k,k'\colon(A,\br{f,g})\to(B^I,\br{d_1,d_2})$ in the slice category $\C/(B\times B)$ (for right-homotopies this is clear, and for left-homotopies it follows Propositions~\ref{prop:fiberwise-homotopy}~and~\ref{prop:slice-total-htpy-compat}). Hence the claim is a special case of the Proposition~\ref{prop:modc-facts}~\ref{item:modc-cyl-choose} (and its dual) since $(A,\br{f,g})$ is cofibrant and $(B^I,\br{d_1,d_2})$ is fibrant.

To get a right-homotopy from a strong right-homotopy, compose $t$ and $\br{\id_{B^I},\id_{B^I}}$ in the diagram (\ref{eq:strong-right-homotopy-diag}) with the weak equivalence $s\colon{}B\to B^I$ to get a diagram as in (\ref{eq:weak-right-homotopy-diag}). Conversely, if $s\colon{}B\to{}B^I$ is a cofibration, then it is a trivial cofibration, and so given a diagram as in (\ref{eq:weak-right-homotopy-diag}), we can find a diagonal filler in the square
\[
  \begin{tikzcd}[column sep=30, row sep=30]
    B\ar[r, "t"]\ar[d, "s"']&B^J\ar[d, "\br{e_1,e_2}"]\\
    B^I\ar[r, "\br{\id_B,\id_B}"']\ar[ru, dotted]&B^I\times_{B\times{}B}B^I
  \end{tikzcd}
\]
(which is a weak equivalence since $s$ and $t$ are) and thus obtain a diagram as in (\ref{eq:strong-right-homotopy-diag}).
\qed

\subsubsection{}\defn
Given morphisms $f,g\colon{}A\to{}B$ in $\C$ with $A$ cofibrant and $B$ fibrant, a left-homotopy $h\colon{}A\times{}I\to{}B$ from $f$ to $g$, and a right-homotopy $k\colon{}A\to{}B^I$ from $f$ to $g$, a \emph{correspondence} between $h$ and $k$ is a morphism $H\colon{}A\times{}I\to{}B^I$ satisfying:
\begin{align*}
  H\partial_1&=k&H\partial_2&=sg\\
  d_1H&=h&d_2H&=g\sigma.
\end{align*}
We say that $h$ and $k$ \emph{correspond} if there exists a correspondence between them.

\subsubsection{}\label{prop:right-left-correspondence}\prop
Given objects $A\in\Ob\C_{\cf}$ and $B\in\Ob\C_\fb$ and morphisms $f,g\colon{}A\to{}B$, the following claims hold, as well as their duals. Together, they say that ``left (or right) homotopic'' is an equivalence relation on left (or right) homotopies, and that correspondence establishes a bijection between left-homotopy classes of left-homotopies and right-homotopy classes of right-homotopies.
\begin{enumerate}[(i)]
\item Right-homotopy is an equivalence relation on right-homotopies from $f$ to $g$.
\item Given a right-homotopy $k\colon{}A\to{}B^I$ from $f$ to $g$ and a cylinder object $A\times{}I$, there is some left-homotopy $h\colon{}A\times{}I\to{}B$ corresponding to $k$.
\item\label{item:rlc-rhtpy-common-cor} Two right-homotopies $k\colon{}A\to{}B^I$ and $k'\colon{}A\to{}B^{I'}$ are right-homotopic if and only if there is some left-homotopy $h\colon{}A\times{}I\to{}B$ to which they both correspond.
\end{enumerate}
\pf
See \cite[II~pp.~2.1-2.5]{quillen-ha}.
\qed

\subsubsection{}\label{prop:right-right-composeroo}\prop
Given $A\in\Ob\C_\cf$, $B\in\Ob\C_\fb$, and a morphism $f\colon{}B^I\to{}B^{I'}$ of path objects making the diagram
\[
  \begin{tikzcd}[column sep=30pt]
    B\ar[r, "s"]\ar[d, "s'"]&B^I\ar[d, "\br{d_1,d_2}"]\ar[ld, "f"]\\
    B^{I'}\ar[r, "\br{d_1',d_2'}"']&B\times{}B
  \end{tikzcd}
\]
commute, any right-homotopy $k\colon{}A\to{}B^I$ is right-homotopic to the right-homotopy $fk\colon{}A\to{}B^{I'}$.
\pf
If $H\colon{}A\times{}I\to{}B^I$ is a correspondence between $k$ and some $h\colon{}A\times{}I\to{}B$, then $fH\colon{}A\times{}I\to{}B^{I'}$ is a correspondence between $fk$ and $h$, and so the claim follows from Proposition~\ref{prop:right-left-correspondence}~\ref{item:rlc-rhtpy-common-cor}.
\qed

\subsubsection{}\defn
Let $f,g,h\colon{}A\to{}B$ be morphisms in $\C$ with $B$ fibrant, and let $k\colon{}A\to{}B^I$ and $k'\colon{}A\to{}B^{I'}$ be right-homotopies from $f$ to $g$ and from $g$ to $h$, respectively. Given any pullback $B^I\times_B{}B^{I'}$ (taken with respect to $d_2\colon{}B^I\to{}B$ and $d_1\colon{}B^{I'}\to{}B$), we have an associated path object
\[
  B\tox{\br{s,s'}}B^I\times_B{}B^{I'}\tox{\br{d_1\pi_1,d_2'\pi_2}}B\times{}B,
\]
(see \cite[\S{}I-1,~Lemma 3]{quillen-ha}) and the induced morphism $\br{k,k'}\colon{}A\to{}B^I\times_BB^{I'}$ is said to be a \emph{composite} of $k$ and $k'$.

\subsection{$\fib{HoF_\fb(\C_\fb)}$-homotopies}\label{subsec:fibc-homotopies}
Henceforth, let $\C$ be a model category.

We now carry out the comparison of the 2-categorical structure on $\C_{\cfb}$ given by Theorems~\ref{thm:mfib-hofb-wedgeq}~and~\ref{thm:2cat} with the one from \cite{quillen-ha}.

\subsubsection{}\label{subsubsec:modc-chtpy}
Strictly speaking, the 2-categorical structure in Theorem~\ref{thm:2cat} is associated to a particular $\wedgeq$-cleavage of $\fib{HoF_\fb(\C_\fb)}$ (see \S\ref{subsec:c-htpy-canon}). Let us fix certain choices in this cleavage -- most importantly, of equality objects -- which will be convenient.

First, for the fiberwise terminal objects, we take the identity morphisms $(A,\id_A)\in\fib{HoF_\fb(\C_\fb)}^A$ for $A\in\Ob\C_\fb$.

By Proposition~\ref{prop:mfib-hofb-wedge}, we have that $\gamma\colon\fib{F_\fb(\C)}\to\fib{HoF_\fb(\C)}$ is a morphism of $\wedge$-fibrations, and so we can (and do) take our cartesian lifts and fiberwise products to be the images under $\gamma$ of the corresponding things in $\fib{F_\fb(\C)}$.

For the equality object $\Eq_B$, we have to choose a cocartesian lift of $\Delta_B\colon{}B\to{}B\times{}B$ for $B\in\Ob\C_\fb$. Choosing a factorization $B\tox{s}B^I\tox{\br{d_1,d_2}}B\times{}B$ of $\Delta_B\colon{}B\to{}B\times{}B$ as a \emph{trivial cofibration} (we will need later on that it is a cofibration) followed by a fibration, we have by Proposition~\ref{prop:mfib-hofb-cocart-crit} that the image of $(s,\Delta_B)\colon(B,B,\id_B)\to(B^I,B\times{}B,\br{d_1,d_2})$ in $\Ho(\C^\to)_\fb$ is cocartesian, and we take this as our chosen equality object.

We fix such a cleavage of $\fib{HoF_\fb(\C_\fb)}$ for the rest of \S\ref{subsec:fibc-homotopies}.

\subsubsection{}\label{prop:char-of-hof-htpies}\prop
Given morphisms $f,g\colon{}A\to{}B$ in $\C_\fb$ with $A$ cofibrant, the $\fib{HoF_\fb(\C_\fb)}$-homotopies from $f$ to $g$ are precisely given by the images under $\gamma$ of the morphisms $(k,\br{f,g})\colon(A,A,\id_A)\to(B^I,B\times{}B,\br{d_1,d_2})$, where $k\colon{}A\to{}B^I$ is a right-homotopy from $f$ to $g$, and $B^I$ is the chosen path object from \S\ref{subsubsec:modc-chtpy}.

Moreover two such homotopies $k,k'\colon{}A\to{}B^I$ map to the same $\fib{HoF_\fb(\C_\fb)}$-homotopy if and only if they are strongly left-homotopic (or equivalently, by Proposition~\ref{prop:right-strong-left-right} and since $B\tox{s}B^I$ is a cofibration, if and only if they are right-homotopic).
\pf
The first statement follows directly from the definitions and Proposition~\ref{prop:hoc-props}~\ref{item:hoc-props-cof-fib-htpy}.

The second statement follows from the definitions and Propositions~\ref{prop:hoc-props}~\ref{item:hoc-props-cof-fib-htpy}~and~\ref{prop:fiberwise-homotopy}.
\qed

\subsubsection{}\label{defn:homotopy-class}\defn
Let $f,g\colon{}A\to{}B$ be morphisms in $\C_\fb$ with $A$ cofibrant.

A \emph{homotopy} from $f$ to $g$ is either a right homotopy or a left homotopy.

We say that two homotopies $h$ and $h'$ are \emph{in the same homotopy class} if $h$ and $h'$ are both left (or both right) homotopies and they are left (or right) homotopic, or if one is a left and one a right homotopy, and they correspond. By Proposition~\ref{prop:right-left-correspondence}, this is an equivalence relation on the (``large'') set of homotopies from $f$ to $g$.

Given a $\fib{HoF_\fb(\C_\fb)}$-homotopy $\alpha\colon\top_A\to\Eq_B$ from $f$ to $g$, we say that a homotopy $h$ from $f$ to $g$ \emph{represents} $\alpha$ if $h$ is in the same homotopy class as a right-homotopy $k$ from $f$ to $g$ whose image under $\gamma$ is $\alpha$.

\subsubsection{}\label{prop:htpies-agree-with-quillen}\prop
With $f,g$ as in Definition~\ref{defn:homotopy-class}, each homotopy from $f$ to $g$ represents a unique $\fib{HoF_\fb(\C_\fb)}$-homotopy, and two homotopies represent the same $\fib{HoF_\fb(\C_\fb)}$-homotopy if and only if they are in the same homotopy class.
\pf
For right-homotopies, this is precisely the content of Proposition~\ref{prop:char-of-hof-htpies}, and for left-homotopies, it follows from Proposition~\ref{prop:right-left-correspondence}.
\qed

\subsubsection{}\label{prop:vert-comp-agrees-with-quillen}\prop
Given morphisms $f,g,h\colon{}A\to{}B$ in $\C_\fb$ with $A$ cofibrant, and $\fib{HoF_\fb(\C_\fb)}$-homotopies $\alpha\colon{}f\to{}g$ and $\beta\colon{}g\to{}h$ represented by right-homotopies $k$ and $k'$ respectively, the composite $\beta\hvcmp\alpha$ is represented by a composite of $k$ and $k'$.
\pf
As the definition of $\beta\hvcmp\alpha$ involves the object $\Eq_B^{12}\wedge\Eq_B^{23}$, let us first identify the latter in the present context.

Let us temporarily write $\abs{P}$ for the domain $X$ of an object $P=(X,A,x)\in\Ob(\C^\to)=\Ob\Ho(\C^\to)$. Recalling that $\Eq_B=(B^I,B\times{}B,\br{d_1,d_2})$, and that both the pullbacks and the fiberwise products in $\fib{F_\fb(\C)}$ (and hence in $\fib{HoF_\fb(\C)}$) are given by pullback squares in $\C$, we have the diagram
\[
  \begin{tikzcd}[column sep={40pt,between origins}, row sep={30pt,between origins}]
    &&\abs{\Eq_B^{12}\wedge\Eq_B^{23}}\ar[ld]\ar[rd]&&\\
    &\abs{\Eq_B^{12}}\ar[ld]\ar[rd]&&\abs{\Eq_B^{23}}\ar[ld]\ar[rd]&\\
    B^I\ar[rd, "\br{d_1,d_2}"']&&
    B^3\ar[ld, "\br{\pi_1,\pi_2}"' near start]\ar[rd, "\br{\pi_2,\pi_3}" near start]&&
    B^I\ar[ld, "\br{d_1,d_2}"]\\
    &B\times{}B\ar[rd, "\pi_2"']&&B\times{}B\ar[ld, "\pi_1"]&\\
    &&B&&
  \end{tikzcd}
\]
of pullback squares in $\C$, which shows that $\Eq_B^{12}\wedge\Eq_B^{23}$ is given by $(B^I\times_BB^I,B^3,\xi)$ for some pullback $B^I\times_BB^I$ along $d_2\colon{}B^I\to{}B$ and $d_1\colon{}B^I\to{}B$, where $\xi=\br{d_1\pi_1,d_2\pi_1,d_2\pi_2}=\br{d_1\pi_1,d_1\pi_2,d_2\pi_2}$.

Next, we have a morphism
\[
  (\id,\br{\pi_1,\pi_3})\colon(B^I\times_BB^I,B^3,\xi)\to
  (B^I\times_BB^I,B\times{}B,\br{d_1\pi_1,d_2\pi_2})
\]
in $(\C^\to)_\fb$, the image of which in $\Ho(\C^\to)_\fb$ is cocartesian by Proposition~\ref{prop:mfib-hofb-cocart-crit}. Let us write $\widetilde\Eq_B$ for the codomain of this morphism.

We then have the following commutative diagram in $\Ho(\C^\to)_\fb$:
\[
  \begin{tikzcd}
    &\Eq_B^{12}\wedge\Eq_B^{23}\ar[d, "\tr_B"]\rac[r, "{\gamma(\id,\br{\pi_1,\pi_3})}"]&[20pt]
    \widetilde\Eq_B\ar[d, "p"]\\
    \top_B\rac[ru, "\brr{\rho_B^{12},\rho_B^{23}}"]\ar[r, "\rho_B^{13}", near end]
    \rac[rr, bend right=10pt, "\rho_B"']&
    \Eq_B^{13}\ar[r, "\ct"]&\Eq_B\\
    B\ar[r, "\Delta_B^3"]&B^3\ar[r, "\br{\pi_1,\pi_3}"]&B\times{}B
  \end{tikzcd}
\]
where $p$ is the unique morphism making the diagram commute. Since the morphisms $\top_B\to\Eq_B$ and $\top_B\to\widetilde\Eq_B$ in the diagram are cocartesian, it follows that $p$ is an isomorphism.

Next, we can find a lift $\hat{p}$ of $p\I$ to $\fib{F(\C)}^{B\times{}B}$ as a filler in the following diagram (in $\C$) since we are assuming that $s$ is a cofibration.
\[
  \begin{tikzcd}
    B\ar[d, "s"']\ar[r, "\br{s,s}"]&
    B^I\times_BB^I\ar[d ,"\br{d_1\pi_1,d_2\pi_2}"]\\
    B^I\ar[r, "\br{d_1,d_2}"']\ar[ru, dashed, "\hat{p}"]&B\times{}B
  \end{tikzcd}
\]

Now consider the commutative diagram
\[
  \begin{tikzcd}
    &\Eq_B^{12}\wedge\Eq_B^{23}\ar[d, "\tr_B"]\rac[r, "{\gamma(\id,\br{\pi_1,\pi_3})}"]&[20pt]
    \widetilde\Eq_B\ar[from=d, "p\I"']\\
    \top_A\ar[ru, "\brr{\cind{\alpha},\cind{\beta}}"]\ar[r, "\cind{\beta\hvcmp\alpha}"']&
    \Eq_B^{13}\ar[r, "\ct"]&\Eq_B\\
    A\ar[r, "\br{f,g,h}"]&B^3\ar[r, "\br{\pi_1,\pi_3}"]&B\times{}B.
  \end{tikzcd}
\]

The outside of this diagram lifts to a diagram in $\C^\to$
\[
  \begin{tikzcd}
    &[-10pt]B^I\times_BB^I\ar[r, "\id"]&[10pt]B^I\times_BB^I\\
    A\ar[ru, "\br{k,k'}"]\ar[rr, "k''"']&&B^I\ar[u, "\hat{p}"']\\[-10pt]
    A\ar[r, "\br{f,g,h}"]&B^3\ar[r, "\br{\pi_1,\pi_3}"]&B\times{}B.
  \end{tikzcd}
\]
commutative up to homotopy, where $k''$ is a lift of $\beta\hvcmp\alpha$, which exists by Proposition~\ref{prop:hoc-props}~\ref{item:hoc-props-cof-fib-htpy} (note that for the sake of readability we are conflating objects in $\C^\to$ with their domains).

Since the diagram commutes up to homotopy the morphisms $\br{k,k'}$ and $\hat{p}k''$ represent the same $\fib{HoF_\fb(\C_\fb)}$-homotopy. But $k''$ and hence, by Proposition~\ref{prop:right-right-composeroo}, $\hat{p}k''$ represents $\beta\hvcmp\alpha$, and $\br{k,k'}$ is by definition a composite of $k$ and $k'$. This proves the claim.
\qed

\subsubsection{}\label{prop:hor-comp-agrees-with-quillen}\prop
Given objects and morphisms
\[
  \begin{tikzcd}
    A\ar[r, "f"]&B\ar[r, "f'", shift left]\ar[r, "g'"', shift right]&C\ar[r, "g"]&D
  \end{tikzcd}
\]
in $\C_\cfb$ and a $\fib{HoF_\fb(\C_\fb)}$-homotopy $\alpha\colon{}f'\to{}g'$ represented by a right-homotopy $k\colon{}B\to{}C^I$ and by a left-homotopy $h\colon{}B\times{}I\to{}C$, the $\fib{HoF_\fb(\C_\fb)}$-homotopies $\alpha\hhcmp\hid_f$ and $\hid_g\hhcmp\alpha$ are represented, respectively, by $kf\colon{}A\to{}C^I$ and $gh\colon{}B\times{}I\to{}D$.
\pf
The $\fib{HoF_\fb(\C_\fb)}$-homotopy $\alpha\hhcmp\hid_{f}$ is given by the composite
\[
  \top_A\tox{\exx_f}\top_B\tox{\alpha}\Eq_C,
\]
and it is immediate that this is represented by $kf\colon{}A\to{}C^I$.

Next, $\hid_g\hhcmp\alpha$ is given by the composite
\[
  \top_B\tox{\alpha}\Eq_C\tox{\widecheck{\hid_g}}\Eq_D,
\]
We can find a lift $\br{r,g\times{}g}$ of $\widecheck{\hid_g}$ to $\C^\to$ as a filler in the following diagram (in $\C$) since we are assuming that $s$ is a cofibration.
\[
  \begin{tikzcd}[column sep=30pt]
    C\ar[d, "s"']\ar[r, "sg"]&
    D^I\ar[d ,"\br{d_1,d_2}"]\\
    C^I\ar[r, "\br{gd_1,gd_2}"']\ar[ru, dashed, "r"]&D\times{}D
  \end{tikzcd}
\]
Hence, $\hid_g\hhcmp\alpha$ is represented by $rk$. Taking some correspondence $H\colon{}B\times{}I\to{}C^I$ between $h$ and $k$, the composite $B\times{}I\tox{H}C^I\tox{r}D^I$ then gives a correspondence between $gh$ and $rk$.
\qed

\subsubsection{}\defn
We define the \emph{Quillen 2-categorical structure on $\C_\cfb$} as follows.

Given morphisms $f,g\colon{}A\to B$ in $\C_\cfb$, the 2-cells $f\to g$ are the homotopy classes of (left or right) homotopies.

Given morphisms $f,g,h\colon{}A\to B$ and 2-cells $f\tox{\alpha}g\tox{\beta}h$ represented by right-homotopies $k\colon{}A\to B^I$ and $k'\colon{}A\to B^{I'}$, the vertical composite of $\alpha$ and $\beta$ is defined to be the homotopy represented by any composite of $k$ and $k'$.

Given objects and morphisms
\[
  \begin{tikzcd}
    A\ar[r, "f"]&B\ar[r, "f'", shift left]\ar[r, "g'"', shift right]&C\ar[r, "g"]&D
  \end{tikzcd}
\]
in $\C_\cfb$ and a 2-cell $\alpha\colon{}f'\to{}g'$ represented by a right-homotopy $k\colon{}B\to{}C^I$ and by a left-homotopy $h\colon{}B\times{}I\to{}C$, the horizontal composites of $\alpha$ with $f$ and $g$ are represented, respectively, by $kf\colon{}A\to{}C^I$ and $gh\colon{}B\times{}I\to{}D$ (this uniquely determines the horizontal composites of arbitrary pairs of 2-cells).

That these operations are well-defined is proven in \cite[\S{}I-2]{quillen-ha}, and that they define a 2-category mostly follows from what is proven there as well.

But in any case, both the well-definedness and the fact that they satisfy the 2-category axioms follows from Propositions~\ref{prop:htpies-agree-with-quillen},~\ref{prop:vert-comp-agrees-with-quillen},~and~\ref{prop:hor-comp-agrees-with-quillen}.

\subsubsection{}\thm\label{thm:model-structure-agrees-with-quillen}
The 2-categorical structure on $\C_{\cfb}$ given by Theorems~\ref{thm:mfib-hofb-wedgeq}~and~\ref{thm:2cat} is isomorphic to the Quillen 2-categorical structure.
\pf
Propositions~\ref{prop:htpies-agree-with-quillen},~\ref{prop:vert-comp-agrees-with-quillen},~and~\ref{prop:hor-comp-agrees-with-quillen} provide the desired isomorphism.
\qed

\bibliographystyle{alpha}
\bibliography{htpy_in_fibs}
\end{document}